\documentclass[11pt,a4paper]{article}
\usepackage [english]{babel}
\usepackage [utf8]{inputenc}
\usepackage{amsmath,amssymb,amsthm,dsfont}
\usepackage{algorithm,placeins,algpseudocode}
\usepackage{xcolor,graphicx}
\usepackage[a4paper,
            left=2.5cm,
            right=2.5cm,
            top=3cm,
            bottom=3cm]{geometry}
\usepackage{geometry}
\usepackage{subcaption}
\usepackage{eqparbox,array}
\usepackage{float}
\usepackage{stmaryrd} 
\usepackage{relsize}
\usepackage{enumitem}   
\usepackage{microtype}  
\definecolor{Mycolor2}{HTML}{4B0082}
\usepackage[colorlinks=true,linkcolor= Mycolor2,citecolor= Mycolor2]{hyperref}

\usepackage{comment}

\newcommand{\RR}{\mathbb{R}}
\newcommand{\NN}{\mathbb{N}}

\newcommand{\RRP}{\mathbb{R}^+_*}
\newcommand{\MC}{\mathcal{M}}

\newcommand{\LC}{\mathcal{L}}

\newcommand{\FC}{\mathcal{F}}
\newcommand{\PC}{\mathcal{P}}

\newcommand{\SCE}{\emph{SCE}}

\newcommand{\A}{(A)}
\newcommand{\Diag}{\mathrm{Diag}}

\algnewcommand{\LongComment}[1]{\hfill// \begin{minipage}[t]{\eqboxwidth{COMMENT\thealgorithm}}#1\strut\end{minipage}}

\newcommand{\E}[1]{\mathbb{E}\left[#1\right]}
\newcommand{\Prob}[1]{\mathbb{P}\left(#1\right)}
\newcommand{\Proc}[1]{\left(#1\right)_{t\geq 0}}

\newcommand{\Seq}[1]{\left(#1\right)_{n\in \mathbb{N}}}

\newcommand{\brac}[1]{\left\langle#1\right\rangle}

\newcommand{\Wnorm}[2]{\left\llbracket #1 \right\rrbracket_{#2}}

\usepackage{todonotes}

\newcommand{\dd}{\mathop{}\!\mathrm{d}}

\newcommand\restr[2]{{
  \left.\kern-\nulldelimiterspace 
  #1 
  \littletaller 
  \right|_{#2} 
  }}

\newcommand{\littletaller}{\mathchoice{\vphantom{\big|}}{}{}{}}

\newtheorem{theorem}{Theorem}[section]
\newtheorem{lemma}[theorem]{Lemma}
\newtheorem{remark}[theorem]{Remark}
\newtheorem{example}[theorem]{Example}
\newtheorem{proposition}[theorem]{Proposition}
\newtheorem{definition}[theorem]{Definition}
\newtheorem{corollary}[theorem]{Corollary}
\theoremstyle{definition}
\newtheorem{assumptionA}{Assumption}

\theoremstyle{plain}
\renewcommand{\A}{(\ref{ass:observable})}

\usepackage{authblk}

\title{Convergence rate of observables in the Marcus-Lushnikov process toward Smoluchowski with singular homogeneous kernels}
\author[1]{Sébastian Baudelet}
\author[2]{Rémi Catellier}
\affil[1]{INRIA Saclay}
\affil[2]{Université Côte d'Azur, CNRS, LJAD, France}
\date{}

\date{}

\begin{document}

\maketitle

\begin{abstract}
We study the approximation of the Smoluchowski coagulation equation by the
Marcus--Lushnikov process with $N$ initial particles for the singular family of kernels $K_\alpha(x,y)=x^{-\alpha}+y^{-\alpha}, \alpha>0$. For regular mean-field particle systems, pathwise fluctuations are naturally
expected on the scale $N^{-1/2}$, whereas sufficiently smooth observables may
exhibit a weak error of order $N^{-1}$ after averaging. The singularity of
$K_\alpha$ at the origin makes it unclear whether this second-order mechanism
survives.

We prove that it does: for a class of sufficiently regular nonlinear
functionals of the empirical measure, we establish a weak error of order
$N^{-1}$ on finite time intervals. The main mechanism is a regularization
produced by the singular dynamics itself. Small clusters coagulate at
increasingly large rates, which generates exponential moments in the singular
variable $x^{-\alpha}$ and yields a rapid escape from the singularity at zero.
We exploit this phenomenon through a Wasserstein-type distance adapted to
$x^{-\alpha}$, together with weighted total variation estimates used to
construct and control the derivative of the Smoluchowski flow.
\end{abstract}



\section{Introduction}

The \emph{Smoluchowski coagulation equation} (SCE), introduced by
Smoluchowski in the early twentieth century
\cite{smoluchowski1927drei}, describes the evolution of a distribution of
cluster masses undergoing binary coagulation. In weak form, it reads
\begin{equation}\label{eq:SCE}
    \begin{aligned}
    \int_{\RRP} f(x)\mu_t(\dd x)
    ={}& \int_{\RRP} f(x)\mu_0(\dd x) \\
    &+\frac12\int_0^t
      \int_{\RRP}\int_{\RRP}
      K(x,y)
      \left[f(x+y)-f(x)-f(y)\right]
      \mu_s(\dd x)\mu_s(\dd y)\,\dd s ,
    \end{aligned}
\end{equation}
for bounded measurable test functions $f$. Here, $\mu_t$ represents the distribution
of cluster masses at time $t$, while the symmetric non-negative function
$K$ is the \emph{coagulation kernel}: two clusters of masses $x$ and $y$
merge at a rate prescribed by $K(x,y)$.

A stochastic counterpart is obtained by considering a finite system of
clusters undergoing the same binary mergers. Starting from $N$ particles,
let
\[
    Z_t^N=\sum_{i=1}^{N_t}\delta_{X_i(t)}
\]
denote the associated finite stochastic coalescence process, in which each
unordered pair of clusters of masses $x$ and $y$ merges at rate $K(x,y)$
and is replaced by a single cluster of mass $x+y$. Here $N_t$ is the number
of clusters at time $t$ and decreases by one at each coagulation event.
The corresponding mean-field rescaling is
\[
    M_t^N:=\frac1N Z_{t/N}^N.
\]
Under this convention, $M^N$ is the usual \emph{Marcus--Lushnikov process}. The factor $1/N$ is the
standard empirical-measure normalization, while the time scaling compensates
for the $O(N^2)$ possible interacting pairs and produces an $O(N)$ number
of coagulation events on macroscopic time intervals.

The convergence of the Marcus--Lushnikov process towards solutions of the
Smoluchowski equation has been studied extensively for broad classes of
coagulation kernels. A foundational result is due to Norris
\cite{norrisSmoluchowskisCoagulationEquation1999}; see also the review of
Aldous \cite{aldousDeterministicStochasticModels1999} for the deterministic
and stochastic theories of coagulation.

\paragraph{Context and objective.}
The starting point of this work is a biological coagulation model introduced
in~\cite{garcia2025multi}. In that setting, the formation of a fungal network
may depend on the simultaneous presence of at least two clusters exceeding a
prescribed mass threshold. This naturally leads to nonlinear observables of
the cluster distribution. The coagulation kernels arising from a scaling limit
of the underlying spatial model take the singular form
\[
    K_\alpha(x,y)=x^{-\alpha}+y^{-\alpha},
    \qquad \alpha>0.
\]

The original biological observable is a hard threshold. If $Z_t^N$ denotes
the unnormalised particle system, the event that at least two clusters have
mass larger than $x_\star$ is
\[
    \left\{Z_t^N([x_\star,\infty))\ge2\right\}.
\]
In terms of the empirical measure $M_t^N=N^{-1}Z_t^N$, this corresponds to
the discontinuous and $N$-dependent functional
\[
    \mathbf 1_{\{
        M_t^N([x_\star,\infty))\ge 2/N
    \}}.
\]
Such hard-threshold observables are not covered by the present analysis.
Their treatment involves a separate difficulty associated with the
discontinuity at the observation threshold. Here we isolate instead the
difficulty created by the singularity of the coagulation dynamics at zero,
and study sufficiently regular nonlinear observables of the empirical
measure.

Two quantitative scales naturally arise in mean-field particle
approximations. Fluctuations of the empirical process, and correspondingly
strong or pathwise approximation errors, typically occur on the scale
$N^{-1/2}$. In the coagulation setting, quantitative estimates and central
limit theorems exhibiting this scale have been obtained for several classes
of kernels; see for instance
\cite{cepedaSmoluchowskisEquationRate2011a,
      kolokoltsov2010central,
      amorim2026scaling}.
For sufficiently regular weak observables, however, averaging cancels the
leading centered fluctuation and an error of order $N^{-1}$ is expected.
This phenomenon is classical in regular mean-field particle systems; see,
for example,
\cite{chassagneux2022weak,bencheikh2022weak}.

The kernel $K_\alpha$ lies outside the standard quantitative frameworks used
for these results. Indeed, the interaction rate becomes arbitrarily large as
soon as either coagulating particle approaches zero. This raises the main
question of the present work:

\medskip
\begin{center}
\emph{Does the classical $N^{-1}$ weak-error mechanism survive the
singularity of $K_\alpha$, and if so, what mechanism makes this possible?}
\end{center}

\paragraph{Main result}

Fast coagulation of small particles make them exiting $0$ at a speed that creates negative $\alpha$-exponential moments see Proposition~\ref{prop:exponential-moments} for the precise statement. This regularization mechanism is at the heart of the theorem below. 
\begin{theorem}[Informal statement]
    Let $\Proc{\mu_t}$ be the solution to the SCE with initial value $\mu$. Let $N\ge 1$ , and let $X_1, \cdots, X_N$ be i.i.d.\ random variables with law $\mu$. Assume that the Marcus-Lushnikov process is initialized by
    \begin{align*}
        M^N_0 = \dfrac{1}{N} \sum\limits_{i = 1}^N \delta_{X_i}.
    \end{align*}
    Then, for all $t \in [0,T]$, we have
    \[
        \left| \mathbb{E}\left[F\left(M^N_t\right)\right] - F(\mu_t) \right| \leq \frac{C(T, F, \mu)}{N},
    \]
    where $F$ is a suitable and sufficiently regular observable, that is, a functional defined on a space of measures and compatible with the singularity of the kernel.
\end{theorem}
At a structural level, the proof can be viewed as a second-order Taylor expansion of the stochastic coagulation generator around the deterministic Smoluchowski flow. The elementary jumps of the particle system perturb the empirical measure by signed measures of the form $\delta_{x+y}-\delta_x-\delta_y$, and controlling the resulting Taylor remainder requires differentiability of the observable in a topology adapted to $K_\alpha$. The natural topology is not the usual weak or total variation topology, but a Wasserstein-type geometry built from the variable $x^{-\alpha}$. A central part of the paper is to show that this geometry is stable under both the Smoluchowski flow and its derivative with respect to the initial measure. Crucially, this geometry is also compatible with the regularization mechanism generated by the singular dynamics: the exponential moments created in the variable $x^{-\alpha}$ provide exactly the control needed to close the corresponding Wasserstein estimates.

Our analysis combines ideas from several sources: the Wasserstein-type distance construction from~\cite{cepedaSmoluchowskisEquationRate2011a}, the existence theory in~\cite{norrisSmoluchowskisCoagulationEquation1999}, the regularity analysis from~\cite{kolokoltsov2010central}, and some methods of~\cite{martiniKolmogorovEquationsSpaces2023}, to treat analysis on spaces of measures. The treatment of nonlinear observables also draws from recent advances in controlling errors between functionals of empirical measures, as in~\cite{delarue2025uniform} to show propagation of chaos. 

\paragraph{Related work} The SCE has been the subject of extensive study from both deterministic and probabilistic perspectives. A seminal review of the field is provided by Aldous~\cite{aldousDeterministicStochasticModels1999}, while a recent survey of the fragmentation side of the theory can be found in~\cite{deaconuProbabilisticRepresentationsFragmentation2023}. Existence and uniqueness for the Smoluchowski Coagulation Equation~\eqref{eq:SCE} have been studied extensively, including under general conditions in~\cite{laurencotClassContinuousCoagulationFragmentation2000,lambExistenceUniquenessResults2004,barikMassconservingSolutionsSmoluchowski2020}. A well-studied phenomenon is the appearance of clusters of infinite size for certain kernels known as gelling~\cite{escobedoGelationCoagulationFragmentation2002,fournier2025gelation}. However, the singular kernel considered in this work is non-gelling and falls outside the scope of most classical treatments, with some exceptions. Notably, in~\cite{fournierExistenceSelfSimilarSolutions2005} some solutions to the $\SCE$ are constructed but the stochastic coalescence is not considered. 

\paragraph{Outline} The rest of the paper is organized as follows. In Section~\ref{section-main-results}, we present the main results along with the key definitions and an overview of the proof strategies. In Section~\ref{section:well-posedness-SCE}, we recall the deterministic SCE estimates, prove the exponential moment generation, and construct the derivative flow. Section~\ref{section:Wass-regularity} contains the core analytical arguments, focusing on Wasserstein regularity estimates and associated technical results. The complete proofs of our main results are then provided in Section~\ref{Section:main-results}.

\section{Definitions, assumptions, and main results}\label{section-main-results}
This section fixes the objects used throughout the paper before stating the assumptions and main results.

We first introduce the $x^{-\alpha}$-weighted measure spaces and the corresponding topologies tailored to the singularity of the kernel. We then define flat derivatives for functions on spaces of measures, which is the calculus used in the Taylor expansions behind the final estimate. Next we formulate the SCE, record the deterministic properties needed later, and construct the Marcus-Lushnikov process. With these objects in place, we state the assumptions on the observable $F$, the main results, and the proof strategy.

\subsection{Notation}

Throughout this work, the parameter $\alpha > 0$ appearing in the kernel $K_\alpha$ is fixed.

The letter $c>0$ denotes a constant depending only on the fixed parameter $\alpha$ (and on $p$ in estimates involving $\psi_p$-functional weights defined later) whose value may change
from line to line. Constants whose dependence is relevant are denoted by
$C(\cdot)$.

Throughout this work, the term \emph{mass} may refer either to the value of a measure on Borel sets (in the measure-theoretic sense), or to the individual atom values in the context of stochastic coalescence. The meaning will be clear from context or explicitly specified when necessary.

\subsection{Measure spaces and topology}

We denote by $\MC^+(\RRP)$ the set of finite positive measures on $\RRP$ and by $\MC(\RRP)$ the set of finite signed measures. For a measure $\mu \in \MC(\RRP)$, we denote by $|\mu| \in \MC^+(\RRP)$ its total variation measure, and by $\|\mu\|$ its total variation norm whose existence comes from the Hahn-Jordan decomposition theorem \cite{rudin1987real}. Given our kernel, we work with measures that have a finite negative $\alpha$-moment, namely measures $\mu \in  \MC^+(\RRP)$ such that:
\begin{align*}
    \int_{\RRP} x^{-\alpha} \mu(\dd x) <  \infty,
\end{align*}
or signed measures $\mu \in  \MC(\RRP)$ such that:
\begin{align*}
    \int_{\RRP} x^{-\alpha} |\mu|(\dd x) <  \infty.
\end{align*}
In full generality let $\psi$ be a measurable, positive function. We define:
\[
\MC_{\psi}(\RRP) = \left\lbrace \mu \in \MC(\RRP) \;\middle|\; \int_{\RRP} \psi(x)\,|\mu|(\dd x) < \infty \right\rbrace,
\]
and $\MC^+_{\psi}(\RRP)$ its restriction to positive measures. This set is the set of finite measures with finite $\psi$-moment. We primarily work with the spaces $\MC_{\psi_p}(\RRP)$ and $\MC^+_{\psi_p}(\RRP)$ where the function $\psi_p$ is defined for all $p \geq 0$ by:
\[
\psi_p(x) = x^{-p\alpha}.
\]
These are the spaces of finite signed and positive measures with finite negative $\alpha p$-moment. Note that if $p \leq q$, then for all $x \in \RRP$:
\[
\psi_p(x) \leq 1 + \psi_q(x),
\]
and therefore $\MC_{\psi_q}(\RRP) \subset \MC_{\psi_p}(\RRP)$, and $\MC^+_{\psi_q}(\RRP) \subset \MC^+_{\psi_p}(\RRP)$. 

\subsubsection*{Total variation}
A natural norm on $\MC_{\psi}(\RRP)$ is the $\psi$-weighted total variation norm:
\[
\|\mu\|_{\psi} = \brac{1 + \psi, |\mu|}.
\]
Note that if $\mu$ is a positive measure then $\mu = |\mu|$ and $\|\mu\|_{\psi} = \brac{1 + \psi,\mu}$. It is well known that the space $\left(\MC_{\psi}(\RRP),\|\cdot\|_{\psi} \right)$ is Banach. This topology is primarily used in Section~\ref{section:well-posedness-SCE} for questions of well-posedness.

\subsubsection*{Wasserstein Distance}
We require a distance that quantifies the convergence rate of sequences of random empirical measures toward their law. The previous weighted-total variation distance, although adapted to the kernel singularity, does not meet this requirement. First we introduce our Wasserstein-type distance and highlight the relation with the classical Wasserstein-$1$.

Let $(E,d)$ be a separable metric space. We denote by $\PC_1(E,d)$ the set of probability measures $\mu$ on $E$ with finite first moment, namely such that for some, hence every, $x_0\in E$,
\[
\int_E d(x,x_0)\,\mu(\dd x)<\infty .
\]
For $\mu,\nu\in\PC_1(E,d)$, the Wasserstein-$1$ distance is defined by
\[ 
W_1^d(\mu,\nu) = \inf\limits_{\pi \in \Pi(\mu,\nu) } \int_{ E \times E } d(x,y) \pi(\dd x, \dd y),
\]
where $\Pi(\mu,\nu)$ is the set of couplings with marginals $\mu$ and $\nu$. One key result is the Kantorovich-Rubinstein duality theorem, which states that
\[
W_1^d(\mu,\nu) = \sup\limits_{\operatorname{Lip}_d(f) \le 1} \brac{f, \mu - \nu},
\]
where the supremum is taken over Lipschitz functions $f$ with Lipschitz constant below $1$ with respect to $d$; for more see \cite{villani2008optimal,santambrogio2015optimal}. When $E \subset \RR$ and $d(x,y)=|x-y|$ we simply write $W_1$ in this case the following long known identity puts the cdfs of $\mu$ and $\nu$ in relation with $W_1$ \cite{vallender1974calculation}:
\[
W_1(\mu,\nu) = \int\left|F_{\mu}(x) - F_{\nu}(x) \right| \dd x.
\]
Our motivation to introduce this distance is that it gives a convergence rate of empirical measures toward their law. The elementary lemma below provides a rate that will suffice for our application.
\begin{lemma}\label{lem:Wasserstein-square-bound}
    Let $p \ge 2$, $q > p$, and take $\mu \in \PC_q(\RRP, |\cdot|)$. Let $\mu^N$ be an empirical measure formed by $N$ i.i.d. atoms with law $\mu$, then:
    \[
    \E{W_1\left(\mu^N, \mu \right)^p}^{1/p} \le \dfrac{C_{p,q}}{\sqrt{N}} \left(\int x^q  \mu(\dd x)\right)^{\tfrac1q}
    \]
\end{lemma}
\begin{proof}
    We rely on the cdf identity denoting $F_N(x) = \mu^N\left((0,x] \right)$ and $F$ the cdf of $\mu$. By Minkowski's inequality:
    \begin{align*}
    \E{W_1\left(\mu^N, \mu \right)^p}^{1/p} &\le \int \E{\left|F_N(y) - F(y) \right|^p}^{\tfrac1p}\dd y.
    \end{align*}
    Take $X$ a random variable with law $\mu$, Rosenthal's inequality \cite{rosenthal1970subspaces} gives 
    \begin{align*}
        \E{\left|F_N(y) - F(y) \right|^p} 
        &\le C_p \left(N^{-p+1}\E{\left|\mathds{1}_{X \le y} - F(y)\right|^p} + N^{-\tfrac{p}{2}}\E{\left|\mathds{1}_{X \le y} - F(y)\right|^2}^{\tfrac{p}{2}} \right).
    \end{align*}
    In the following $C_p$ may change value from line to line. The first term satisfies
    \[
    \E{\left|\mathds{1}_{X \le y} - F(y)\right|^p} \le F(y)(1-F(y)) .
    \]
    The second term satisfies:
    \[
    \E{\left|\mathds{1}_{X \le y} - F(y)\right|^2} = F(y)(1-F(y))
    \]
    The leading order term in $N$ for $p \ge 2$ is $N^{-\tfrac{p}{2}}$. Using the bound $F(y)(1-F(y) \le 1/2$ we get:
    \[
    \E{\left|F_N(y) - F(y) \right|^p} 
    \le C_p N^{-\tfrac{p}{2}} F(y)(1-F(y)) .
    \]
    
    By Markov's inequality:
    \[
    1 - F(y) = \Prob{X > y} \le \dfrac{\E{X^q}}{y^q}.
    \]
    Using also the bound $F(y)(1-F(y)) \le \tfrac14$ we get:
    \[
    \E{W_1\left(\mu^N, \mu \right)^p}^{1/p} \le \dfrac{C_p}{\sqrt{N}} \int \min\left(\dfrac{1}{4^{\tfrac1p}} ,\dfrac{\E{X^q}^{\frac{1}{p}}}{y^{\tfrac{q}{p}}} \right)\dd y \le \dfrac{C_{p,q}}{\sqrt{N}} \E{X^q}^{1/q}.
    \]
    
\end{proof}

We now define a Wasserstein-type distance adapted to our kernel. The estimates that appear in the following sections naturally generate finite sums such as $\psi_1+\psi_2$ and $\psi_1+\psi_p$. We therefore call a weight $\psi$ \emph{admissible} if
\[
    \psi=\sum_{i\in I}\beta_i\psi_i,
\]
where $I\subset(0,\infty)$ is finite and $\beta_i>0$. Any admissible weight is strictly decreasing, and therefore
\[
d_{\psi}(x,y) = | \psi(x) - \psi(y) |
\]
defines a distance on $\RRP$. Note that, for $p<q$, the moment control $\psi_p\leq 1+\psi_q$ does not imply that the distances $d_{\psi_p+\psi_q}$ and $d_{\psi_q}$ are equivalent. This is why we keep track of admissible weights.

We can now introduce the corresponding Wasserstein distance on $\PC_\psi(\RRP)$, the set of probability measures on $\RRP$ with finite $\psi$-moment:
\[
W_\psi(\mu,\nu) = \inf\limits_{\pi \in \Pi(\mu,\nu) } \int_{ \RRP \times \RRP } d_{\psi}(x,y) \pi(\dd x, \dd y).
\]
It is finite since $d_\psi(x,y) \le \psi(x) + \psi(y)$. It is connected to the classical $W_1$ with the following identity:
\[
W_\psi(\mu,\nu) = W_1 (\psi_\# \mu, \psi_\# \nu ),
\]
where $\psi_\#\mu$ is the pushforward measure defined for all $f$ by:
\[
\brac{f,\psi_\#\mu} = \brac{f \circ \psi, \mu} .
\]
This identity allows us to transfer Lemma \ref{lem:Wasserstein-square-bound} to our distance directly:

\begin{corollary}\label{cor:wass-square-bound}
    Let $\psi$ be an admissible weight, $p \ge 2$, $q > p$, and take $\mu \in \PC_{\psi^q}(\RRP, |\cdot|)$. Let $\mu^N$ be an empirical measure formed by $N$ i.i.d. atoms with law $\mu$, then:
    \[
    \mathbb{E}\left[W_\psi \left(\mu^N, \mu \right)^p \right]^{\tfrac1p} \le \dfrac{C_{p,q}} {\sqrt{N}} \left(\int \psi(x)^q \mu(\dd x)\right)^{1/q}.
    \]
\end{corollary}

We would like to extend this distance to $\MC_{\psi}(\RRP)$. For that purpose we rely on a duality formulation. First, for $f : \RRP \to \RR$, set
\[
\mathrm{Lip}_{d_\psi}(f) := \sup_{\substack{x,y\in\RRP\\x\neq y}}
\frac{|f(x)-f(y)|}{d_{\psi}(x,y)}
\]
and
\[
\Wnorm{f}{\psi}
:=
\max\left\{
\sup_{x\in\RRP}\frac{|f(x)|}{1+\psi(x)}
,
\mathrm{Lip}_{d_\psi}(f)
\right\}.
\]
We introduce the normed class
\[
\FC_{\psi}
=
\left\{
f:\RRP\to\RR
\;\middle|\;
\Wnorm{f}{\psi}<\infty
\right\}.
\]

\begin{definition}
We define the $\psi$-Wasserstein distance on $\MC_{\psi}\left(\RRP \right)$ by:
\[
W_{\psi}(\mu, \nu) = \sup_{\Wnorm{f}{\psi}\leq 1} \brac{f, \mu - \nu}.
\]
\end{definition}
The control $|f| \lesssim 1 + \psi$ ensures that the supremum is not equal to infinity for measures that do not have the same total mass. For two probability measures $\mu,\nu \in \PC_\psi(\RRP)$, this dual
definition agrees with the optimal-transport definition above. Indeed, by
Kantorovich--Rubinstein duality,
\[
\inf_{\pi\in\Pi(\mu,\nu)}
\int_{\RRP\times\RRP} d_\psi(x,y)\,\pi(\dd x,\dd y)
=
\sup_{\mathrm{Lip}_{d_\psi}(f)\le1}
\langle f,\mu-\nu\rangle .
\]
The additional growth condition in $\Wnorm{\cdot}{\psi}$ does not change the
supremum. Indeed since $\psi(x)\to0$ as $x\to\infty$, choose $x_\psi\in\RRP$ such
that $\psi(x_\psi)\le1$. For any $f$ with $\mathrm{Lip}_{d_\psi}(f)\le1$, the
function $g=f-f(x_\psi)$ satisfies
\[
|g(x)|\le d_\psi(x,x_\psi)\le 1+\psi(x),
\qquad
\mathrm{Lip}_{d_\psi}(g)\le1.
\]
Moreover $\langle g,\mu-\nu\rangle=\langle f,\mu-\nu\rangle$, since
$\mu$ and $\nu$ have the same total mass. Hence the two definitions coincide
on $\PC_\psi(\RRP)$.

For admissible weights $\psi,\chi$ and signed measures with the required moments, we shall use
\begin{equation}\label{eq:sum-weights}
    W_{\psi+\chi}(\mu,\nu)
    \leq W_\psi(\mu,\nu)+W_\chi(\mu,\nu)
    \leq 2W_{\psi+\chi}(\mu,\nu).
\end{equation}
Indeed, a test function $f\in\FC_{\psi+\chi}$ of norm at most one has a limit $f(\infty)$ with $|f(\infty)|\leq1$. Split its almost-everywhere derivative into
$f'|\psi'|/(|\psi'|+|\chi'|)$ and $f'|\chi'|/(|\psi'|+|\chi'|)$, and assign the constant $f(\infty)$ to the first part, note that $|\omega'| = -\omega'$ for every admissible weight $\omega$. Integrating from infinity gives two functions of norms at most one in $\FC_\psi$ and $\FC_\chi$. The reverse inequality follows by inclusion of the two unit balls in that of $\FC_{\psi+\chi}$.

\subsection{Derivatives on spaces of measures}
Our analysis relies heavily on differentiating functionals defined over spaces of measures. To this end, we employ the notion of flat derivatives, originating from optimal transport theory~\cite{santambrogio2015optimal} and widely used in the context of mean-field games~\cite{cardaliaguet2019master, carmona2018probabilistic} and related works~\cite{martiniKolmogorovEquationsSpaces2023}.
The definitions and chain rule we provide are formally classical; however, the mathematical rigor of the proofs depends on the choice of topology on the space of measures. Since our framework is based on a modified Wasserstein distance, we introduce the corresponding definitions in this setting. For every admissible weight $\psi$, we equip $\MC^+_\psi(\RRP)$ and $\MC_\psi(\RRP)$ with the topology induced by $W_\psi$, and $\RRP$ with the topology induced by $d_\psi$.

\begin{definition}\label{def:flat-derivative}
    Let $\psi$ be an admissible weight and let $F: \MC^+_{\psi}(\RRP) \to \RR$ be a continuous function. We say that $F$ has a flat derivative if there exists a function $\delta_{\mu}F :\MC^+_{\psi}(\RRP) \times \RRP \to \RR $ continuous in the product topology such that for all $\mu,\nu \in \MC^+_{\psi}(\RRP) $:
    \[
    x \mapsto \sup_{\theta \in [0,1]} |\delta_\mu F(\theta \mu + (1-\theta)\nu;x)| \in L^1(\mu+\nu),
    \]
    and
    \[ F(\mu) - F(\nu) = \int_0^1 \int_{\RRP} \delta_\mu F(\theta \mu + (1-\theta)\nu;x)\left(\mu - \nu\right)(\dd x)\dd \theta.\]
    We denote by $C^1(\MC_{\psi}^+(\RRP);\RR)$ the space of differentiable functions with a flat derivative.
\end{definition}
We also define flat derivatives for measure-valued maps. In this setting, the source and target may carry different admissible weights $\psi$ and $\chi$.
\begin{definition}
     Let $\psi,\chi$ be admissible weights and let $\varphi: \MC^+_{\psi}(\RRP) \to \MC^+_{\chi}(\RRP)$ be a continuous function. We say that $\varphi$ has a flat derivative if there exists a function $\delta_{\mu}\varphi :\MC^+_{\psi}(\RRP) \times \RRP \to \MC_{\chi}(\RRP)$ continuous in the product topology such that for all $\mu,\nu \in \MC^+_{\psi}(\RRP)$ and for all measurable $f$ such that $|f| \leq 1 + \chi$:
     \begin{itemize}
         \item The function
         \[x \mapsto \sup_{\theta \in [0,1]} |\langle f, \delta_\mu\varphi(\theta \mu + (1-\theta) \nu;x) \rangle| \in L^1(\mu+\nu).\]
        \item The following equality holds:
            \[ \left\langle f,\varphi(\mu) - \varphi(\nu) \right\rangle= \int_0^1 \int_{\RRP} \left\langle f,\delta_\mu \varphi(\theta \mu + (1-\theta)\nu,x)\right\rangle\left(\mu - \nu\right)(\dd x)\dd \theta.\]
     \end{itemize}
\end{definition}
We state the chain rule below.
\begin{proposition}\label{prop:chain-rule}
    Let $\psi,\chi,\omega$ be admissible weights with $1+\omega\leq C(1+\chi)$. Let $F : \MC_{\chi}^+(\RRP) \to \RR$ and $\varphi : \MC^+_{\psi}(\RRP) \to \MC_{\chi}^+(\RRP)$ have flat derivatives. Assume that $\delta_\mu\varphi$ is also jointly continuous with values in $(\MC_\omega,W_\omega)$ and that:
    \begin{itemize}
        \item For all $\eta \in \MC^+_{\chi}(\RRP)$,
        $\Wnorm{\delta_\mu F(\eta;\cdot)}{\omega}\leq u(\|\eta\|_\chi)$.
        \item For all $(\lambda,x)\in\MC^+_\psi(\RRP)\times\RRP$ and $f\in\FC_\omega$,
        \[
            |\brac{f,\delta_\mu\varphi(\lambda;x)}|
            \leq\Wnorm{f}{\omega}v(\|\lambda\|_\psi)(1+\psi(x)).
        \]
    \end{itemize}
    Here $u,v$ are non-decreasing positive functions on $\RR^+$. Then $F\circ\varphi$ has a flat derivative and
    \[
        \delta_\mu(F\circ\varphi)(\mu;x)
        =\brac{\delta_\mu F(\varphi(\mu);\cdot),\delta_\mu\varphi(\mu;x)}.
    \]
\end{proposition}
\begin{proof}
    Postponed to Appendix~\ref{appendix:chain-rule}.
\end{proof}

\subsection{The Smoluchowski coagulation equation}
We call a solution to the Smoluchowski coagulation equation ($\SCE$) of kernel $K_\alpha$ and initial value $\mu$, a map
\begin{align*}
    t \mapsto \mu_t : [0,T) \to \MC^+(\RRP)
\end{align*}
where $T \in (0,\infty]$, such that for all bounded measurable functions $f$:
\begin{equation}\tag{\ref{eq:SCE}}
\begin{aligned}
    \int_{\RRP} f(x)\mu_t (\dd x) =& \int_{\RRP} f(x)\mu (\dd x)\\
     &+ \dfrac12\int_0^t \int_{\RRP}\int_{\RRP} 
        K_{\alpha}(x,y)\left[f(x+y) - f(x) - f(y)\right] \mu_s(\dd x)\mu_s(\dd y) \dd s.
\end{aligned}
\end{equation}
With bracket notation the equation can be written:
\begin{align*}
    \dfrac{d}{dt}\langle f,\mu_t \rangle = \dfrac12\langle K_\alpha f , \mu_t \otimes \mu_t\rangle.
\end{align*}
where $K_\alpha f(x,y) = K_{\alpha}(x,y)\left(f(x+y) - f(x) - f(y) \right)$. Results on the well-posedness of this equation for the considered kernel are already known; see~\cite{norrisSmoluchowskisCoagulationEquation1999}, where the proof relies on the positivity of the solution. We use this well-posedness theory as an input and record only the estimates needed below.

All the proof of the following Lemmas and Propositions are done in Section \ref{section:well-posedness-SCE}.

\begin{lemma}\label{lem:well_posedness_smol_eq}
    Let $p \geq 2$ and $\mu \in \MC_{\psi_p}^+ (\RRP)$. For all $T \in (0,\infty]$, there exists a unique solution to the Smoluchowski Coagulation Equation \eqref{eq:SCE} with kernel $K_\alpha$; this solution is a continuous map from $[0,T)$ to $\left(\MC_{\psi_p}^+ (\RRP),\|\cdot\|_{\psi_p} \right)$.
\end{lemma}

We write $\varphi_t(\mu)$ for the solution at time $t$ with initial value $\mu$. The singularity of $K_\alpha$ has a second consequence which is specific to
this kernel. Small particles coagulate at a high rate, and this depletion near
the origin produces exponential moments in the singular variable
$x^{-\alpha}$. This regularization is the main deterministic mechanism used
in the Wasserstein estimates.
\begin{proposition}[Generation of Exponential Moments]\label{prop:exponential-moments}
    Let $\psi$ be a positive non-increasing function and let $\mu \in \MC^+_{\psi}(\RRP) \cap \MC^+_{\psi_2}(\RRP)$. Then, for all $t \geq 0$,
    \begin{align*}
        \int_{\RRP} \psi(x)\, e^{x^{-\alpha}\int_0^t \brac{1,\varphi_s(\mu)}\dd s}\, \varphi_t(\mu)(\dd x) \leq \brac{\psi,\mu}.
    \end{align*}
\end{proposition}

For the final result we need to know how functionals $F : \MC_{\psi_p}(\RRP) \to \RR$ evaluated in $\mu_t$ the solution to the $\SCE$ at time $t \geq 0$ evolve through time.
\begin{definition}
    For a functional $F$ with a flat derivative, we define the Smoluchowski coagulation operator at every $\mu$ in its domain for which the following integral is absolutely convergent:
    \begin{align*}
        \LC F(\mu) = \dfrac12 \langle K_\alpha \delta_\mu F(\mu;\cdot) , \mu\otimes \mu\rangle.
    \end{align*}
\end{definition}
The notation $\delta_\mu F$ is the flat derivative of $F$, defined in Definition~\ref{def:flat-derivative}. The coagulation operator drives this equation through the result of the next proposition.
\begin{proposition}
    Let $p\geq1$ and let $(\mu_t)_{t\in[0,T)}$ be a solution of the Smoluchowski Coagulation Equation \eqref{eq:SCE} with $\mu_0\in\MC^+_{\psi_{p+1}}(\RRP)$. For every $F \in C^1(\MC_{\psi_p}^+\left(\RRP \right);\RR)$ such that $|\delta_{\mu}F(\mu;z)| \leq C_F(\|\mu\|_{\psi_p})(1 + \psi_p(z))$, with $C_F$ non-decreasing, and every $t\in[0,T)$ (with a right derivative at zero), we have:
    \begin{align*}
        \dfrac{d}{dt} F( \mu_t) = \LC F(\mu_t).
    \end{align*}
\end{proposition}
\begin{proof}
    Let $t,s$ be two times. Then
    \begin{align*}
        F(\mu_t) - F(\mu_s) = \int_0^1 \int_{\RRP} \delta_\mu F(\left[\mu_t,\mu_s \right]^\theta;z) \left(\mu_t - \mu_s \right)(\dd z) \dd \theta.
    \end{align*}
    By hypothesis on $F$, the functional part of its derivative $\delta_\mu F(\mu,\cdot)$ is dominated by $C_F(\|\mu\|_{\psi_p})(1 + \psi_p)$, with $C_F$ non-decreasing. Hence, for all $s \geq 0$,
    \begin{align*}
        |\delta_\mu F(\left[\mu_t,\mu_s \right]^\theta;z)| \leq C_F(\|\mu_0\|_{\psi_p}) (1 + \psi_p(z)). 
    \end{align*}
    The $(p+1)$-moment assumption makes the SCE vector field continuous in the $\psi_p$-weighted total variation norm: the absolute integrand is bounded by $c(1+\psi_{p+1}(x))(1+\psi_{p+1}(y))$. Thus $t\mapsto\mu_t$ is continuously differentiable in this norm. The regularity of $\delta_\mu F$ in the measure variable and dominated convergence theorem give
    \begin{align*}
        \dfrac{\dd}{\dd t} F(\mu_t) = \int_{\RRP} \delta_\mu F(\mu_t;z) \mu'_t(\dd z) = \LC F(\mu_t).
    \end{align*}
\end{proof}
    
    \medskip
    A central aim in our analysis is to investigate how the solution $\varphi_t(\mu)$ depends on the initial condition $\mu$. Under stronger integrability assumptions, we show that the map $\mu \mapsto \varphi_t(\mu)$ is differentiable in a suitable sense. In particular we show in Section~\ref{section:well-posedness-SCE} that:
\begin{lemma}\label{lem:derivative-existence-and-eq}
    Let $T > 0$ and $p \geq 2$. For all $t \in [0,T)$, the function $ \mu \mapsto \varphi_t(\mu)$ has a flat derivative $\delta_\mu\varphi_t : \MC^+_{\psi_{p + 1}}(\RRP) \times \RRP \to \MC_{\psi_p}(\RRP) $. Furthermore, for all $z \in \RRP$ and $\mu \in \MC^+_{\psi_{p+1}}(\RRP)$, $\Proc{\delta_\mu\varphi_t(\mu;z)}$ is the unique solution to the equation defined for all bounded measurable functions $f$ by:
    \begin{align*}
        \brac{f,\delta_\mu\varphi_t(\mu;z)} = f(z) + \int_0^t \brac{K_{\alpha} f,\varphi_s(\mu)\otimes \delta_\mu\varphi_s(\mu;z)} \dd s.
    \end{align*}
\end{lemma}
We emphasize that for $\mu \in \MC^+_{\psi_{p+1}}(\RRP)$ and $z\in\RRP$, $\delta_\mu\varphi_t(\mu;z)$ is a signed measure. The SCE estimates and the proof of Lemma~\ref{lem:derivative-existence-and-eq} are given in Section~\ref{section:well-posedness-SCE}. 

\subsection{Construction of the Marcus-Lushnikov process}
We begin with a pathwise construction using exponential clocks. Assume that the initial system consists of $N$ particles with masses $x_1,\ldots,x_N$, and represent its state by the point measure
\[ 
    Z_0=\sum_{i=1}^N \delta_{x_i}. 
\]
More generally, suppose that the current state of the system is 
\[ 
    \mu=\sum_{i=1}^n \delta_{x_i}. 
\] 
To each pair $1\leq i<j\leq n$, we associate an independent exponential clock with parameter $K_\alpha(x_i,x_j)$. The waiting time until the next coagulation is therefore exponentially distributed with parameter 
\[ 
    \Lambda(\mu) := \sum_{1\leq i<j\leq n} K_\alpha(x_i,x_j). 
\] 
Conditionally on a coagulation occurring, the pair $(i,j)$ is selected with probability 
\[ 
    \frac{K_\alpha(x_i,x_j)}{\Lambda(\mu)}. 
\]
The selected particles merge into a single particle of mass $x_i+x_j$, and the system jumps from $\mu$ to 
\[ 
    \mu+\delta_{x_i+x_j}-\delta_{x_i}-\delta_{x_j}. 
\]
The procedure is then repeated from the new configuration, using new independent exponential clocks. If $\Lambda(\mu)=0$, the process remains in the state $\mu$. Since every coagulation decreases the number of particles by one, a system starting from $N$ particles undergoes at most $N-1$ jumps. This construction therefore defines a non-explosive càdlàg jump process, which we denote by $\Proc{M_t}$. Let 
\[ 
    \MC^{+,p}(\RRP) := \left\{ \sum_{i=1}^n \delta_{x_i} \;\middle|\; n\in\NN,\quad (x_1,\ldots,x_n)\in(\RRP)^n \right\} 
\] 
denote the space of finite point measures on $\RRP$. 

\begin{proposition}\label{prop:gen-SC} 
The process $\Proc{Z_t}$ is a Markov process with values in $\MC^{+,p}(\RRP)$. For every function $F\colon\MC^{+,p}(\RRP)\to\RR$ in the domain of its infinitesimal generator and every 
\[
\mu=\sum_{k=1}^n \delta_{x_k}, 
\] one has 
\begin{align*}
    L F(\mu) &= \sum_{1\leq i<j\leq n} K_\alpha(x_i,x_j) \left[ F\left( \mu+\delta_{x_i+x_j}-\delta_{x_i}-\delta_{x_j} \right) -F(\mu) \right]
\end{align*}
Equivalently, 
\begin{align*} 
    L F(\mu) &= \frac12 \int_{\RRP\times\RRP} K_\alpha(x,y) \left[ F\left( \mu+\delta_{x+y}-\delta_x-\delta_y \right) -F(\mu) \right] \left( \mu^{\otimes 2}-\Diag_{\#}\mu \right)(\dd x,\dd y), 
\end{align*} 
where 
\[ 
    \Diag\colon\RRP\longrightarrow\RRP\times\RRP, \qquad \Diag(x)=(x,x).
\] 
\end{proposition} 
\begin{proof} 
From a configuration $\mu=\sum_{i=1}^n\delta_{x_i}$, the transition 
\[ 
    \mu \longmapsto \mu+\delta_{x_i+x_j}-\delta_{x_i}-\delta_{x_j} 
\] 
occurs at rate $K_\alpha(x_i,x_j)$, which gives the first expression. Moreover,
\[ 
    \mu^{\otimes 2}-\Diag_{\#}\mu = \sum_{\substack{1\leq i,j\leq n\\ i\neq j}} \delta_{(x_i,x_j)}. 
\] 
The integral therefore sums over ordered pairs of distinct particles. Since both $K_\alpha$ and the coagulation operation are symmetric in $x$ and $y$, every unordered pair is counted twice, which yields the factor $1/2$.
\end{proof}
In order to compare the Marcus-Lushnikov process with the solution of the Smoluchowski coagulation equation, we introduce the rescaled process 
\[ 
    M_t^N:=\frac1N Z_{t/N}, \qquad t\geq0. 
\] 

Its state space is contained in 
\[ 
    \MC_N^{+,p}(\RRP) := \left\{ \frac1N\sum_{i=1}^n\delta_{x_i} \;\middle|\; 0\leq n\leq N,\quad (x_1,\ldots,x_n)\in(\RRP)^n \right\}. 
\] 
\begin{proposition}\label{prop:gen-rescaled-SC} 
    The process $\Proc{M_t^N}$ is a Markov process with infinitesimal generator \begin{multline*} 
        L^N F(\mu) = \frac N2 \int_{\RRP\times\RRP} K_\alpha(x,y) \left[ F\left( \mu+\frac1N \left(\delta_{x+y}-\delta_x-\delta_y\right) \right) -F(\mu) \right] \\ {}\times \left( \mu^{\otimes 2} -\frac1N\Diag_{\#}\mu \right)(\dd x,\dd y). \end{multline*} 
\end{proposition}
\begin{proof} 
    Let 
    \[ 
        \mu=\frac1N\sum_{i=1}^n\delta_{x_i}.
    \] 
    Because of the time rescaling $t\mapsto t/N$, the transition associated with the pair $(i,j)$ occurs at rate $K_\alpha(x_i,x_j)/N$. Hence \begin{align*} 
        L^N F(\mu) &= \frac1N \sum_{1\leq i<j\leq n} K_\alpha(x_i,x_j) \left[ F\left( \mu+\frac1N \left( \delta_{x_i+x_j} -\delta_{x_i} -\delta_{x_j} \right) \right) -F(\mu) \right]. 
    \end{align*} 
    Furthermore, 
    \[ 
        \mu^{\otimes 2} -\frac1N\Diag_{\#}\mu = \frac1{N^2} \sum_{\substack{1\leq i,j\leq n\\ i\neq j}} \delta_{(x_i,x_j)}.
    \] 
    Rewriting the preceding sum as an integral over ordered pairs gives the stated formula. 
\end{proof} 

\begin{remark} 
    When the initial masses are integers, the rescaled process is a version of the Marcus--Lushnikov process \cite{marcus1968stochastic,lushnikov1978coagulation}. It was introduced as a stochastic particle approximation of the Smoluchowski coagulation equation.
\end{remark}

\begin{proposition}\label{prop:SC_gen_differentiable}
    Let $\psi$ be an admissible weight and let $F \in C^1(\MC_{\psi}^+\left(\RRP \right);\RR)$. For all integers $N\geq1$ and $\mu \in \MC_N^{+,p}(\RRP)$, we have:
    \begin{multline*}
        L^N F(\mu) = \dfrac{1}{2}\int_0^1 \int_{\RRP \times \RRP} K_{\alpha}(x,y) \left\langle  \delta_\mu F\left(\mu + \frac{\theta}{N}\left(\delta_{x + y} - \delta_x - \delta_y \right),\cdot\right),\delta_{x+y} - \delta_x - \delta_y\right\rangle \\
        \left(\mu^{\otimes 2} - \dfrac1N\Diag_{\# }\mu \right)(\dd x,\dd y) \dd \theta.
    \end{multline*}
\end{proposition}
\begin{proof}
    For every pair \((x,y)\) contributing to 
    \[ 
        \mu^{\otimes 2}-\frac1N\Diag_{\#}\mu, 
    \] the measures 
    \[ 
        \mu+\frac{\theta}{N} \left(\delta_{x+y}-\delta_x-\delta_y\right), \qquad \theta\in[0,1], 
    \] 
    are positive. The result then follows by the definition of the flat derivative.
\end{proof}

\subsection{Main result}
In this section, we state the main assumptions and results. The first theorem provides a bound on the difference between observables when the Marcus-Lushnikov process is started from a general deterministic empirical measure. The bound depends on the Wasserstein distance between the two initial measures. While the result does not establish weak convergence—due to the restricted class of functionals considered—it is particularly relevant in applications where we are primarily interested in expected values of functionals of the process.

The second result is a corollary, where the Marcus-Lushnikov is initialized with a random empirical measure distributed according to a law $\mu$, and the Smoluchowski coagulation equation is started from the same law $\mu$. A particularly interesting feature of this result is that, unlike pathwise bounds which typically yield a convergence rate of $N^{-\frac{1}{2}}$ (see \cite{cepedaSmoluchowskisEquationRate2011a}), this result achieves a rate of $N^{-1}$.

\subsubsection*{Assumptions on \texorpdfstring{$F$}{F}}

The purpose of the assumptions below is to encode the differentiability and Lipschitz regularity of the observable in the geometry naturally associated with the singular kernel. 
\begin{assumptionA}\label{ass:observable}
Let $F : \MC^+_{\psi_2}(\RRP) \to \RR$ admit a flat derivative $\delta_\mu F : \MC^+_{\psi_2}(\RRP) \times \RRP \to \RR$, such that there exists a constant $c_F > 0$ such that for all $\mu, \tilde{\mu} \in \MC^+_{\psi_2}(\RRP)$ and $z, z_1, z_2 \in \RRP$, the following conditions hold:
\begin{itemize}
    \item $|\delta_\mu F(\mu; z)| \leq c_F(1 + \|\mu\|_{\psi_2}) (1 + \psi_1(z))$,
    \item $|\delta_\mu F(\mu; z_1) - \delta_\mu F(\mu; z_2)| \leq c_F (1 + \|\mu\|_{\psi_2}) d_{\psi_1}(z_1, z_2)$,
    \item $|\delta_\mu F(\mu; z) - \delta_\mu F(\tilde{\mu}; z)| \leq c_F(1 + \|\mu + \tilde\mu\|_{\psi_2})(1 + \psi_1(z)) W_{\psi_1}(\mu, \tilde{\mu})$,
    \item $|\delta_\mu F(\mu; z_1) - \delta_\mu F(\tilde{\mu}; z_1) - (\delta_\mu F(\mu; z_2) - \delta_\mu F(\tilde{\mu}; z_2))| \leq c_F  (1 + \|\mu + \tilde\mu\|_{\psi_2})W_{\psi_1}(\mu, \tilde{\mu}) d_{\psi_1}(z_1, z_2)$.
\end{itemize}

The origin of these assumptions lies in the Lipschitz regularity properties satisfied by the derivative of the solution to the $\SCE$; see Proposition~\ref{prop:wass-derivative}. Additional assumptions are motivated by the regularity described in Proposition~\ref{prop:wass-final}, and require that the map $z \mapsto \delta_\mu F(\mu, z)$ is differentiable with:
\begin{align*}
    \left| \frac{\dd}{\dd z} \delta_\mu F(\mu; z_1) - \frac{\dd}{\dd z} \delta_\mu F(\mu; z_2) \right| \leq c_F (1 + \|\mu\|_{\psi_2})d_{{|\psi_1'|}}(z_1, z_2).
\end{align*}

The constant $c_F$ is fixed as part of Assumption~$\A$, we denote it by $\Wnorm{\delta_\mu F}{\A}$.
\end{assumptionA}

\begin{remark}
    To better understand the nature of Assumption~$\A$, consider a linear
    functional of the form
    \[
        F(\mu)=\langle f,\mu\rangle.
    \]
    In this case, the flat derivative is independent of $\mu$ and is given by
    \[
        \delta_\mu F(\mu;z)=f(z).
    \]
    Hence the assumptions involving the dependence on the measure variable are
    automatically satisfied. The remaining conditions reduce to weighted
    growth and regularity assumptions on $f$. In the natural variable
    $r=z^{-\alpha}$, writing
    \[
        g(r)=f(r^{-1/\alpha}),
    \]
    for $g\in C^2((0,\infty))$, these conditions amount, up to multiplicative constants, to
    \[
        |g(r)|\lesssim 1+r,
        \qquad
        |g'(r)|\lesssim 1,
        \qquad
        r|g''(r)|\lesssim 1.
    \]
    Indeed,
    \[
        \frac{f''(z)}{\alpha(\alpha+1)z^{-\alpha-2}}
        =g'(r)+\frac{\alpha}{\alpha+1}r g''(r).
    \]
    The same characterization holds almost everywhere when $g'$ is locally Lipschitz.
    Thus, for linear observables, Assumption~$\A$ corresponds to a weighted
    $C^2$ regularity condition in the variable $z^{-\alpha}$. The same philosophy applies to the measure variable.
\end{remark}

\begin{remark}
    Consider a cylindrical function:
    \[
    F = \Psi\left(\brac{f_1, \mu}, \cdots, \brac{f_n,\mu} \right).
    \]
    Assume that for all $1 \le i \le n$ the function $g_i(r) = f_i(r^{-1/\alpha})$ satisfies:
    \[
    |g_i(r)|\lesssim 1+r,
        \qquad
        |g_i'(r)|\lesssim 1,
        \qquad
        r|g_i''(r)|\lesssim 1.
    \]
    And that
    \[
    \left|\nabla \Psi (\mathbf{x}) \right|_{1} \lesssim 1 + \left|\mathbf{x}\right|_{1}, \qquad \left|\nabla \Psi (\mathbf{x})  - \nabla \Psi (\mathbf{y})  \right|_{1} \lesssim |\mathbf{x} - \mathbf{y}|_{1},
    \]
    where $|\mathbf{x}|_{1} = |x_1| + \cdots + |x_n|$ for $\mathbf{x} \in \RR^n$. Then $F$ satisfies assumption $\A$.
\end{remark}
\subsubsection*{Main result}

Let $\mu^N$ be an empirical probability measure. We denote by $\Proc{\Phi_t(\mu^N)} := \Proc{M^N_t}$ the Marcus-Lushnikov process started from $\mu^N$ with kernel $K_\alpha$, and by $\Proc{\varphi_t(\mu)}$ the unique solution to the Smoluchowski coagulation equation started from $\mu$ with the same kernel $K_\alpha$. We state our first main theorem which deals with a deterministic initial measure for the Marcus-Lushnikov process.

\begin{theorem}\label{thm:main-result}
    Let $\mu^N \in \MC_N^{+,p}(\RRP)$ and $\mu \in \MC^+_{\psi_3}(\RRP)$, and let $F$ satisfy Assumption~$\A$. Then for all $t \in [0,T]$:
    \begin{multline*}
        \left| \mathbb{E}\left[F\left(\Phi_t\left(\mu^N\right)\right)\right] - F\left(\varphi_t(\mu)\right) - \int_{\RRP} \left\langle \delta_\mu F\left(\varphi_t \left(\mu\right);\cdot \right), \delta_\mu \varphi_t\left(\mu;z \right) \right\rangle \left(\mu^N - \mu\right)(\dd z) \right| \\
        \leq \Wnorm{\delta_\mu F}{\A} \left( \frac{ (1 + \|\mu^N\|_{\psi_3})C(T,\mu^N)}{N} + C(T,\mu + \mu^N)\left(W_{\psi_1}\left(\mu,\mu^N\right) + W_{\psi_2}\left(\mu,\mu^N\right) \right)^2\right).
    \end{multline*}
    where the constant is of the form
    \[
       C(T,\nu) =  c(1 + T)( 1 + \|\nu\|_{\psi_2})(1 + T\|\nu\|_{\psi_3}) e^{c T\|\nu\|_{\psi_2}(1 + T\|\nu\|)}.
    \]
\end{theorem}

The main appeal of this result is that, apart from the additional moment
needed to construct the flat derivative flow, it relies on the natural
integrability scale for the Smoluchowski equation. In particular, the assumption $\mu \in \MC_{\psi_2}(\RRP)$ is sufficient for the existence and uniqueness of solutions to the Smoluchowski coagulation equation (see \cite{norrisSmoluchowskisCoagulationEquation1999}). The stronger assumption $\mu \in \MC_{\psi_3}(\RRP)$ we require seems reasonable for the well-posedness of the derivative, which is necessary in our analysis.

To obtain the convergence rate of order $1/N$ from the estimates above, we need more moments on the initial law. The reason is twofold, first as an assumption required in order to use the Wasserstein rate result, second in controlling the random term $C(T,\mu^N)$.

\begin{theorem}\label{thm:main-result-2}
    Let $\varepsilon > 0$ and $\mu \in \PC_{\psi_{7 + \varepsilon}}(\RRP)$, let $\mu^N$ be an empirical measure with $N$ i.i.d. atoms of law $\mu$. Fix a time horizon $T > 0$, assume that $F$ satisfies Assumption~$\A$. Then for all $N \ge 1$ and for all $t \in [0, T]$, we have:
    \[
    \left| \mathbb{E}\left[F\left(\Phi_t\left(\mu^N\right)\right)\right] - F\left(\varphi_t(\mu)\right) \right| \leq \Wnorm{\delta_\mu F}{\A} \dfrac{c_{\varepsilon} \|\mu\|_{\psi_{7 + \varepsilon}}\|\mu\|_{\psi_2}e^{cT(1 + T)\|\mu\|_{\psi_2}}}{N} .
    \]
\end{theorem}

\subsection{Outline and proof strategies}
The proof of our main result relies on several key intermediate results, which naturally fall into two categories:

\textbf{1. Deterministic SCE mechanisms and differentiability.} The starting point is the existence and uniqueness of solutions to the $\SCE$, stated in Lemma~\ref{lem:well_posedness_smol_eq}, together with the exponential moment generation in Proposition~\ref{prop:exponential-moments}. Then, Lemma~\ref{lem:derivative-existence-and-eq} establishes that for all $t \geq 0$, the map $\mu \mapsto \varphi_t(\mu)$ is differentiable on $\MC^+_{\psi_3}(\RRP)$, with derivative denoted by $\delta_\mu \varphi_t$. Moreover, for all $\mu \in \MC^+_{\psi_3}(\RRP)$ and $z \in \RRP$, the flow $\Proc{\delta_\mu \varphi_t(\mu;z)}$ solves a linearized form of the $\SCE$. These results are collected in Section~\ref{section:well-posedness-SCE}.

\textbf{2. Lipschitz regularity.} The second group of results concerns stability of the $\SCE$ with respect to Wasserstein distances. In particular, both the solution $\varphi_t(\mu)$ and its linear derivative $\delta_\mu \varphi_t$ satisfy Lipschitz bounds with respect to the distance $W_{\psi_1 + \psi_2}$. These are proved in Propositions~\ref{prop:wass-sce}, \ref{prop:wass-derivative}, and \ref{prop:wass-final}, and are detailed in Section~\ref{section:Wass-regularity}.

\medskip

\textbf{Proof strategy overview.} The proof of the main convergence result is split into two theorems, reflecting two distinct sources of error:

\begin{itemize}
    \item In Theorem~\ref{thm:samesame}, we compare the Marcus-Lushnikov process with the deterministic solution to the $\SCE$ started from the same empirical measure $\mu^N$. This isolates the stochastic fluctuations and yields a $1/N$ convergence rate.
    \item In Theorem~\ref{thm:difdif}, we compare two deterministic trajectories: the solution to the $\SCE$ started from $\mu^N$ and the one started from the true law $\mu$. This captures the deterministic approximation error due to discretization of the initial condition.
\end{itemize}

We begin with the comparison between the Marcus-Lushnikov process and the deterministic $\SCE$ solution. 

\begin{theorem}\label{thm:samesame}
    Let $\mu^N \in \MC_{N}^{+,p}(\RRP)$ be a deterministic empirical measure and let $F$ satisfy Assumption~$\A$. For all $t \in [0,T]$, we have:
    \begin{align*}
        \left| \mathbb{E}\left[F\left(\Phi_t(\mu^N)\right)\right] - F\left(\varphi_t(\mu^N)\right) \right| 
        \leq \dfrac{\Wnorm{\delta_\mu F}{\A}(1 + \|\mu^N\|_{\psi_3}) C\left(T,\mu^N\right)}{N}.
    \end{align*}
    The constant is the same as in Theorem~\ref{thm:main-result}.
\end{theorem}

The proof is in Section \ref{sec:thm:samesame}, it is based on the following Duhamel identity, whose use was inspired by \cite{kolokoltsov2010central}, which is central to controlling the difference in expectations:
\begin{align}\label{eq:semi-group-relation}
    \mathbb{E}\left[F(\Phi_t(\mu^N))\right] - F(\varphi_t(\mu^N)) 
    = \int_0^t \mathbb{E}\left[\left(L^N - \LC\right)F\circ \varphi_s\left(\Phi_{t-s}(\mu^N)\right)\right]\, \mathrm{d}s.
\end{align}
Here, $L^N$ and $\LC$ denote the infinitesimal generators of the Marcus-Lushnikov process and the $\SCE$, respectively. The operator $L^N - \LC$ acts on the time-evolved test function $\mu \mapsto F(\varphi_s(\mu))$. 
\medskip

We now turn to the deterministic comparison between two solutions to the $\SCE$ started from different initial data.

\begin{theorem}\label{thm:difdif}
    Let $\mu,\nu \in \MC^+_{\psi_3}(\RRP)$ and let $F$ satisfy Assumption~$\A$. Then for all $t \in [0,T]$:
    \begin{multline*}
        \left| F\left(\varphi_t(\mu)\right) - F\left(\varphi_t(\nu)\right)
        - \int_{\RRP}\left\langle\delta_\mu F\left(\varphi_t(\nu);\cdot\right), \delta_\mu \varphi_t(\nu;z)
        \right\rangle
        (\mu-\nu)(\dd z)
        \right| \\
        \leq
        \Wnorm{\delta_\mu F}{\A}
        C(T,\mu+\nu)
        \left(W_{\psi_1}(\mu,\nu)+W_{\psi_2}(\mu,\nu)\right)^2 .
    \end{multline*}
    The constant $C$ is the same as in Theorem \ref{thm:main-result}.
\end{theorem}
The proof is in Section \ref{sec:thm:difdif}. Theorem~\ref{thm:main-result} follows by combining Theorems~\ref{thm:samesame} and~\ref{thm:difdif}. The proof of Theorem~\ref{thm:main-result-2} combines these estimates with a fourth-moment bound for empirical Wasserstein distances. Then the convergence rate of Corollary~\ref{cor:wass-square-bound} gives the required order; finally a localization argument allows to conclude with the initial law-moments assumptions. 
\section{The Smoluchowski coagulation equation and its derivative}\label{section:well-posedness-SCE}

This section collects the deterministic estimates for the Smoluchowski coagulation equation that will be used throughout the rest of the paper. We rely on the existing well-posedness theory for the kernel $K_\alpha$, and record only the a priori bounds, stability estimates, and cutoff approximation facts needed below. We then prove the generation of exponential moments in the variable $x^{-\alpha}$, a key ingredient for the Wasserstein analysis of Section~\ref{section:Wass-regularity}. Finally, we construct the flat derivative of the SCE flow and identify it as the unique signed-measure solution of the linearized equation. The common tool throughout is the use of weighted total variation estimates, which are particularly well suited to the derivative equation since the derivative evolves as a signed measure.
\subsection{A priori bounds and exponential regularization for the SCE}
We first record  deterministic facts on the nonlinear SCE needed below. The well-posedness statement is Lemma~\ref{lem:well_posedness_smol_eq}; the estimates in this subsection are the pieces of the proof that will be reused for the derivative equation and for the Wasserstein analysis.


\subsubsection{A priori bounds for solutions of the SCE}

In this section we assume the well-posedness of the Smoluchowski coagulation equation in $C([0,T); \MC^+_{\psi_2}(\RRP))$ for a finite time horizon $T > 0$. We recall that the solution of the $\SCE$ is a positive measure. The first lemma gives a stability result for solutions of the $\SCE$ and holds for all kernels.
\begin{lemma}\label{lem:sub-add-bound-smol}
    Let $\psi$ be sub-additive, meaning that $\psi(x+y)\leq \psi(x)+\psi(y)$ for all $x,y\in\RRP$. Assume $\mu \in \MC^+_{\psi}(\RRP) \cap \MC^+_{\psi_2}(\RRP)$ and denote by $\Proc{\mu_t}$ the solution of the SCE started from $\mu$. Then for all $t\in[0,T]$, $\mu_t \in \MC^+_{\psi}(\RRP)$ and:
    \begin{align*}
        \|\mu_t\|_{\psi} \leq \|\mu\|_{\psi}.
    \end{align*}
\end{lemma}
\begin{proof}
    The proof is straightforward, since $\mu_t$ is non-negative for all $t\geq 0$: $\|\mu_t\|_{\psi} = \brac{1 + \psi,\mu_t}$. Furthermore:
    \begin{align*}
        \|\mu_t\|_{\psi}
        &= \|\mu_0\|_{\psi} + \dfrac{1}{2}\int_0^t \int_{\RRP\times \RRP} K_{\alpha}(x,y) \left(-1 + \psi(x+y) - \psi(x) - \psi(y)\right) \mu_s(\dd x) \mu_s (\dd y )\dd s.
    \end{align*}
    The sub-additivity of $\psi$ yields the desired result.
\end{proof}
Note also that any non-increasing positive function is sub-additive, and positive linear combinations of sub-additive functions are also sub-additive. Therefore the weighted spaces used in the introduction are stable under the $\SCE$ flow. The next lemma shows that the solution of the $\SCE$ is Lipschitz with respect to its initial value in a specific weighted total variation norm.
\begin{lemma}\label{lem:Lip-smol-TV}
    Let $p \ge 1$, and $\mu,\tilde \mu \in \MC^+_{\psi_{p +1}}(\RRP)$ and denote $\Proc{\mu_t},\Proc{\tilde \mu_t}$ the respective solutions of the $\SCE$ started from $\mu,\tilde \mu$. We have:
    \begin{align*}
        \sup\limits_{t \in [0,T]}\| \mu_t - \tilde \mu_t\|_{\psi_p} \leq e^{cT\left(\|\mu\|_{\psi_{p+1}} + \|\tilde \mu\|_{\psi_{p+1}}\right)}\|\mu- \tilde \mu\|_{\psi_p}.
    \end{align*}
\end{lemma}
\begin{proof}
    For lighter notation in this proof we denote $1 + \psi_p$ by $\psi$. From Lemma \ref{lem:ODE_total_variation_kolo}:
    \begin{multline*}
        \| \mu_t - \tilde \mu_t\|_{\psi_p} = \|\mu- \tilde \mu\|_{\psi_p} \\
        + \dfrac12\int_0^t \int_{\RRP}\int_{\RRP} K_\alpha(x,y) \left[\psi\sigma_s(x+y) - \psi\sigma_s(x) - \psi\sigma_s(y) \right]\left(\mu_s + \tilde \mu_s\right)(\dd x) \left(\mu_s - \tilde \mu_s\right)(\dd y) \dd s.
    \end{multline*}
    where $\sigma_s$ is a density (or Radon-Nikodym derivative) of $\mu_s - \tilde \mu_s$ with respect to $|\mu_s - \tilde \mu_s|$. The function $\sigma_s$ is defined on $\RRP$ and takes values in $\left\lbrace -1,0,1 \right\rbrace $. Since $\Proc{\mu_t}$ and $\Proc{\tilde \mu_t}$ are positive and $|\sigma_t| \leq 1$ we have for the first term:
    \begin{multline*}
        \int_0^t \int_{\RRP}\int_{\RRP} K_\alpha(x,y) \psi\sigma_s(x+y)\left(\mu_s + \tilde \mu_s\right)(\dd x) \left(\mu_s - \tilde \mu_s\right)(\dd y) \dd s\\
        \leq \int_0^t \int_{\RRP}\int_{\RRP} K_\alpha(x,y) \psi(x+y)\left(\mu_s + \tilde \mu_s\right)(\dd x) \left|\mu_s - \tilde \mu_s\right|(\dd y) \dd s.
    \end{multline*}
    For the second term with the same arguments:
    \begin{multline*}
        -\int_0^t \int_{\RRP}\int_{\RRP} K_\alpha(x,y) \psi\sigma_s(x)\left(\mu_s + \tilde \mu_s\right)(\dd x) \left(\mu_s - \tilde \mu_s\right)(\dd y) \dd s\\
        \leq \int_0^t \int_{\RRP}\int_{\RRP} K_\alpha(x,y) \psi(x)\left(\mu_s + \tilde \mu_s\right)(\dd x) \left|\mu_s - \tilde \mu_s\right|(\dd y) \dd s.
    \end{multline*}
    Finally for the third term we use the fact that $\mu_t - \tilde \mu_t = \sigma_t\left| \mu_t - \tilde \mu_t \right|$ and we get:
    \begin{multline*}
        -\int_0^t \int_{\RRP}\int_{\RRP} K_\alpha(x,y) \psi\sigma_s(y)\left(\mu_s + \tilde \mu_s\right)(\dd x) \left(\mu_s - \tilde \mu_s\right)(\dd y) \dd s \\
        = -\int_0^t \int_{\RRP}\int_{\RRP} K_\alpha(x,y) \psi(y)\left(\mu_s + \tilde \mu_s\right)(\dd x) \left|\mu_s - \tilde \mu_s\right|(\dd y) \dd s.
    \end{multline*}
    In the end by using that $\psi(x+y) - \psi(y) \leq \psi(x)$:
    \begin{align*}
        \| \mu_t - \tilde \mu_t\|_{\psi_p} &\leq \|\mu- \tilde \mu\|_{\psi_p} + \int_0^t \int_{\RRP}\int_{\RRP} K_\alpha(x,y)\psi(x)\left(\mu_s + \tilde \mu_s\right)(\dd x) \left|\mu_s - \tilde \mu_s\right|(\dd y)\dd s\\
        &\leq \|\mu- \tilde \mu\|_{\psi_p} + c\left(\|\mu\|_{\psi_{p+1}} + \|\tilde \mu\|_{\psi_{p+1}}\right) \int_0^t \|\mu_s - \tilde \mu_s \|_{\psi_p}\dd s.
    \end{align*}
    We conclude thanks to Gronwall's lemma.
\end{proof}


\subsubsection{Well-posedness input and cutoff approximation}

We use the well-posedness theory of Norris~\cite{norrisSmoluchowskisCoagulationEquation1999}, together with the weighted estimates above. We will, however, use the standard cutoff approximation repeatedly, so we record it explicitly.

For $\varepsilon>0$, set $\mu^\varepsilon=\mathds{1}_{[\varepsilon,\infty)}\mu$, and let $\Proc{\mu_t^\varepsilon}$ solve
\begin{equation}\label{eq:cutoff_SCE}
    \brac{f, \mu^\varepsilon_t} = \brac{f,\mu^\varepsilon} + \dfrac12\int_0^t \brac{K_\alpha f,\left(\mu^\varepsilon_s\right)^{\otimes 2}} \dd s.
\end{equation}

\begin{proposition}[Cutoff approximation]\label{prop:well-posedness-cutoff-smol}
    Let $p\geq 2$ and $\mu\in\MC^+_{\psi_p}(\RRP)$. For every $\varepsilon>0$, the cutoff equation~\eqref{eq:cutoff_SCE} has a unique positive solution in $C([0,\infty);\MC^+_{\psi_p}(\RRP))$, with support included in $[\varepsilon,\infty)$. Moreover, for every finite $T>0$,
    \[
        \lim_{\varepsilon\to0}\sup_{t\in[0,T]}
        \|\mu_t^\varepsilon-\mu_t\|_{\psi_{p-1}}=0,
    \]
    where $\Proc{\mu_t}$ is the unique solution of the SCE with kernel $K_\alpha$ and initial condition $\mu$.
\end{proposition}

\begin{proof}
    For fixed $\varepsilon>0$, the kernel is bounded on $[\varepsilon,\infty)^2$, so the cutoff equation is a classical bounded-kernel coagulation equation. This gives existence and uniqueness of a positive solution. The support remains in $[\varepsilon,\infty)$, since the equation only removes particles at existing sizes and creates particles of size $x+y$ with $x,y\geq\varepsilon$.
    
    The cutoff convergence is obtained by applying Lemma~\ref{lem:Lip-smol-TV} to $\mu^\varepsilon$ and $\mu$, since
    \[
        \|\mu^\varepsilon-\mu\|_{\psi_{p-1}}
        \leq 2\varepsilon^\alpha\|\mu\|_{\psi_p}\qquad(0<\varepsilon\leq1).
    \]
\end{proof}

\begin{proof}[Proof of Lemma~\ref{lem:well_posedness_smol_eq}]
    Existence and uniqueness follow from the well-posedness theory for this kernel in~\cite{norrisSmoluchowskisCoagulationEquation1999}. The preservation of the weighted space follows from Lemma~\ref{lem:sub-add-bound-smol} applied to $\psi_p$. To check continuity in the $\psi_p$-norm, apply the same lemma to the non-increasing function $(1+\psi_p)\mathds1_{(0,a)}$. The resulting weighted tails of $\mu_t^\varepsilon$ and $\mu_t$ are bounded by the corresponding tail of $\mu$, uniformly in time and in $\varepsilon$. On $[a,\infty)$, the $\psi_p$-norm is bounded by a constant times the $\psi_{p-1}$-norm. The cutoff convergence therefore holds in the $\psi_p$-norm too, and continuity follows from that of the cutoff solutions.
\end{proof}

\subsubsection{Generation of exponential moments}\label{section:sce-exponential}

We now prove Proposition~\ref{prop:exponential-moments}. This result is specific to the kernel $K_\alpha(x,y)=x^{-\alpha}+y^{-\alpha}$: the singular coagulation rate removes small particles fast enough to produce exponential integrability in the variable $x^{-\alpha}$. This property is central to the Wasserstein stability of the solutions to the SCE and its derivative proved in Section~\ref{section:Wass-regularity}. It exploits an idea previously used to show that solutions to the $\SCE$ are positive in \cite{norrisSmoluchowskisCoagulationEquation1999} via the definition of the function:
\[
h_t(x)=\exp\left(\int_0^t\int K_\alpha(x,y)\,\mu_s(dy)\,ds\right).
\]
Here and below in this subsection, $\Proc{\mu_t}$ denotes the solution to the SCE started from $\mu$.
\begin{lemma}[Rescaled SCE and exponential moments]\label{lem:exponential-moments}
Let $\psi$ be positive and non-increasing, and assume the initial datum $\mu \in \MC^+_{\psi_2}(\RRP) \cap \MC^+_{\psi}(\RRP)$. Then, for every finite $T>0$, $\nu_t:=h_t\mu_t$ is well-defined in the weighted space
$\MC^+_{\psi}(\RRP)$, belongs to $C([0,T];\MC_{\psi}(\RRP))$, and solves
\[
\langle f,\nu_t\rangle
=
\langle f,\mu\rangle
+
\frac12\int_0^t
\iint K_\alpha(x,y)f(x+y)G_s(x,y)\,\nu_s(dx)\nu_s(dy)\,ds,
\]
where
\[
G_s(x,y)
=
\frac{h_s(x+y)}{h_s(x)h_s(y)}.
\]
Moreover,
\[
\langle \psi,h_t\mu_t\rangle
\le
\langle \psi,\mu\rangle
\exp\left(
\int_0^t \langle \psi_1,\mu_s\rangle\,ds
\right).
\]
\end{lemma}

A formal computation gives the inequality immediately. Since the cancellation involved in that computation may contain infinite terms, we use the cutoff approximation
$
\mu^\varepsilon = \mathds{1}_{[\varepsilon,\infty)}\mu
$
and prove the estimate first for $h_t^\varepsilon\mu_t^\varepsilon$. The family $h_t^\varepsilon\mu_t^\varepsilon$ is Cauchy in
$C([0,T];\MC^+_{\psi+\psi_1}(\RRP))$, its limit is $h_t\mu_t$, and the estimate passes to the limit. The details are given in Appendix~\ref{appendix:proof-exponential-moments}.
\begin{proof}
    See Appendix~\ref{appendix:proof-exponential-moments}.
\end{proof}

\begin{proof}[Proof of Proposition~\ref{prop:exponential-moments}]
    Lemma~\ref{lem:exponential-moments} gives
    \begin{align}\label{proof:eq:exponential-moment-prop-2}
        \int_{\RRP} \psi(x)\, e^{\int_0^t \int_{\RRP} K_{\alpha}(x,y) \varphi_s(\mu)(\dd y)\dd s}\, \varphi_t(\mu)(\dd x) \leq \brac{\psi,\mu} \, e^{\int_0^t \brac{\psi_1,\varphi_s(\mu)} \dd s}.
    \end{align}
    Since
    \[
        \int_{\RRP}K_\alpha(x,y)\varphi_s(\mu)(\dd y)
        =
        x^{-\alpha}\brac{1,\varphi_s(\mu)}
        +
        \brac{\psi_1,\varphi_s(\mu)},
    \]
    the common factor $e^{\int_0^t \brac{\psi_1,\varphi_s(\mu)}\dd s}$ cancels on both sides of~\eqref{proof:eq:exponential-moment-prop-2}, yielding the estimate.
\end{proof}


\subsection{Derivative}
In this section, we study the differentiability of the solution to the Smoluchowski coagulation equation (SCE) with respect to its initial condition. Let $ T > 0 $ be a fixed time horizon. The solution to the SCE can be viewed as a mapping
\[
\varphi: \MC^+_{\psi_p}(\RRP) \to C\left([0,T); \MC^+_{\psi_p}(\RRP)\right),
\]
where we denote this map by $ \mu \mapsto \Proc{\varphi_t(\mu)} $. Our goal is to understand the \emph{flat derivative} of $ \varphi_t $, which we denote by $ \delta_\mu \varphi_t(\mu, z) $.

A formal differentiation of the Smoluchowski coagulation equation suggests that $ t \mapsto \delta_\mu \varphi_t(\mu, z) $ satisfies the following weak formulation:
\[
\brac{f, \delta_\mu \varphi_t(\mu, z)} = f(z) + \int_0^t \brac{K_\alpha f, \varphi_s(\mu) \otimes \delta_\mu \varphi_s(\mu, z)} \dd s,
\]
for all bounded measurable test functions $ f $. In other words, for all $ \mu \in \MC^+_{\psi_p}(\RRP) $ and $z\in\RRP$, the family $ \Proc{\delta_\mu \varphi_t(\mu, z)} $ should solve the equation:
\begin{equation} \label{eq:derivative-initial-general}
    \brac{f, m_t} = \brac{f, \delta_z} + \int_0^t \brac{K_\alpha f, \varphi_s(\mu) \otimes m_s} \dd s.
\end{equation}

This section is organized as follows. First we give a priori bounds on equation \eqref{eq:derivative-initial-general}. Then we show that this equation is well-posed. Finally, we show that for all $t \in [0,T)$, $\varphi_t$ is differentiable and that $t \mapsto \delta_\mu \varphi_t(\mu,z)$ is the unique solution of \eqref{eq:derivative-initial-general}.

\subsubsection{A priori bounds for solutions of the linearized SCE}

The estimates below are the replacement for positivity in the linearized equation. They control signed solutions through their total variation and will be used both to construct the derivative flow and to justify the stability estimates with respect to the initial measure. We denote by $\Proc{\mu_t}$ the solution of the $\SCE$ with initial value $\mu$ and kernel $K_\alpha$.
\begin{lemma}\label{lem:bound_derivative}
    Let $p\ge1$, $\mu\in\MC^+_{\psi_{p+1}}$, and $m\in\MC_{\psi_p}$. Any solution $\Proc{m_t}$ of the linearized equation
    \[
        \brac{f,m_t}=\brac{f,m}+\int_0^t\brac{K_\alpha f,\mu_s\otimes m_s}\dd s
    \]
    on the time interval $[0,T]$ satisfies
    \[ \sup\limits_{t \in [0,T]}\| m_t\|_{\psi_p} \leq  \| m\|_{\psi_p} \exp{\left( cT\| \mu\|_{\psi_{p+1}}\right)}.\]
\end{lemma}
\begin{proof}
    Lemma \ref{lem:ODE_total_variation_kolo} applied to the signed equation gives
    \begin{align*}
        \| m_t\|_{\psi_p} = \| m\|_{\psi_p} + \int_0^t \brac{K_\alpha (1 + \psi_p)\sigma_s, \mu_s\otimes m_s} \dd s
    \end{align*}
    where $\sigma_t$ is a density of $m_t$ with respect to $|m_t|$ and takes values in $\lbrace -1,0,1 \rbrace$. The integrand is bounded as in the proof of Lemma \ref{lem:Lip-smol-TV}; using also Lemma \ref{lem:sub-add-bound-smol},
    \begin{align*}
        \brac{K_\alpha (1 + \psi_p)\sigma_s, \mu_s\otimes m_s}  \leq c\brac{1 + \psi_{p+1},\mu} \| m_s\|_{\psi_p}.
    \end{align*}
    Gronwall's lemma gives the estimate.
\end{proof}
We next compare two initial values in Equation \eqref{eq:derivative-initial-general} driven by two different $\SCE$ solutions.
\begin{lemma}\label{lem:pseudo-lip-derivative-eq}
     Let $p \ge 2$, let $\Proc{\mu_t}$ and $\Proc{\tilde{\mu}_t}$ be two solutions of the $\SCE$ with respective initial measures $\mu,\tilde \mu \in \MC^+_{\psi_{p + 1}}$. We denote by $\Proc{m_t}$ and $\Proc{\tilde{m}_t}$ two solutions of the linearized equation driven respectively by $\Proc{\mu_t}$ and $\Proc{\tilde{\mu}_t}$, with initial values $m,\tilde m\in\MC_{\psi_p}(\RRP)$. The following holds:
    \begin{multline*}
        \sup\limits_{t \in [0,T]} \| m_t - \tilde{m}_t \|_{\psi_{p - 1}} \leq \|m- \tilde{m} \|_{\psi_{p - 1}} C_1\left(T;\| \mu\|_{\psi_p} + \| \tilde{\mu}\|_{\psi_p}\right)\\
        + \|\mu - \tilde{\mu} \|_{\psi_p}C_2\left(T; \|m\|_{\psi_p}+ \|\tilde{m}\|_{\psi_p};\|\mu\|_{\psi_{p +1}} + \|\tilde{\mu}\|_{\psi_{p +1}}\right),
    \end{multline*}
    with constants of the form:
    \begin{align*}
        C_1 &= e^{cT\left(\| \mu\|_{\psi_p} + \|\tilde{ \mu}\|_{\psi_p}\right)}\\
        C_2 &= cT(\|m\|_{\psi_p}+ \|\tilde{m}\|_{\psi_p})e^{cT\left(\| \mu\|_{\psi_{p+1}} + \|\tilde{ \mu}\|_{\psi_{p+1}}\right)}.
    \end{align*}
\end{lemma}
\begin{proof}
    For lighter notation set $\psi=1+\psi_{p-1}$. Lemma \ref{lem:ODE_total_variation_kolo} gives a density $\sigma_t$ such that $m_t - \tilde{m}_t = \sigma_t|m_t - \tilde{m}_t|$ and
    \begin{multline*}
        \| m_t - \tilde{m}_t \|_{\psi_{p-1}} = \|m-\tilde{m} \|_{\psi_{p-1}} + \int_0^t \brac{K_\alpha\psi\sigma_s, \mu_s\otimes m_s - \tilde{\mu}_s\otimes \tilde{m}_s} \dd s \\
        = \|m-\tilde{m} \|_{\psi_{p-1}} + \dfrac12\int_0^t \brac{K_\alpha\psi\sigma_s, (\mu_s - \tilde{\mu}_s)\otimes (m_s + \tilde{m}_s) } \\
        + \dfrac12\brac{K_\alpha\psi\sigma_s, (\mu_s + \tilde{\mu}_s)\otimes (m_s - \tilde{m}_s)} \dd s.
    \end{multline*}
    For the first term we have for all $s \in [0,t]$:
    \begin{align*}
        \left| \brac{K_\alpha\psi\sigma_s, (\mu_s - \tilde{\mu}_s)\otimes (m_s + \tilde{m}_s)}\right| \leq c\| \mu_s - \tilde{\mu}_s\|_{\psi_p} \left(\| m_s\|_{\psi_p} + \| \tilde{m}_s\|_{\psi_p}\right).
    \end{align*}
    Lemmas \ref{lem:bound_derivative} and \ref{lem:Lip-smol-TV} imply
    \begin{multline*}
        \left| \brac{K_\alpha\psi\sigma_s, (\mu_s - \tilde{\mu}_s)\otimes (m_s + \tilde{m}_s) }\right| \\
        \leq c(\|m\|_{\psi_p}+ \|\tilde{m}\|_{\psi_p})e^{cT\left(\|\mu\|_{\psi_{p +1}} + \|\tilde{\mu}\|_{\psi_{p +1}}\right)} \|\mu - \tilde{\mu}\|_{\psi_p}.
    \end{multline*}
    For the second term, the definition of $\sigma_s$ gives, for all $s \in [0,t]$,
    \begin{align*}
         \brac{K_\alpha\psi\sigma_s, (\mu_s + \tilde{\mu}_s)\otimes (m_s - \tilde{m}_s)} 
        &\leq c\left(\| \tilde{\mu}_s\|_{\psi_p} + \| \mu_s\|_{\psi_p}\right)\|m_s - \tilde{m}_s \|_{\psi_{p-1}} \\
        &\leq c\left(\| \tilde{\mu}\|_{\psi_p} + \| \mu\|_{\psi_p}\right)\|m_s - \tilde{m}_s \|_{\psi_{p-1}}.
    \end{align*}
    The result then follows by Gronwall's lemma.
\end{proof}

\subsubsection{Well-posedness of the linearized derivative's equation}
The goal of the following is to show that the equation followed by the derivative of the $\SCE$ solutions is well-posed.
\begin{lemma}\label{lem:well_posedness_derivative_eq}
    Let $p \geq 2$ and denote by $\Proc{\mu_t}$ the solution of the Smoluchowski coagulation equation started from $\mu \in \MC^+_{\psi_{p + 1}}\left(\RRP\right)$. For all $T > 0$ and all $m \in \MC_{\psi_p}\left(\RRP\right)$, the equation defined for all bounded measurable $f$ by
    \begin{equation}\label{eq:derivative-eq-1}
        \brac{f,m_t} = \brac{f,m} + \int_0^t \brac{K_\alpha f,\mu_s\otimes m_s} \dd s.
    \end{equation}
    has a unique solution in $C\left([0,T];\MC_{\psi_p}\left(\RRP\right)\right)$.
\end{lemma}

\begin{remark}
    Our proof requires us to take the solution of the Smoluchowski equation in $\MC^+_{\psi_{p+1}}(\RRP)$ to ensure the well-posedness of \eqref{eq:derivative-eq-1} in $\MC_{\psi_{p}}(\RRP)$. This is reminiscent of the fact that this equation is the one verified by its derivative.
\end{remark}
Our strategy is to use a cut-off of the initial measure. We consider $\Proc{\mu^\varepsilon}$, the solution to the $\SCE$ with initial value $\mu^\varepsilon = \mu\mathds{1}_{[\varepsilon,\infty)}$, and the cut-off initial value $m^\varepsilon = m\mathds{1}_{[\varepsilon,\infty)}$ for \eqref{eq:derivative-eq-1}. We then consider the family of equations below and prove that the corresponding solutions converge.
\begin{equation}\label{eq:cut-off-derivative-equation}
    \brac{f,m^\varepsilon_t} = \brac{f,m^\varepsilon} + \int_0^t \brac{K_\alpha f,\mu^\varepsilon_s\otimes m^\varepsilon_s} \dd s.
\end{equation}
We start by showing that \eqref{eq:cut-off-derivative-equation} is well-posed on a bounded subset of $\MC_{\psi_p}\left(\RRP\right)$. For $M \geq 0$ define:
\begin{align*}
    \MC^{M,\varepsilon}_{\psi_p}\left(\RRP\right) = \left\lbrace m\in \MC_{\psi_p}\left(\RRP\right): \|m\|_{\psi_p} \leq M,\ \mathrm{supp}(m) \subset [\varepsilon,\infty) \right\rbrace.
\end{align*}
The proof is a Banach contraction argument, but here the linearity of the equation gives a local time depending on $\varepsilon$ and $\mu$, not on the radius $M$.
\begin{proposition}\label{prop:well-posedness-derivative-cutoff-bounded}
    Let $\mu \in \MC^+_{\psi_{p+1}}\left(\RRP\right)$, let $\Proc{\mu^\varepsilon_t}$ be the cutoff SCE solution started from $\mu^\varepsilon$, and let $m^\varepsilon \in \MC^{M/2,\varepsilon}_{\psi_p}\left(\RRP\right)$. There exists a time $T^\varepsilon>0$, independent of $M$, such that Equation~\eqref{eq:cut-off-derivative-equation} has a unique solution in $C\left([0,T^\varepsilon];\MC^{M,\varepsilon}_{\psi_p}\left(\RRP\right)\right)$. Consequently the cutoff equation has a unique solution in $C([0,T];\MC_{\psi_p}(\RRP))$ for every finite $T>0$.
\end{proposition}

\begin{proof}
    Let $\psi=1+\psi_p$, and define, for $\eta \in C([0,T];\MC^{M,\varepsilon}_{\psi_p})$,
    \[
        \brac{f,\Gamma^\varepsilon(\eta)_t}
        =
        \brac{f,m^\varepsilon}
        +
        \int_0^t \brac{K_\alpha f,\mu_s^\varepsilon\otimes \eta_s}\dd s .
    \]
    Since $\mu_s^\varepsilon$ and $\eta_s$ are supported in $[\varepsilon,\infty)$ and the coalescence operator increases mass, the image also has support in $[\varepsilon,\infty)$. Moreover, for $|f|\leq \psi$,
    \[
        \left|\brac{K_\alpha f,\mu_s^\varepsilon\otimes \eta_s}\right|
        \leq
        c\varepsilon^{-\alpha}\|\mu_s^\varepsilon\|_{\psi_p}\|\eta_s\|_{\psi_p}
        \leq
        c\varepsilon^{-\alpha}\|\mu\|_{\psi_p}M,
    \]
    where we used Lemma~\ref{lem:sub-add-bound-smol} for $\mu_s^\varepsilon$ and the defining bound $\sup_s\|\eta_s\|_{\psi_p}\leq M$. Therefore $\Gamma^\varepsilon$ maps the ball into itself as soon as
    $cT\varepsilon^{-\alpha}\|\mu\|_{\psi_p}\leq 1/2$.

    The same estimate applied to $\eta^1-\eta^2$ gives
    \[
        \sup_{t\in[0,T]}
        \|\Gamma^\varepsilon(\eta^1)_t-\Gamma^\varepsilon(\eta^2)_t\|_{\psi_p}
        \leq
        cT\varepsilon^{-\alpha}\|\mu\|_{\psi_p}
        \sup_{t\in[0,T]}\|\eta^1_t-\eta^2_t\|_{\psi_p}.
    \]
    Thus $\Gamma^\varepsilon$ is a contraction for $T=T^\varepsilon$ small enough. The space $C([0,T^\varepsilon];\MC^{M,\varepsilon}_{\psi_p})$ is complete, since $\left( \MC_{\psi_p}(\RRP), \|\cdot\|_{\psi_p}\right)$ is Banach and the support and radius conditions are closed. Banach's contraction theorem gives a unique local solution.

    The choice of $T^\varepsilon$ is independent of $M$. Restarting the equation and using Lemma~\ref{lem:bound_derivative} extends the cutoff solution to every finite time horizon.
\end{proof}
It remains to pass from the cutoff equation to the singular equation.
\begin{proof}[Proof of Lemma \ref{lem:well_posedness_derivative_eq}]
    Fix $T>0$ and take $\varepsilon_n=1/n$. Let $\Proc{m^{\varepsilon_n}_t}$ be the cutoff solutions constructed above. We first show that this sequence is Cauchy in $C([0,T];\MC_{\psi_{p-1}}(\RRP))$. By Lemma~\ref{lem:pseudo-lip-derivative-eq}, for $\varepsilon,\rho>0$,
    \begin{multline*}
        \sup_{t\in[0,T]}\|m^\varepsilon_t-m^\rho_t\|_{\psi_{p-1}}
        \leq
        C(T,\|\mu\|_{\psi_p})\|m^\varepsilon-m^\rho\|_{\psi_{p-1}}
        +
        C(T,\|m\|_{\psi_p},\|\mu\|_{\psi_{p+1}})
        \|\mu^\varepsilon-\mu^\rho\|_{\psi_p}.
    \end{multline*}
    The two terms on the right vanish as $\varepsilon,\rho\to0$, because truncation at the origin is continuous from $\MC_{\psi_p}$ to $\MC_{\psi_{p-1}}$, and from $\MC_{\psi_{p+1}}$ to $\MC_{\psi_p}$. Hence the cutoff solutions converge to some $\Proc{m_t}\in C([0,T];\MC_{\psi_{p-1}}(\RRP))$.

   We now identify the limit. By homogeneity, it suffices to consider bounded measurable $f$ with $\|f\|_\infty\leq1$. Then
    \begin{multline*}
        \left|\brac{f,m_t} - \brac{f,m}
        - \int_0^t \brac{K_\alpha f,\mu_s\otimes m_s}\dd s\right|\\
        \leq
        \|m_t-m^\varepsilon_t\|_{\psi_0}
        +\|m-m^\varepsilon\|_{\psi_0}
        +c\int_0^t
        \left(
        \|\mu_s-\mu_s^\varepsilon\|_{\psi_1}\|m_s\|_{\psi_1}
        +
        \|\mu_s^\varepsilon\|_{\psi_1}\|m_s-m_s^\varepsilon\|_{\psi_1}
        \right)\dd s .
    \end{multline*}
    The first two terms vanish by the convergence just proved and by truncation of $m$. The integral vanishes by Proposition~\ref{prop:well-posedness-cutoff-smol}, Lemma~\ref{lem:bound_derivative}, and the same Cauchy estimate applied in the $\psi_1$-norm. Therefore $\Proc{m_t}$ solves~\eqref{eq:derivative-eq-1}.

    Uniqueness follows from Lemma~\ref{lem:pseudo-lip-derivative-eq} with the same SCE solution $\Proc{\mu_t}$ and the same initial value. Finally, Lemma~\ref{lem:bound_derivative} applied to the cutoff solutions and then passed to the limit gives
    \[
        \sup_{t\in[0,T]}\|m_t\|_{\psi_p}
        \leq
        \|m\|_{\psi_p}e^{cT\|\mu\|_{\psi_{p+1}}},
    \]
    To obtain continuity in this stronger norm, set $w_a=(1+\psi_p)\mathds{1}_{(0,a)}$. The sign estimate used in Lemma~\ref{lem:bound_derivative}, applied first to the cutoff solutions, yields
    \[
        \sup_{t\leq T}\brac{w_a,|m_t^\varepsilon|}
        \leq\brac{w_a,|m|}
        +cT\sup_{s,\varepsilon}\|m_s^\varepsilon\|_{\psi_1}
        \int_{(0,a)}(1+\psi_{p+1})\,\dd\mu.
    \]
    Indeed, $w_a(x+y)\leq w_a(y)$ cancels the loss in the signed variable, and the remaining term is bounded by $\brac{w_a K_\alpha,\mu_s^\varepsilon\otimes|m_s^\varepsilon|}$. The corresponding weighted tails of $\mu_s^\varepsilon$ do not increase. The right-hand side tends to zero as $a\downarrow0$, uniformly in the cutoff. On $[a,\infty)$, the $\psi_p$-norm is bounded by a constant depending on $a$ times the $\psi_{p-1}$-norm. Therefore the cutoff solutions converge uniformly in time in the $\psi_p$-norm as well. This proves that the limit belongs to $C([0,T];\MC_{\psi_p}(\RRP))$.
\end{proof}

\subsubsection{Existence of the derivative}
In this part $p\geq3$, and the initial measures are positive elements of $\MC^+_{\psi_p}(\RRP)$. The index is shifted by one relative to Lemma~\ref{lem:derivative-existence-and-eq}: the derivative takes values in $\MC_{\psi_{p-1}}(\RRP)$. We prove that the solution map admits a flat derivative and identify this derivative with the linearized equation.
\begin{lemma}[Restated Lemma~\ref{lem:derivative-existence-and-eq}]
    Let $T > 0$. For all $t \in [0,T)$, the map $ \mu \mapsto \varphi_t(\mu)$ has a flat derivative denoted by $\delta_\mu\varphi_t$. Moreover, for all $z \in \RRP$ and $\mu \in \MC^+_{\psi_p}(\RRP)$, the process $\Proc{\delta_\mu\varphi_t(\mu,z)}$ is the unique solution to Equation~\eqref{eq:derivative-initial-general} with initial value $\delta_z$.
\end{lemma}

For a base measure $\lambda\in\MC^+_{\psi_p}(\RRP)$, we denote by $\Proc{\psi_t^\lambda(m)}$ the solution of the linearized equation along $\Proc{\varphi_t(\lambda)}$, with initial value $m\in\MC_{\psi_{p-1}}(\RRP)$:
\[
    \brac{f,\psi_t^\lambda(m)}
    =
    \brac{f,m}
    +\int_0^t \brac{K_\alpha f,\varphi_s(\lambda)\otimes \psi_s^\lambda(m)}\dd s .
\]
The following elementary linearity property will be used to turn directional derivatives into flat derivatives.
\begin{proposition}\label{prop:derivative-prop-fubini}
    For all $t\geq0$, $m\in\MC_{\psi_{p-1}}(\RRP)$, and bounded measurable $f$,
    \[
        \brac{f,\psi_t^\lambda(m)}
        =
        \int_{\RRP}\brac{f,\psi_t^\lambda(\delta_z)}m(\dd z).
    \]
\end{proposition}
\begin{proof}
    It is enough to justify Fubini. By Lemma~\ref{lem:bound_derivative}, for $s\leq T$,
    \[
        \|\psi_s^\lambda(\delta_z)\|_{\psi_1}
        \leq (1+\psi_1(z))e^{cT\|\lambda\|_{\psi_2}}.
    \]
    Hence, for fixed bounded $f$,
    \[
        \int_{\RRP}\int_0^t
        \brac{|K_\alpha f|,\varphi_s(\lambda)\otimes|\psi_s^\lambda(\delta_z)|}
        \dd s\, |m|(\dd z)
        \leq
        C(T,\lambda)\|f\|_\infty\|m\|_{\psi_1}<\infty.
    \]
    Fubini then shows that
    \[
        \tilde m_t:=\int_{\RRP}\psi_t^\lambda(\delta_z)m(\dd z)
    \]
    satisfies the same linearized equation as $\psi_t^\lambda(m)$, with initial value $m$. Uniqueness in Lemma~\ref{lem:well_posedness_derivative_eq} gives the claim.
\end{proof}

\begin{proof}[Proof of Lemma \ref{lem:derivative-existence-and-eq}]
    For $\theta\in[0,1]$, write
    \[
        \lambda_\theta=[\mu,\nu]^\theta=\theta\mu+(1-\theta)\nu.
    \]
    Let $\mu,\nu\in\MC^+_{\psi_p}(\RRP)$. We first show that, for all bounded measurable $f$, the map
    \begin{align*}
        (0,1) &\to \RR \\
        \theta &\mapsto \brac{f,\varphi_t(\lambda_\theta)}
    \end{align*}
    is differentiable, with derivative $\brac{f,\psi_t^{\lambda_\theta}(\mu-\nu)}$. Let $h$ be such that $\theta+h\in(0,1)$, and set
    \[
        r_t^h
        =
        h^{-1}\left(\varphi_t(\lambda_{\theta+h})-\varphi_t(\lambda_\theta)\right)
        -
        \psi_t^{\lambda_\theta}(\mu-\nu).
    \]
    Subtracting the two SCEs and the linearized equation gives
    \[
        \brac{f,r_t^h}
        =
        \int_0^t \brac{K_\alpha f,\varphi_s(\lambda_\theta)\otimes r_s^h}\dd s
        +
        \frac{h^{-1}}{2}\int_0^t
        \brac{K_\alpha f,\left(\varphi_s(\lambda_{\theta+h})-\varphi_s(\lambda_\theta)\right)^{\otimes2}}\dd s .
    \]
    Lemma~\ref{lem:Lip-smol-TV} bounds the quadratic remainder by
    $|h|C(T,\mu+\nu)\|\mu-\nu\|_{\psi_1}$. Therefore, by Lemma~\ref{lem:ODE_total_variation_kolo},
    \[
        \|r_t^h\|_{\psi_1}
        \leq
        c\int_0^t \|\varphi_s(\lambda_\theta)\|_{\psi_2}\|r_s^h\|_{\psi_1}\dd s
        +
        |h| C(T,\mu+\nu)\|\mu-\nu\|_{\psi_1}.
    \]
    Gronwall's lemma gives $\sup_{t\leq T}\|r_t^h\|_{\psi_1}\to0$ as $h\to0$, proving the differentiability in the direction $\mu-\nu$.

    Integrating along the segment gives
    \begin{align*}
        \brac{f,\varphi_t(\mu) - \varphi_t(\nu)} &= \int_0^1 \partial_\theta \brac{f,\varphi_t(\lambda_\theta)} \dd \theta \\
        &= \int_0^1 \brac{f,\psi_t^{\lambda_\theta}(\mu-\nu)} \dd \theta \\
        &= \int_0^1 \int_{\RRP} \brac{f,\psi_t^{\lambda_\theta}(\delta_z)}(\mu - \nu)(\dd z) \dd \theta,
    \end{align*}
    where the last equality follows from Proposition~\ref{prop:derivative-prop-fubini}. Thus
    $\delta_\mu\varphi_t(\lambda;z)=\psi_t^\lambda(\delta_z)$, and it solves the linearized equation. The identity extends from bounded tests to $|f|\leq1+\psi_{p-1}$ by truncation and the bound
    \[
        \|\psi_t^\lambda(\delta_z)\|_{\psi_{p-1}}
        \leq(1+\psi_{p-1}(z))e^{cT\|\lambda\|_{\psi_p}}.
    \]
    We record how to obtain the joint continuity required by the definition. The estimates in Section~\ref{section:Wass-regularity} use only the linearized equation and the segment identity just proved. Propositions~\ref{prop:wass-derivative} and~\ref{prop:wass-final} therefore give joint continuity of $(\lambda,z)\mapsto\psi_t^\lambda(\delta_z)$ in $W^1_{\psi_1}$ for the source topology $W_{\psi_p}$. For a convergent sequence of positive source measures, their $\psi_p$-tails are uniformly integrable. The tail estimate in the preceding proof, with index $p-1$, then makes the $\psi_{p-1}$-tails of the derivative measures uniformly integrable as well; the initial atoms $z$ stay in a compact subset of $(0,\infty)$. On a bounded interval in the variable $r=x^{-\alpha}$, the unit ball of $\FC_{\psi_{p-1}}$ can be uniformly approximated by smooth functions with a common finite $\FC^1_{\psi_1}$-norm. Truncation in $r$ and this approximation upgrade the continuity to $W_{\psi_{p-1}}$. The same argument, using the non-increasing weighted tails of the positive SCE solutions and Proposition~\ref{prop:wass-sce}, proves continuity of the flow itself. This completes the flat-derivative assertion, including its topology.
\end{proof}

\subsection{Rescaled derivative}
The Wasserstein estimates of the next section use the same exponential change
of variables for the derivative equation as for the nonlinear SCE. We record
only the identity needed there; the quantitative estimates for this rescaled
flow are proved later in Proposition~\ref{lem:wass-rescaled-derivative}.

Let $\mu\in\MC^+_{\psi_3}(\RRP)$, write $\Proc{\mu_t}:=\Proc{\varphi_t(\mu)}$,
and let $\Proc{m_t}$ be the solution of the linearized equation
\[
    \brac{f,m_t}
    =
    \brac{f,m_0}
    +
    \int_0^t \brac{K_\alpha f,\mu_s\otimes m_s}\dd s .
\]
Define
\[
    h_t(x)=\exp\left(\int_0^t\brac{K_\alpha(x,\cdot),\mu_s}\dd s\right),
    \qquad
    G_t(x,y)=h_t(x+y)h_t^{-1}(x)h_t^{-1}(y).
\]
For $m_0 \in \MC_{\psi_2}$, the rescaled signed
flow $\Proc{h_t m_t}$ is well-defined and satisfies, for every bounded
measurable $f$,
\[
    \brac{f,h_t m_t}
    =
    \brac{f,m_0}
    +
    \int_0^t
    \int_{\RRP}
    \left(A_s f(x)+B_s f(x)\right)
    h_s m_s(\dd x)\dd s,
\]
where
\[
    A_t f(x)
    =
    \int_{\RRP}
    K_\alpha(x,y)f(x+y)G_t(x,y)h_t\mu_t(\dd y),
\]
and
\[
    B_t f(x)
    =
    -h_t^{-1}(x)
    \int_{\RRP}
    K_\alpha(x,y)f(y)h_t\mu_t(\dd y).
\]

Indeed, first perform the computation for cutoff initial data. In that case
the kernel is bounded on the support, $h_t$ is bounded on the support of the
linearized solution, and the identity follows by differentiating
$\brac{f h_t,m_t}$. The term coming from $\partial_t h_t$ cancels the loss term
in which $f$ is evaluated at the variable carried by $m_t$, leaving precisely
the two operators $A_t$ and $B_t$ above. Passing to the singular equation uses
the same cutoff approximation as in Lemma~\ref{lem:exponential-moments},
together with the well-posedness of the linearized equation from
Lemma~\ref{lem:well_posedness_derivative_eq}. More explicitly, writing $n_t^\varepsilon=h_t^\varepsilon m_t^\varepsilon$, the sign estimate for the rescaled equation and the positive exponential-moment bound give
\[
    \sup_{t\leq T}\|n_t^\varepsilon\|_{\psi_2}
    \leq\|m_0\|_{\psi_2}
    \exp\!\left(cT e^{T\|\mu\|_{\psi_1}}\|\mu\|_{\psi_3}\right).
\]
Indeed, the gain is controlled using $\psi_2(x+y)\leq\psi_2(x)\wedge\psi_2(y)$; the $B_t$ term requires the $\psi_3$-moment of $h_t\mu_t$. The displayed bound is uniform in the cutoff and gives uniform integrability in the weaker weight $\psi_1$. On $[a,\infty)$, the weights $h_t^\varepsilon$ converge uniformly and are uniformly bounded. Combined with the strong convergence of $m_t^\varepsilon$, this identifies the limit as $h_tm_t$ and justifies passage in the rescaled equation for bounded tests. The same argument with index $p$ and the moment $\psi_{p+1}$ applies where higher weights are used below. 

\section{Wasserstein regularity of solutions and derivative}\label{section:Wass-regularity}
In this section, we establish Wasserstein bounds for the solution to the $\SCE$ and for its derivative. The Wasserstein distance we consider is particularly well suited to the considered kernel $K_\alpha$. An equivalent distance was used in the works of Fournier and Cepeda~\cite{cepedaSmoluchowskisEquationRate2011a}, as well as in~\cite{fournierDistanceCoagulation2006,fournierStochasticCoalescenceHomogeneouslike2009} for kernels satisfying $K(x,y) \le  (x + y)^{-\alpha}$. It also plays a central role in our final result, as it provides a convergence rate of random empirical measures.

This section is very technical due to the singular structure of $K_\alpha$, which prevents us from obtaining Wasserstein bounds through direct methods. An illustration of this difficulty is given in Example~\ref{exp:wasserstein-is-tricky}. We first state the main results in Section~\ref{section:wass-main-results}, then introduce the rescaling strategy and its technical estimates in Section~\ref{section:wass-outline}. The argument relies on the generation of exponential moments proved in Section~\ref{section:sce-exponential}. We next control the rescaled derivative and the rescaled solutions in Sections~\ref{section:wass-rescaled-derivative} and~\ref{section:wass-rescaled-solutions}; the final subsection proves Propositions~\ref{prop:wass-derivative} and~\ref{prop:wass-final}.

\subsection{Main results}\label{section:wass-main-results}

Recall that for $\mu \in \MC^+_{\psi_2}(\RRP)$, the notation $\Proc{\varphi_t(\mu)}$ denotes the solution to the Smoluchowski coagulation equation (SCE) with initial condition $\mu$ and kernel $K_\alpha$. And for $\mu \in \MC^+_{\psi_3}(\RRP)$, we denote by $\Proc{\delta_\mu\varphi_t(\mu)}$ its flat derivative flow. The latter satisfies the following stability properties for the Wasserstein-type distance $W_{\psi_1}$.

\begin{proposition}[Wasserstein Derivative Bound]\label{prop:wass-derivative}
    Let $\mu \in \MC^+_{\psi_3}(\RRP)$ and fix $T > 0$. Then there exists a constant $C(T,\mu)$ of the form
    \[
    C(T,\mu) = ce^{cT\|\mu\|_{\psi_2}\left(1 + T\|\mu\|\right)},
    \]
    such that for all $t \in [0, T]$, the following bounds hold:
    \begin{enumerate}
        \item For all $z \in \RRP$,
        \[
        W_{\psi_1}\left( \delta_\mu \varphi_t(\mu; z), 0 \right) \leq (1 + \psi_1(z)) C(T,\mu).
        \]
        \item For all $z_1, z_2 \in \RRP$,
        \[
        W_{\psi_1}\left( \delta_\mu \varphi_t(\mu; z_1), \delta_\mu \varphi_t(\mu; z_2) \right) \leq d_{\psi_1}(z_1, z_2) C(T,\mu).
        \]
    \end{enumerate}
\end{proposition}

\begin{remark}
    Equivalently, for any $f \in  \FC_{\psi_1}$, the map $z \mapsto \brac{f, \delta_\mu \varphi_t(\mu; z)}$ is controlled in $\FC_{\psi_1}$ by $C(T,\mu)\Wnorm{f}{\psi_1}$.
\end{remark}
\begin{proof}
    The proof is given in Section~\ref{section:wass-rescaled-derivative}. Its central ingredients are the rescaling strategy and the technical estimates developed in Section~\ref{section:wass-outline}.
\end{proof}

The next result establishes the Lipschitz regularity of the SCE flow in the Wasserstein distance. Its proof is a direct consequence of Proposition~\ref{prop:wass-derivative} and the definition of the flat derivative.

\begin{proposition}[Wasserstein Smoluchowski Bound]\label{prop:wass-sce}
    Let $T > 0$ and $\mu, \tilde{\mu} \in \MC^+_{\psi_3}(\RRP)$. Then, for all $t \in [0,T]$,
    \begin{align*}
        W_{\psi_1} (\varphi_t(\mu),\varphi_t(\tilde{\mu})) \leq W_{\psi_1} (\mu,\tilde{\mu})\, C(T,\mu + \tilde{\mu}),
    \end{align*}
    where the constant $C$ is of the form
    \[
        C(T,m) = ce^{cT\|m\|_{\psi_2}\left(1 + T\|m\|\right)}.
    \]
\end{proposition}

\begin{proof}
    Let $f \in \FC_{\psi_1}$ with $\Wnorm{f}{\psi_1}\leq 1$. By the definition of the flat derivative, we have:
    \begin{align*}
        \brac{f,\varphi_t(\mu)-\varphi_t(\tilde{\mu})} = \int_0^1 \int_{\RRP} \brac{f, \delta_\mu\varphi_t(\theta\mu + (1-\theta)\tilde{\mu};z)} \left(\mu-\tilde{\mu}\right)(\dd z)\, \dd \theta.
    \end{align*}
    From Proposition~\ref{prop:wass-derivative}, the integrand is controlled in $\FC_{\psi_1}$ by $C(T,\mu + \tilde{\mu})$, hence:
    \[
    \int_{\RRP} \brac{f, \delta_\mu\varphi_t(\theta\mu + (1-\theta)\tilde{\mu};z)} \left(\mu-\tilde{\mu}\right)(\dd z) \le C(T,\mu + \tilde{\mu})W_{\psi_1}(\mu,\tilde{\mu}).
    \]
    The result follows by taking the supremum over $\Wnorm{f}{\psi_1}\leq 1$.
\end{proof}

Next, we provide a key comparison estimate between derivatives corresponding to different initial conditions. This result is central to establishing the quadratic error term that leads to the~$1/N$-rate. To proceed, we introduce a first-order version of the weighted Wasserstein distance. For $p \geq 1$, set
\begin{align*}
    \FC^1_{\psi_p}
    &:=
    \left\{f \in \FC_{\psi_p}, \quad f \text{ differentiable },\quad \sup\limits_{x \neq y} \dfrac{\left|f'(x) - f'(y) \right|}{d_{|\psi_p'|}(x,y)} < \infty \right\}, \\
    \Wnorm{f}{\psi_p,1} &:=
    \max\left\{
        \Wnorm{f}{\psi_p},
        \sup\limits_{x \neq y} \dfrac{\left|f'(x) - f'(y) \right|}{d_{|\psi_p'|}(x,y)}
    \right\}.
\end{align*}
We define
\begin{align*}
    W^1_{\psi_p}(\mu, \nu)
    :=
    \sup_{\Wnorm{f}{\psi_p,1}\leq 1} \brac{f,\mu-\nu}.
\end{align*}
Thus $W^1_{\psi_p}\leq W_{\psi_p}$ and note that $\psi_p \in \FC^1_{\psi_p}$. This stronger test class is only used for the Proposition below.

\begin{proposition}[Wasserstein Comparison of Derivatives]\label{prop:wass-final}
    Let $\mu, \tilde{\mu} \in \MC_{\psi_3}^+(\RRP)$, and fix $T > 0$. Then there exists a constant $C(T,\mu + \tilde{\mu})$ such that:
    \begin{enumerate}
        \item For all $z \in \RRP$,
        \begin{multline*}
        W^1_{\psi_1}\left( \delta_\mu \varphi_t(\mu; z), \delta_\mu \varphi_t(\tilde{\mu}; z)\right)  \\
        \leq \left( W_{\psi_1 + \psi_2}(\mu,\tilde{\mu})(1 + \psi_1(z) )  +W_{\psi_1}(\mu,\tilde{\mu})(1 + \psi_2(z)) \right) C(T,\mu + \tilde{\mu}) .
        \end{multline*}
        \item For all $z_1,z_2 \in \RRP$,
        \begin{multline*}
            W^1_{\psi_1}\left( \delta_\mu \varphi_t(\mu; z_1) - \delta_\mu \varphi_t(\mu; z_2), \delta_\mu \varphi_t(\tilde{\mu}; z_1)-\delta_\mu \varphi_t(\tilde{\mu}; z_2)\right) \\
             \leq \left( W_{\psi_1 + \psi_2}(\mu,\tilde{\mu})d_{\psi_1}(z_1,z_2)  +W_{\psi_1}(\mu,\tilde{\mu})d_{\psi_1+ \psi_2}(z_1,z_2)  \right) C(T,\mu + \tilde{\mu}) .
        \end{multline*}
    \end{enumerate}
    The constant $C$ is of the form:
    \begin{align*}
       C(T,\nu) =  c(1 + T)(1 + T\|\nu\|_{\psi_3}) e^{c T\|\nu\|_{\psi_2}(1 + T\|\nu\|)}.
    \end{align*}
\end{proposition}
\begin{proof}
    The proof is given in Section~\ref{section:wass-final} and combines all the estimates developed in this section.
\end{proof}

\subsection{Strategy and technical estimates}\label{section:wass-outline}\label{section:wass-technical}

This subsection has two purposes. First, we illustrate why the singularity of $K_\alpha$ prevents a direct Wasserstein argument. We then introduce the weighted rescaling that overcomes this obstruction and collect the regularity properties of its coefficients in a technical lemma used throughout the section.

We begin by contrasting our kernel with the simpler one treated in~\cite{cepedaSmoluchowskisEquationRate2011a}; the rescaling based on exponential moments will then provide the remedy.

\begin{example}\label{exp:wasserstein-is-tricky}
    Let $\Proc{\mu_t}, \Proc{\tilde{\mu}_t}$ be two solutions of the $\SCE$ starting from $\mu$ and $\tilde{\mu}$, respectively, and with kernel $K'_\alpha(x,y) = (x+y)^{-\alpha}$.
   
    This kernel, unlike $K_\alpha$, is only singular at the origin $(0,0)$. For $f \in \FC_{\psi_1}$, the difference of two solutions can be written in the form
    \begin{align*}
        \langle f, \mu_t - \tilde{\mu}_t \rangle
        = \langle f, \mu - \tilde{\mu} \rangle
        + \int_0^t \langle A_s f, \mu_s - \tilde{\mu}_s \rangle \dd s,
    \end{align*}
    where
    \begin{align*}
        A_s f(x) = \dfrac12\int_{\RRP} K'_\alpha(x,y) \left[f(x+y) - f(x) - f(y)\right] \left(\mu_s + \tilde{\mu}_s\right)(\dd y).
    \end{align*}
    The key estimate is immediate since $K'_{\alpha} \le x^{-\alpha} \wedge y^{-\alpha}$:
    \[
        \Wnorm{A_s f}{\psi_1}
        \lesssim \|\mu+\tilde{\mu}\|_{\psi_1}\Wnorm{f}{\psi_1}.
    \]
    Taking the supremum over $\Wnorm{f}{\psi_1}\leq 1$ and applying Gronwall's lemma gives the desired Lipschitz bound in $W_{\psi_1}$ of $\mu \mapsto \varphi_t(\mu)$. The whole difficulty below is that this one-line estimate fails for $K_\alpha$.
\end{example}

To recover analogous estimates for the singular kernel, we use the generation of exponential moments. This phenomenon may be viewed as an exponentially fast escape of mass from a neighborhood of the singularity at $0$. It motivates the weight
\[
    h_t(x) := \exp\left( \int_0^t \int_{\RRP} K_\alpha(x,y) \mu_s(\dd y) \dd s \right).
\]
Lemma~\ref{lem:exponential-moments} shows that the \emph{rescaled solution} $\Proc{h_t\mu_t}$ belongs to the same space as $\mu$ and is the unique solution of
\[
    \langle f, h_t \mu_t \rangle
    =
    \langle f, \mu \rangle
    + \frac{1}{2}\int_0^t
    \iint_{\RRP \times \RRP}
    K_\alpha(x,y) f(x+y) h_s(x+y) \mu_s(\dd x) \mu_s(\dd y)\dd s
\]
for every bounded measurable test function $f$. We similarly call $\Proc{h_t \delta_\mu \varphi_t(\mu)}$ the \emph{rescaled derivative}. The rescaling removes the obstruction identified above at the cost of more involved weight estimates.

The technical lemma below provides the common estimates used for the rescaled derivative, the rescaled solution, and the final comparison of derivatives. Since $h_0=1$, the original and rescaled solutions have the same initial state.

To derive regularity results for $\Proc{\mu_t}$, we use the identity
\[
    \mu_t = h_t^{-1} (h_t \mu_t),
\]
and observe that for any $f \in \FC_{\psi_1}$, the function $f h_t^{-1}$ also belongs to $\FC_{\psi_1}$.

To verify the Lipschitz part of the norms below, we repeatedly use the characterization
\[
f \in \operatorname{Lip}_1(d_{\psi_1})
\quad\Longleftrightarrow\quad
f \text{ is locally absolutely continuous and }
|f'(x)| \leq \alpha x^{-1-\alpha}\ \text{a.e.}
\]

Let $\Proc{\tilde{\mu}_t}$ be the solution of the Smoluchowski coagulation equation with initial state $\tilde\mu$, and define $\tilde h_t$ analogously to $h_t$. We will also use the weight ratio
\begin{align*}
        G_{t}(x,y) &= h_t(x+y)h_t^{-1}(x)h_t^{-1}(y)\\
         &= e^{\int_0^t \int_{\RRP} \left(K_\alpha(x+y,z) - K_\alpha(x,z) -K_\alpha(y,z)\right)\mu_s(\dd z) \dd s}
\end{align*}
and define $\tilde G_t$ analogously. The next lemma collects the regularity estimates needed for the Wasserstein arguments. Throughout, $h_t^{-1}$ and $\tilde h_t^{-1}$ denote the pointwise inverse weights,
\[
h_t^{-1}(x)=\frac{1}{h_t(x)},
\qquad
\tilde h_t^{-1}(x)=\frac{1}{\tilde h_t(x)},
\]
and not inverse functions.

\begin{lemma}\label{lem:technical_G_h}
    Let $p \ge 1$, the following estimates hold uniformly on bounded time intervals.
    \begin{enumerate}[label=(\roman*)]
        \item The inverse weights satisfy for all $f \in \FC_{\psi_p}$ 
        \[
            \Wnorm{f h_t^{-1}}{\psi_1 + \psi_p} \lesssim (1+t\|\mu\| )\Wnorm{f}{\psi_p},
        \]
        and for all $f \in \FC_{\psi_1}$:
        \[
            \Wnorm{f(h_t^{-1}-\tilde h_t^{-1})}{\psi_1+\psi_2 }
            \lesssim (1+t\|\mu+\tilde\mu\|)\int_0^t W_{\psi_1}(\mu_s,\tilde\mu_s)\dd s \Wnorm{f}{\psi_1}.
        \]
        \item For every fixed $y\in\RRP$, and $f \in \FC_{\psi_p}$
        \[
            \Wnorm{f(\cdot + y)K_\alpha(\cdot,y)G_t(\cdot,y)}{\psi_1}
            \lesssim (1 + \psi_p(y))e^{-\int_0^t \brac{\psi_1,\mu_s}\dd s}\Wnorm{f}{\psi_p},
        \]
        \item For every fixed $y\in\RRP$, and $f \in \FC_{\psi_p}$
        \begin{multline*}
            \Wnorm{f(\cdot + y)K_\alpha(\cdot,y)\left(G_t(\cdot,y)-\tilde G_t(\cdot,y)\right)}{\psi_1+\psi_2}\\
            \lesssim (1 + t\|\mu+\tilde\mu\|) (1+\psi_2(y) + \psi_p(y) )\int_0^t W_{\psi_1}(\mu_s,\tilde\mu_s)\dd s \Wnorm{f}{\psi_p}.
        \end{multline*}
    \end{enumerate}
    The same bounds hold after exchanging the variables $x$ and $y$.
\end{lemma}

\begin{proof}[Proof of Lemma~\ref{lem:technical_G_h}]
    Write
    \[
    \gamma_t = \int_0^t \brac{1,\mu_s}\dd s,
    \qquad
    \tilde\gamma_t = \int_0^t \brac{1,\tilde\mu_s}\dd s.
    \]
    Then
    \[
    h_t^{-1}(x)= e^{-\int_0^t \brac{K_\alpha(x,\cdot), \mu_s} \dd s} =  e^{-x^{-\alpha}\gamma_t}e^{-\int_0^t\brac{\psi_1,\mu_s}\dd s}.
    \]
    Recall that solutions to the $SCE$ are positive and therefore
    \[
    \gamma_t \ge 0\ ,\qquad h_t \ge 1.
    \]
    Weights depending on $\Proc{\tilde{\mu}_t}$ are defined similarly.

    $(i)$ Let $f \in \FC_{\psi_p}$, since $h_t^{-1} \le e^{-\int_0^t\brac{\psi_1,\mu_s}\dd s}$ we get the bound 
    \[
    \left| f h_t^{-1} \right| \le (1 + \psi_p) e^{-\int_0^t\brac{\psi_1,\mu_s}\dd s}\Wnorm{f}{\psi_p}.
    \]
    A direct differentiation gives
    \begin{align*}
        \dd_x(fh_t^{-1})(x) &= \alpha x^{-1-\alpha}f(x)\gamma_t h_t^{-1}(x) + f'(x)h_t^{-1}(x)  \\
        &= -\psi_1'(x) f(x)\gamma_t h_t^{-1}(x) + f'(x)h_t^{-1}(x).
    \end{align*}
    We have:
    \begin{align*}
        |\dd_x(fh_t^{-1})(x)|
        &\le (\alpha x^{-1-\alpha}\gamma_t(1 + x^{-p\alpha}) + p\alpha x^{-p\alpha-1})e^{-x^{-\alpha}\gamma_t}e^{-\int_0^t\brac{\psi_1,\mu_s}\dd s}\Wnorm{f}{\psi_p} \\
        &\lesssim \left[ \gamma_t x^{-1-\alpha} +  x^{-p\alpha-1}\left(1 + \gamma_t x^{-\alpha}e^{-\gamma_t x^{-\alpha}} \right) \right] e^{-\int_0^t\brac{\psi_1,\mu_s}\dd s}\Wnorm{f}{\psi_p} \\
        &\lesssim \left(|\left(\psi_1 + \psi_p\right)'(x)|\right) (1 + \gamma_t)\left(1 + \gamma_t x^{-\alpha}e^{-\gamma_t x^{-\alpha}} \right) e^{-\int_0^t\brac{\psi_1,\mu_s}\dd s}\Wnorm{f}{\psi_p}.
    \end{align*}
    Lemma~\ref{lem:sub-add-bound-smol} gives $\gamma_t\leq t\|\mu\|$, from $ue^{-u}\le 1, u \ge 0$ we get the first announced bound:
    \[
    \Wnorm{f h_t^{-1}}{\psi_1 + \psi_p}\lesssim \left( 1+t\|\mu\| \right) e^{-\int_0^t\brac{\psi_1,\mu_s}\dd s}\Wnorm{f}{\psi_p}.
    \]
    Let $f \in \FC_{\psi_1}$, for the difference of two inverse weights, the inequality
    $|e^{-a}-e^{-b}|\leq |a-b|$ for $a,b > 0$ gives
    \[
        \left|h_t^{-1}(x)-\tilde h_t^{-1}(x) \right| \le \int_0^t \left|\int_{\RRP} K_\alpha(x,y) (\mu_s - \tilde{\mu}_s)(\dd y) \right|\dd s 
    \]
    The function $y \mapsto K_\alpha(x,y)$ has $\FC_{\psi_1}$-norm controlled by $1 + \psi_1(x)$ giving:
    \begin{equation}\label{proof:technical:diff_ht}
    \left|h_t^{-1}(x)-\tilde h_t^{-1}(x) \right| \lesssim (1 + \psi_1(x))\int_0^t W_{\psi_1}(\mu_s,\tilde{\mu}_s) \dd s.
    \end{equation}
    Multiplication by $f$ finally gives the bound:
    \[
    |f(x)|\left|h_t^{-1}(x)-\tilde h_t^{-1}(x) \right|
    \lesssim (1 + \psi_2(x))\int_0^t W_{\psi_1}(\mu_s,\tilde{\mu}_s) \dd s\,\Wnorm{f}{\psi_1}.
    \]
    Differentiating gives
    \begin{equation}\label{eq:proof:formdfht}
    \dd_x\left(f(h_t^{-1}-\tilde h_t^{-1})\right)(x)
    = f'(x)\left(h_t^{-1}(x)-\tilde h_t^{-1}(x)\right)
    + f(x)\dd_x\left(h_t^{-1}-\tilde h_t^{-1}\right)(x).
    \end{equation}
    We already treated the first term on the RHS, it remains to bound the weights-derivative difference:
    \[
    \dd_x\left(h_t^{-1}-\tilde h_t^{-1}\right)(x) = \alpha x^{-1-\alpha} \left( \gamma_t h_t^{-1}(x) - \tilde \gamma_t \tilde h_t^{-1}(x)\right)
    \]
    Here we rely on the following inequality, for $a,b \ge 0$ and $c \ge d \ge 0$:
    \[
    \left| a e^{-c} - b e^{-d} \right| \le |a -b | + b e^{- d} |c - d |.
    \]
    Assuming with no loss of generality that for fixed $t \ge 0$ and $x \in \RRP$:
    \[
    \int_0^t \brac{K_\alpha(x,\cdot),\tilde{\mu}_s}  \dd s \ge \int_0^t \brac{K_\alpha(x,\cdot),\mu_s}  \dd s
    \]
    we get:
    \[
    \left| \dd_x\left(h_t^{-1}-\tilde h_t^{-1}\right)(x)\right| \le \alpha x^{-1-\alpha}\left| \gamma_t - \tilde{\gamma}_t\right| + \alpha x^{-1-\alpha}\gamma_t h_t^{-1}(x) \int_0^t \big| \brac{K_\alpha(x,\cdot),\mu_s - \tilde \mu_s }\big|\dd s
    \]
    Using $\left| \gamma_t - \tilde{\gamma}_t\right| \le \int_0^t W_{\psi_1}(\mu_s,\tilde \mu_s)\dd s$ we obtain:
    \[
    \left| \dd_x\left(h_t^{-1}-\tilde h_t^{-1}\right)(x)\right| \lesssim x^{-1-\alpha} \left(1 + \gamma_t h_t^{-1}(x)(1 + x^{-\alpha}) \right) \int_0^t W_{\psi_1}(\mu_s,\tilde \mu_s)\dd s.
    \]
    Using the damping provided by $\gamma_t x^{-\alpha} e^{-\gamma_t x^{-\alpha}}$ we end with
    \begin{equation*}
        \left| \dd_x\left(h_t^{-1}-\tilde h_t^{-1}\right)(x)\right| \lesssim x^{-1-\alpha} \left(1 + t\|\mu+\tilde\mu\|\right) \int_0^t W_{\psi_1}(\mu_s,\tilde \mu_s)\dd s.
    \end{equation*}
    By using the bound above with inequality \eqref{proof:technical:diff_ht} into the product formula \eqref{eq:proof:formdfht} we finally find:
    \[
    \left| \dd_x\left(f(h_t^{-1}-\tilde h_t^{-1})\right)(x) \right|
    \lesssim \left(x^{-1- \alpha} + x^{-1- 2\alpha} \right)
    \left(1 + t\|\mu+\tilde\mu\|\right)\int_0^t W_{\psi_1}(\mu_s,\tilde \mu_s)\dd s\,\Wnorm{f}{\psi_1}.
    \]
    Together with the growth bound above, this proves
    \[
    \Wnorm{f(h_t^{-1}-\tilde h_t^{-1})}{\psi_1+\psi_2}
    \lesssim \left(1+t\|\mu+\tilde\mu\|\right)
    \int_0^t W_{\psi_1}(\mu_s,\tilde\mu_s)\dd s\,\Wnorm{f}{\psi_1}.
    \]

    Assertions~$(ii)$--$(iii)$ are proved in Appendix~\ref{appendix:proof-technical-G-h}. They rely on the same mechanism: the non-positive exponent in $G_t$ provides a damping factor that compensates for the singular nature of $K_\alpha$.

\end{proof}

\subsection{Control of the rescaled derivative}\label{section:wass-rescaled-derivative}

\begin{proposition}[Wasserstein bounds for the rescaled derivative]\label{lem:wass-rescaled-derivative}
    Let $T > 0$ be a fixed horizon of time. The following assertions hold for all $t \in [0,T]$.
    \begin{enumerate}[label=(\roman*)]
        \item Assume that $\mu \in \MC^+_{\psi_3}(\RRP)$. For all $f\in\FC_{\psi_1}$ and all $z \in \RRP$:
        \begin{align*}
            \brac{f,h_t\delta_\mu \varphi_t(\mu;z)}
            &\leq C_1(T,\mu)(1 + \psi_1(z))\Wnorm{f}{\psi_1},
        \end{align*}
        where
        \[
        C_1(T,\mu) = e^{cT\|\mu\|_{\psi_2}\left(1 + T\|\mu\|\right)}.
        \]
        \item For all $f\in\FC_{\psi_1}$ and all $z_1,z_2 \in \RRP$:
            \begin{align*}
                \brac{f,h_t\delta_\mu \varphi_t(\mu;z_1) - h_t\delta_\mu \varphi_t(\mu;z_2)}
                \leq C_1(T,\mu)d_{\psi_1}(z_1,z_2)\Wnorm{f}{\psi_1} .
            \end{align*}
        \item For all $f\in\FC_{\psi_p}$ and all $z \in \RRP$:
            \begin{align*}
                \brac{f,h_t\delta_\mu \varphi_t(\mu;z)} \leq \left(1 + \psi_p(z) + C_2(T,\mu )(1 + \psi_1(z))\right)\Wnorm{f}{\psi_p},
            \end{align*}
             where 
            \begin{align*}
                C_2(T,\mu ) = cT \|\mu\|_{\psi_{p+1}}e^{cT\|\mu\|_{\psi_{2}}\left(1 + T\|\mu\|\right)}.
            \end{align*}
        \item For all $f\in\FC_{\psi_p}$ and all $z_1,z_2 \in \RRP$:
            \begin{align*}
                \brac{f,h_t\delta_\mu \varphi_t(\mu;z_1) - h_t\delta_\mu \varphi_t(\mu;z_2)}
                \leq \left(d_{\psi_p}(z_1,z_2) + C_2(T,\mu ) d_{\psi_1}(z_1,z_2)\right)\Wnorm{f}{\psi_p}.
            \end{align*}
    \end{enumerate}
\end{proposition}

\begin{proof}
    We denote the solution to the SCE $\Proc{\mu_t} = \Proc{\varphi_t(\mu)}$. We prove $(i)$ and $(ii)$ together leveraging the linearity of the derivative's equation. We work with a process $\Proc{m_t}$, the unique solution of:
    \begin{align*}
        \brac{f,m_t} = \brac{f,m_0} + \int_0^t \brac{K_\alpha f,\mu_s\otimes m_s} \dd s.
    \end{align*}
    Taking either $m_0 = \delta_z$ or $\delta_{z_1} - \delta_{z_2}$ yields either $(i)$ or $(ii)$ by uniqueness. The rescaled equation tested against $f$ is:
    \[
        \brac{f,h_tm_t} = \brac{f,m_0} + \int_0^t \int_{\RRP} \left(A_sf(x) + B_sf(x)\right) h_sm_s(\dd x) \dd s.
    \]
    with
    \begin{subequations}
    \begin{align}
        A_tf(x) &= \int_{\RRP} K_\alpha(x,y)f(x+y)G_t(x,y)  h_t\mu_t(\dd y) \label{proof:eq:wass-rescaled-SCE-1}\\ 
        B_tf(x) &= -\int_{\RRP}K_\alpha(x,y)h_t^{-1}(x) f(y) h_t\mu_t(\dd y)  \label{proof:eq:wass-rescaled-SCE-2}.
    \end{align}
    \label{proof:eq:wass-rescaled-SCE}
    \end{subequations}
    Recall the notation $\gamma_t = \int_0^t \brac{1 ,\mu_s} \dd s$. For \eqref{proof:eq:wass-rescaled-SCE-1} we have:
    \begin{align*}
    \Wnorm{A_tf}{\psi_1}
    &\le \int_{\RRP} \Wnorm{f(y+\cdot)K_\alpha(y,\cdot)G_t(y,\cdot)}{\psi_1}  h_t\mu_t(\dd y)  \\
    &\lesssim \int_{\RRP} (1 + \psi_1(y)) e^{ y^{-\alpha} \gamma_t} \mu_t(\dd y)\,\Wnorm{f}{\psi_1} \\
    &\le \|\mu\|_{\psi_1}\Wnorm{f}{\psi_1}.
    \end{align*}
    where we used Lemma \ref{lem:technical_G_h}-$(ii)$ for the second inequality and Proposition~\ref{prop:exponential-moments} for the last.

    For \eqref{proof:eq:wass-rescaled-SCE-2}, since $\Wnorm{K_\alpha(\cdot,y)}{\psi_1} \le 1 + \psi_1(y)$ Lemma \ref{lem:technical_G_h}-$(i)$ and Proposition~\ref{prop:exponential-moments} give
    \begin{align*}
        \Wnorm{B_tf}{\psi_1}
        &\le \int_{\RRP}\Wnorm{K_\alpha(\cdot,y)h_t^{-1}}{\psi_1} |f(y)| h_t\mu_t(\dd y) \\
        &\lesssim (1 + t\|\mu\|)\int_{\RRP}(1 + \psi_1(y)) |f(y)| e^{-\int_0^t \brac{\psi_1,\mu_s} \dd s}h_t\mu_t(\dd y)\\
        &\le (1 + t\|\mu\|)\|\mu\|_{\psi_2} \Wnorm{f}{\psi_1}.
    \end{align*}
   
    By a slight abuse of notation, write $W_{\psi_1}(h_tm_t) := W_{\psi_1}(h_tm_t,0)$. We have shown that
    \[
    \brac{f,h_tm_t} \le W_{\psi_1}(m_0)\Wnorm{f}{\psi_1}
    + c(1 + t\|\mu\|)\|\mu\|_{\psi_2} \int_0^t W_{\psi_1}(h_sm_s) \dd s\,\Wnorm{f}{\psi_1}.
    \]
    Taking the supremum over $\Wnorm{f}{\psi_1}\le1$ and using Gronwall's lemma gives
    \[
    W_{\psi_1}(h_tm_t)
    \le
    W_{\psi_1}(m_0)
    e^{cT\|\mu\|_{\psi_2}(1+T\|\mu\|)}.
    \]
    Finally,
    \[
    W_{\psi_1}(\delta_z,0)=1+\psi_1(z),
    \qquad
    W_{\psi_1}(\delta_{z_1}-\delta_{z_2},0)=d_{\psi_1}(z_1,z_2),
    \]
    which gives $(i)$ and $(ii)$.

    For $(iii)$ and $(iv)$ we follow the same line of computations. Take $f \in \FC_{\psi_p}$ we find
    \[
    \Wnorm{A_tf}{\psi_1} \lesssim \|\mu\|_{\psi_p}\Wnorm{f}{\psi_p}\qquad\  \text{and}\qquad\ \Wnorm{B_tf}{\psi_1} \lesssim (1 + t\|\mu\|)\|\mu\|_{\psi_{p + 1}} \Wnorm{f}{\psi_p}.
    \]
    This yields the estimate after taking the supremum over $\Wnorm{f}{\psi_p} \le 1$:
    \[
    W_{\psi_p}(h_tm_t) \le W_{\psi_p}(m_0) + c(1 + t\|\mu\|)\|\mu\|_{\psi_{p + 1}}\int_0^t W_{\psi_1}(h_s m_s) \dd s. 
    \]
    Injecting the bound found previously for $(i)$ and $(ii)$ we find for all $t \in [0,T]$:
    \[
    W_{\psi_p}(h_tm_t) \le W_{\psi_p}(m_0) + W_{\psi_1}(m_0)cT(1 + T\|\mu\|)\|\mu\|_{\psi_{p + 1}}e^{cT\|\mu\|_{\psi_2}(1+T\|\mu\|)}.
    \]
    To obtain the final constant form observe that $e^{cT\|\mu\|_{\psi_2}} = e^{cT\brac{1 + \psi_2, \mu}} = e^{cT(\|\mu\| + \brac{\psi_2,\mu})}$ therefore we can absorb the factor $T\|\mu\|$ into the exponential.
\end{proof}

\begin{proof}[Proof of Proposition \ref{prop:wass-derivative}]
    Using the same linearity reduction as above, let $m_t$ solve the derivative equation with $m_0=\delta_z$ or $m_0=\delta_{z_1}-\delta_{z_2}$. For $f\in\FC_{\psi_1}$, Lemma~\ref{lem:technical_G_h}--$(i)$ and Proposition~\ref{lem:wass-rescaled-derivative} give
    \begin{align*}
        |\brac{f,m_t}|
        &=|\brac{fh_t^{-1},h_tm_t}| \\
        &\leq c(1+T\|\mu\|)
        e^{cT\|\mu\|_{\psi_2}(1+T\|\mu\|)}
        W_{\psi_1}(m_0)\Wnorm{f}{\psi_1}.
    \end{align*}
    Taking the supremum over $\Wnorm{f}{\psi_1}\leq1$ and using the identities for $W_{\psi_1}(m_0)$ from the preceding proof yields the result after absorbing the linear term into the exponential.
\end{proof}

\subsection{Control of rescaled solutions}\label{section:wass-rescaled-solutions}

\begin{proposition}[Wasserstein bounds for rescaled solutions]\label{prop:wass-rescaled-SCE}
    Let $T>0$, and let $\Proc{\mu_t}$ and $\Proc{\tilde\mu_t}$ be the solutions to the $\SCE$ with initial conditions $\mu,\tilde\mu\in\MC^+_{\psi_3}(\RRP)$. The following statements hold:
    
    \begin{enumerate}[label=(\roman*)]
    \item For every $t\in[0,T]$:
    \[
            W_{\psi_1}\left(h_t\mu_t,\tilde h_t\tilde\mu_t\right)\leq C_1(T,\mu + \tilde{\mu}) W_{\psi_1}(\mu,\tilde\mu),
        \]
    where
    \[
     C_1(T,m) = e^{
            c_1T\|m\|_{\psi_2}(1+T\|m\|)
        }.
    \]
    \item For every $t\in[0,T]$:
     \[
            W_{\psi_2}\left(h_t\mu_t,\tilde h_t\tilde\mu_t\right)
            \leq
            W_{\psi_2}(\mu,\tilde\mu)
            +C_2(T,m)
            W_{\psi_1}(\mu,\tilde\mu),
    \]
        where
    \[
         C_2(T,m) =  e^{
            c_2T\|m\|_{\psi_2}(1+T\|m\|)
        }.
     \]
     \end{enumerate}
\end{proposition}

This proposition can be interpreted as follows. The first point shows that the rescaled solutions of the $\SCE$ are Lipschitz with respect to the distance $W_{\psi_1}$. The second point shows that the dynamic variation of the distance $W_{\psi_2}$ is controlled by the initial $W_{\psi_2}$ distance and by the propagated $W_{\psi_1}$ estimate. This is another argument in favor of using $W_{\psi_1}$ as the main distance in our problem.

\begin{proof}
    We first prove $(i)$. Set $m:=\mu+\tilde\mu$. Let $f\in\FC_{\psi_1}$ and introduce the functions
    \begin{align*}
        H_tf(x,y)
        &:=
        K_\alpha(x,y)f(x+y)G_t(x,y),
        \\
        \tilde H_tf(x,y)
        &:=
        K_\alpha(x,y)f(x+y)\tilde G_t(x,y).
    \end{align*}
    The difference between the two rescaled solutions satisfies:
    \begin{align*}
        \brac{f,h_t\mu_t-\tilde h_t\tilde\mu_t}
        =
        \brac{f,\mu-\tilde\mu}
        +
        \frac12\int_0^t
        \left(
        \brac{H_sf,(h_s\mu_s)^{\otimes2}}
        -
        \brac{\tilde H_sf,(\tilde h_s\tilde\mu_s)^{\otimes2}}
        \right)\dd s.
    \end{align*}
    We decompose the integrand as
    \begin{equation}
    \begin{aligned}
            \brac{H_tf,(h_t\mu_t)^{\otimes2}}
            -
            \brac{\tilde H_tf,(\tilde h_t\tilde\mu_t)^{\otimes2}}
            =&
            \brac{H_tf,h_t\mu_t\otimes(h_t\mu_t-\tilde h_t\tilde\mu_t)}
            \\
            &+
            \brac{\tilde H_tf,(h_t\mu_t-\tilde h_t\tilde\mu_t)\otimes\tilde h_t\tilde\mu_t}
            \\
            &+
            \brac{H_tf-\tilde H_tf,h_t\mu_t\otimes\tilde h_t\tilde\mu_t}.
    \end{aligned}
    \label{proof:eq:h_tmu_t}
    \end{equation}
    For the first term, write
    \[
        A_tf(x)
        :=\int_{\RRP}H_tf(x,y)h_t\mu_t(\dd y).
    \]
    The function $A_tf$ is exactly the one treated in the preceding proof. Hence Lemma~\ref{lem:technical_G_h}--$(ii)$ and Proposition~\ref{prop:exponential-moments} give
    \[
        \Wnorm{A_tf}{\psi_1}
        \lesssim\|\mu\|_{\psi_1}\Wnorm{f}{\psi_1}.
    \]
    Consequently,
    \begin{align*}
         \brac{H_tf,h_t\mu_t\otimes(h_t\mu_t-\tilde h_t\tilde\mu_t)} 
        &= \brac{A_tf, h_t\mu_t-\tilde h_t\tilde\mu_t } \\
        &\lesssim \|\mu\|_{\psi_1}
        W_{\psi_1}\left(h_t\mu_t,\tilde h_t\tilde\mu_t\right)
        \Wnorm{f}{\psi_1}.
    \end{align*}
    By symmetry the second term satisfies:
    \[
    \brac{\tilde H_tf,(h_t\mu_t-\tilde h_t\tilde\mu_t)\otimes\tilde h_t\tilde\mu_t}
        \lesssim \|\tilde{\mu}\|_{\psi_1}
        W_{\psi_1}\left(h_t\mu_t,\tilde h_t\tilde\mu_t\right)
        \Wnorm{f}{\psi_1}.
    \]
    For the third term, Lemma~\ref{lem:technical_G_h}--$(iii)$ and Proposition~\ref{prop:exponential-moments} give
    \[
        \left|\brac{H_tf-\tilde H_tf,h_t\mu_t\otimes\tilde h_t\tilde\mu_t}\right|
        \lesssim
        (1+T\|m\|)\|m\|_{\psi_2}^2 e^{T\|m\|_{\psi_1}}\int_0^tW_{\psi_1}(\mu_s,\tilde\mu_s)\dd s\,
        \Wnorm{f}{\psi_1}.
    \]
    Proposition~\ref{prop:wass-sce} then controls the time integral. Hence, for $t\in[0,T]$,
    \[
        \left|\brac{H_tf-\tilde H_tf,h_t\mu_t\otimes\tilde h_t\tilde\mu_t}\right|
        \leq
        cT(1+T\|m\|)^2\|m\|_{\psi_2}^2
        e^{cT\|m\|_{\psi_2}(1+T\|m\|)}
        W_{\psi_1}(\mu,\tilde\mu)\Wnorm{f}{\psi_1}.
    \]
    Combining these estimates and taking the supremum over $\Wnorm{f}{\psi_1}\leq1$ yields
    \begin{multline*}
        W_{\psi_1}\left(h_t\mu_t,\tilde h_t\tilde\mu_t\right)
        \leq
        \left(
        1+c T^2 (1+T\|m\|)^2\|m\|_{\psi_2}^2
        e^{cT\|m\|_{\psi_2}(1+T\|m\|)}
        \right)
        W_{\psi_1}(\mu,\tilde\mu)
        \\
        +c\|m\|_{\psi_1}
        \int_0^tW_{\psi_1}\left(h_s\mu_s,\tilde h_s\tilde\mu_s\right)\dd s.
    \end{multline*}
    Gronwall's lemma proves $(i)$ after absorbing the extra polynomial terms into the exponential.

    We now prove $(ii)$. Let $f\in\FC_{\psi_2}$. The decomposition~\eqref{proof:eq:h_tmu_t} is unchanged, and the same function $A_tf$ occurs. Lemma~\ref{lem:technical_G_h}--$(ii)$ gives
    \[
        \Wnorm{A_tf}{\psi_1}
        \lesssim\|\mu\|_{\psi_2}\Wnorm{f}{\psi_2}.
    \]
    Therefore,
    \[
         \brac{H_tf,h_t\mu_t\otimes(h_t\mu_t-\tilde h_t\tilde\mu_t)}  \lesssim \|\mu\|_{\psi_2}
        W_{\psi_1}\left(h_t\mu_t,\tilde h_t\tilde\mu_t\right)
        \Wnorm{f}{\psi_2}.
    \]
    And by symmetry
    \[
         \brac{\tilde H_tf,(h_t\mu_t-\tilde h_t\tilde\mu_t)\otimes\tilde h_t\tilde\mu_t}
        \lesssim \|\tilde{\mu}\|_{\psi_2}
        W_{\psi_1}\left(h_t\mu_t,\tilde h_t\tilde\mu_t\right)
        \Wnorm{f}{\psi_2}.
    \]
    Moreover, Lemma~\ref{lem:technical_G_h}--$(iii)$ and the same moment estimates give
    \[
        \left|\brac{H_tf-\tilde H_tf,h_t\mu_t\otimes\tilde h_t\tilde\mu_t}\right|
        \leq
        cT (1+T\|m\|)^2 \|m\|_{\psi_2}^2 e^{cT\|m\|_{\psi_2}(1+T\|m\|)}
        W_{\psi_1}(\mu,\tilde\mu)\Wnorm{f}{\psi_2}.
    \]
    Inserting these bounds in the rescaled equation and taking the supremum over $\Wnorm{f}{\psi_2}\leq1$ yields
    \begin{multline*}
        W_{\psi_2}\left(h_t\mu_t,\tilde h_t\tilde\mu_t\right)
        \leq
        W_{\psi_2}(\mu,\tilde\mu) + 
        C_2(T,\mu + \tilde{\mu})
        W_{\psi_1}(\mu,\tilde\mu)
        \\
        +c\|m\|_{\psi_{2}}
        \int_0^tW_{\psi_1}\left(h_s\mu_s,\tilde h_s\tilde\mu_s\right)\dd s.
    \end{multline*}
    Inserting $(i)$ in the last integral gives $(ii)$ after absorbing the extra polynomial terms into the exponential factor.
\end{proof}

\subsection{Proof of Proposition \ref{prop:wass-final} }\label{section:wass-final}
We start by adding technical estimates for $h_t$ and $G_t$ in the stronger norm $\Wnorm{\cdot}{\psi_1,1}$. The new point is the mixed derivative estimate in $(iii)$. It is needed when a two-variable function is paired against a product of signed measures: after using a Wasserstein estimate in one variable, the resulting function of the other variable must still have the required Lipschitz regularity. This is the only place in the argument where the second-order control in $\FC_{\psi_1}^{1}$ is used.

\begin{lemma}\label{lem:technical_G_h_dos}
    The following estimates hold uniformly on bounded time intervals.
    \begin{enumerate}[label=(\roman*)]
        \item The inverse weights satisfy for all $f \in \FC_{\psi_1}^1$ 
        \[
            \Wnorm{f h_t^{-1}}{\psi_1,1} \lesssim (1+t\|\mu\| ) \Wnorm{f}{\psi_1}^1.
        \]

        \item For every fixed $y\in\RRP$, and $f \in \FC_{\psi_1}^1$
        \[
            \Wnorm{f(\cdot + y)K_\alpha(\cdot,y)G_t(\cdot,y)}{\psi_1,1}
            \lesssim (1 + \psi_1(y))e^{-\int_0^t \brac{\psi_1,\mu_s}\dd s}\Wnorm{f}{\psi_1,1}.
        \]

        \item For every $x, y\in\RRP$ and $f \in \FC_{\psi_1}^{1}$,
        \[
            |\partial_{xy}\left(f(\cdot + \cdot)K_\alpha G_t\right)(x,y)|
            \lesssim (1+t\|\mu\|)y^{-1-\alpha} \left( x^{-1-\alpha} + x^{-1 - 2\alpha} \right)\Wnorm{f}{\psi_1}^1.
        \]
    \end{enumerate}
\end{lemma}

\begin{proof}
    Second derivatives below are understood almost everywhere; $f'$ is locally Lipschitz. We reuse the notations $\gamma_t$ and $\mathbf \Lambda_t$:
    \[
    \gamma_t = \int_0^t \brac{1,\mu_s} \dd s\ ,\qquad \mathbf{\Lambda_t} = \gamma_t \left( x^{-\alpha} + y^{-\alpha} - (x + y)^{-\alpha} \right).
    \]

    $(i)$ The direct and first derivative bound are already treated in Lemma~\ref{lem:technical_G_h}'s proof. A second differentiation gives
    \[
    \dd_{xx} fh_t^{-1}(x) = f''(x)h_t^{-1}(x) + 2f'(x)\dd_x h_t^{-1}(x) + f(x)\dd_{xx}h_t^{-1}(x).
    \]
    We have:
    \begin{align*}
        \left| \dd_{xx} fh_t^{-1}(x)\right| \lesssim x^{-2 -\alpha}\left(1+\gamma_t+\gamma_t x^{-\alpha}+\gamma_t^2 x^{-\alpha}+\gamma_t^2 x^{-2\alpha}\right)h_t^{-1}(x)\Wnorm{f}{\psi_1,1}.
    \end{align*}
    Using the damping provided by $h_t^{-1}$ yields the first bound.

    $(ii)$ A second differentiation of the weight $K_\alpha G_t$ gives,
    \[
    \partial_{xx} K_\alpha G_t(x,y) = \partial_{xx} K_\alpha(x,y) G_t(x,y) + 2 \partial_x K_\alpha(x,y) \partial_x G_t(x,y) + K_\alpha(x,y) \partial_{xx} G_t(x,y).
    \]
    We have:
    \begin{align*}
        \left| \partial_{xx} K_\alpha G_t(x,y)\right| 
        &\lesssim x^{-2 - \alpha} \left(1 + \gamma_t x^{-\alpha} + \gamma_t K_\alpha(x,y) \left(1 + \gamma_t x^{-\alpha} \right) \right) e^{- \mathbf \Lambda_t}\\
        &\lesssim x^{-2 - \alpha}.
    \end{align*}
    For the product it gives,
    \begin{multline*}
        \partial_{xx} f(\cdot + y) K_\alpha G_t(x,y) = f''(x+y) K_\alpha G_t(x,y)\\
         + 2 f'(x+y) \partial_x K_\alpha G_t(x,y) + f(x+y) \partial_{xx} K_\alpha G_t(x,y).
    \end{multline*}
    Using the second derivative control above as well as the first derivative control and the bound previously established in the proof of Lemma~\ref{lem:technical_G_h}-$(ii)$, we get:
    \begin{align*}
        \left|\partial_{xx} f(\cdot + y) K_\alpha G_t(x,y) \right| \lesssim x^{-2 - \alpha} (1 + \psi_1(y))e^{-\int_0^t\brac{\psi_1,\mu_s}\dd s}\Wnorm{f}{\psi_1,1}.
    \end{align*}

    $(iii)$ Differentiating $\partial_x(K_\alpha G_t)$ with respect to $y$ gives
\begin{multline*}
    \partial_{xy}\left(K_\alpha G_t\right)(x,y) = \left[\vphantom{{A^2}^2}-\partial_y \mathbf \Lambda_t (x,y)\left(\vphantom{{A^2}} \partial_x K_\alpha(x,y) - K_\alpha(x,y )\partial_x \mathbf \Lambda_t(x,y) \right) \right.\\
    -\left. \vphantom{{A^2}^2} \partial_y K_\alpha(x,y)\partial_x \mathbf \Lambda_t(x,y) - K_\alpha(x,y)\partial_{xy} \mathbf \Lambda_t(x,y)\right] e^{-\mathbf \Lambda_t(x,y)-\int_0^t\brac{\psi_1,\mu_s}\dd s}.
\end{multline*}
In addition,
\[
|\partial_y\mathbf \Lambda_t(x,y)|\lesssim \gamma_t y^{-1-\alpha},
\qquad
|\partial_{xy}\mathbf \Lambda_t(x,y)|
=\alpha(\alpha+1)\gamma_t(x+y)^{-\alpha-2}.
\]
Using $K_\alpha(x,y)(x+y)^{-\alpha-2}\leq2(xy)^{-\alpha-1}$, we find:
\[
\left| \partial_{xy}\left(K_\alpha G_t\right)(x,y) \right|
 \lesssim \gamma_t (x y)^{-1 - \alpha}\left(1 + \gamma_t K_\alpha(x,y)\right) e^{-\mathbf \Lambda_t(x,y)-\int_0^t\brac{\psi_1,\mu_s}\dd s}.
\]
The damping provided by $ e ^{- \mathbf \Lambda_t(x,y)} \le e ^{- c_\alpha \gamma_t K_\alpha(x,y ) }$ gives
\[
\left| \partial_{xy}\left(K_\alpha G_t\right)(x,y) \right|
 \lesssim (1 + t\|\mu\|) (xy)^{-1 - \alpha} .
\]

We now treat the complete bound:
\begin{multline*}
    \partial_{xy} f(\cdot + \cdot)K_\alpha G_t(x,y) =  f'(x+y) \left( \partial_x K_\alpha G_t(x,y) + \partial_y K_\alpha G_t(x,y)\right)\\
    + f''(x+ y) K_\alpha G_t(x,y)+ f(x+y) \partial_{xy} K_\alpha G_t(x,y).
\end{multline*}
Since $f \in \FC_{\psi_1}^{1}$ the first and second term can be easily controlled thanks to the estimate for $\partial_x(K_\alpha G_t)$ obtained in~$(ii)$ by 
\[
(xy)^{-1-\alpha}\Wnorm{f}{\psi_1,1}.
\]
In the third term distributing $\psi_1(x + y)$ forces us to add an extra exponent $-\alpha$ that we arbitrarily give to the $x$-variable:
\[
    \left| f(x+y) \partial_{xy} K_\alpha G_t(x,y)\right| \lesssim
    (1+t\|\mu\|)y^{-1-\alpha} \left( x^{-1-\alpha} + x^{-1 - 2\alpha} \right)\Wnorm{f}{\psi_1}.
\]

\end{proof}

\begin{proof}[Proof of Proposition~\ref{prop:wass-final}]
    We prove the two bounds simultaneously. Let $m_t$ and $\tilde m_t$ be the solutions of the derivative equation driven respectively by $\mu_t:=\varphi_t(\mu)$ and $\tilde\mu_t:=\varphi_t(\tilde\mu)$, with the same initial datum $m_0$. The choices $m_0=\delta_z$ and $m_0=\delta_{z_1}-\delta_{z_2}$ give the two estimates. Set
    \begin{align*}
        \Delta
        &:= W_{\psi_1 + \psi_2}(\mu,\tilde\mu),
        &
        E_0
        &:= W_{\psi_1 + \psi_2}(m_0,0).
    \end{align*}
    We write $h_t(x)=\exp\left(\int_0^t \brac{K_\alpha(x,\cdot),\mu_s}\dd s\right)$ and use the analogous notation $\tilde h_t$.

    First remove the weights. For $f\in \FC_{\psi_1}^{1}$
    \begin{equation}\label{eq:proof:final:step1}
        \brac{f,m_t-\tilde m_t}
        =\brac{f h_t^{-1},h_t m_t - \tilde h_t \tilde m_t}
        +\brac{f(h_t^{-1}-\tilde h_t^{-1}),\tilde h_t\tilde m_t}.
    \end{equation}
    By Lemma~\ref{lem:technical_G_h_dos}-$(i)$ the first term is bounded by
    \[
        \left|\brac{fh_t^{-1}, h_tm_t-\tilde h_t\tilde m_t}\right|
        \lesssim \left(1 + t\|\mu\|\right)W^1_{\psi_1}\left(h_tm_t , \tilde h_t\tilde m_t\right)\Wnorm{f}{\psi_1,1}.
    \]
    For the second term, Lemma~\ref{lem:technical_G_h}--$(i)$ controls the difference of the inverse weights by
    \begin{align*}
        \Wnorm{f(h_t^{-1}-\tilde h_t^{-1})}{\psi_1+\psi_2}
        &\lesssim (1 + t\|\mu + \tilde\mu\|) \int_0^t W_{\psi_1}(\mu_s,\tilde\mu_s)\dd s\, \Wnorm{f}{\psi_1,1}.
    \end{align*}
    Proposition~\ref{prop:wass-sce} controls the time integral, while Lemma~\ref{lem:wass-rescaled-derivative} controls $\tilde h_t\tilde m_t$ in the weights $\psi_1$ and $\psi_2$. Therefore
    \[
        \left|\brac{f(h_t^{-1}-\tilde h_t^{-1}),\tilde h_t\tilde m_t}\right|
        \leq  C_1(T,\mu+\tilde\mu)W_{\psi_1}(\mu,\tilde\mu)E_0
        \Wnorm{f}{\psi_1,1}.
    \]
    with 
    \[
    C_1(T,\nu) = cT (1 +  T\|\nu\|_{\psi_3})  e^{c T\|\nu\|_{\psi_2}(1 + T\|\nu\|)}.
    \]
    Taking the supremum over $\Wnorm{f}{\psi_1}^1\le 1$ in \eqref{eq:proof:final:step1} gives
    \begin{equation}\label{eq:unweight-final-wass}
        W^1_{\psi_1}(m_t,\tilde m_t)
        \leq  c(1 + t\|\mu\|) W^1_{\psi_1}(h_tm_t,\tilde h_t\tilde m_t) 
        + C_1(T,\mu+\tilde\mu)W_{\psi_1}(\mu,\tilde\mu)E_0.
    \end{equation}

    It remains to compare the weighted derivative measures. Define
    \begin{align*}
        H_sf(x,y) := K_\alpha(x,y)f(x+y)G_s(x,y),
        \qquad
        I_sf(x,y) := K_\alpha(x,y)f(x)h_s^{-1}(y),
    \end{align*}
    with tildes denoting the same quantities built from $\tilde\mu$. Subtracting the two rescaled derivative equations and using $h_0 m_0 = \tilde h_0 m_0$ gives
    \begin{align*}
        \brac{f,h_t m_t - \tilde h_t \tilde m_t}
        = \int_0^t \left(R_s^{(1)}+R_s^{(2)}+R_s^{(3)}\right)\dd s,
    \end{align*}
    where
    \begin{align*}
        R_s^{(1)}
        &:= \brac{H_sf-I_sf,h_s\mu_s\otimes \left(h_s m_s - \tilde h_s \tilde m_s \right)},\\
        R_s^{(2)}
        &:= \brac{H_sf-I_sf, \left(h_s\mu_s-\tilde h_s\tilde\mu_s\right)\otimes\tilde h_s\tilde m_s},\\
        R_s^{(3)}
        &:= \brac{(H_s-\tilde H_s)f+(\tilde I_s-I_s)f,\tilde h_s\tilde\mu_s\otimes\tilde h_s\tilde m_s}.
    \end{align*}

   For the first term, Lemma~\ref{lem:technical_G_h_dos}--$(i)$, $(ii)$
    and Proposition~\ref{prop:exponential-moments} give
    \[
        \Wnorm{y\mapsto\brac{(H_s-I_s)f(\cdot,y),h_s\mu_s}}{\psi_1,1}
        \lesssim (1+s\|\mu\|)\|\mu\|_{\psi_2}
        \Wnorm{f}{\psi_1,1}.
    \]
    Here the common factor
    $e^{-\int_0^s\brac{\psi_1,\mu_r}\dd r}$ in the inverse-weight
    estimates is retained before applying the exponential-moment bound.
    Consequently,
    \[
        |R_s^{(1)}|
        \lesssim (1+s\|\mu+\tilde\mu\|)\|\mu+\tilde\mu\|_{\psi_2}
        W^1_{\psi_1}(h_sm_s,\tilde h_s\tilde m_s)
        \Wnorm{f}{\psi_1,1}.
    \]
    
    For the second term, define
    \[
        q_s(x):=\brac{(H_s-I_s)f(x,\cdot),\tilde h_s\tilde m_s}.
    \]
    Lemma~\ref{lem:technical_G_h}--$(i)$, $(ii)$ gives
    \begin{align*}
        |(H_s-I_s)f(x,y)|
        &\lesssim (1+\psi_1(y))(1+\psi_2(x))\Wnorm{f}{\psi_1,1},\\
        |\partial_y((H_s-I_s)f)(x,y)|
        &\lesssim (1+s\|\mu\|)y^{-1-\alpha}
        (1+\psi_2(x))\Wnorm{f}{\psi_1,1}.
    \end{align*}
    Moreover, the preceding technical estimates,
    Lemma~\ref{lem:technical_G_h_dos}--$(iii)$ and direct
    differentiation of $I_sf$ give
    \begin{align*}
        |\partial_x((H_s-I_s)f)(x,y)|
        &\lesssim (1+s\|\mu\|)(1+\psi_1(y))
        (x^{-1-\alpha}+x^{-1-2\alpha})\Wnorm{f}{\psi_1,1},\\
        |\partial_{xy}((H_s-I_s)f)(x,y)|
        &\lesssim (1+s\|\mu\|)y^{-1-\alpha}
        (x^{-1-\alpha}+x^{-1-2\alpha})\Wnorm{f}{\psi_1,1}.
    \end{align*}
    Pairing these bounds in the $y$-variable controls both $q_s$ and
    its derivative in the $x$-variable. Thus
    \[
        \Wnorm{q_s}{\psi_1+\psi_2}
        \lesssim (1+s\|\mu\|)
        W_{\psi_1}(\tilde h_s\tilde m_s,0)\Wnorm{f}{\psi_1,1}.
    \]
    Propositions~\ref{prop:wass-rescaled-SCE}
    and~\ref{lem:wass-rescaled-derivative}, with \eqref{eq:sum-weights},
    therefore yield
    \begin{equation*}
        |R_s^{(2)}|
        \leq W_{\psi_1+\psi_2}(h_s\mu_s,\tilde h_s\tilde\mu_s)
        \Wnorm{q_s}{\psi_1+\psi_2}
        \leq C_2(T,\mu+\tilde\mu)\Delta
        W_{\psi_1}(m_0,0)\Wnorm{f}{\psi_1,1},
    \end{equation*}
    where we may take:
    \[
    C_2(T,\nu) = c e^{cT\|\nu\|_{\psi_2}(1 + T\|\nu\|)}.
    \]

    Finally, split $R_s^{(3)}$ according to the two differences of
    coefficients. Lemma~\ref{lem:technical_G_h}--$(iii)$ controls
    $H_s-\tilde H_s$, whereas part $(i)$ controls
    $\tilde I_s-I_s$. In both cases the weight in the signed variable
    is $\psi_1+\psi_2$, and the positive variable requires only
    the moment $\psi_2$. The exponential-moment estimate gives
    \begin{equation*}
        |R_s^{(3)}|
        \lesssim (1+T\|\mu+\tilde\mu\|)\|\mu+\tilde\mu\|_{\psi_2}
        e^{cT\|\mu+\tilde\mu\|_{\psi_2}}
        \times W_{\psi_1+\psi_2}(\tilde h_s\tilde m_s,0)
        \int_0^sW_{\psi_1}(\mu_r,\tilde\mu_r)\dd r\,
        \Wnorm{f}{\psi_1,1}.
    \end{equation*}
    Propositions~\ref{prop:wass-sce}
    and~\ref{lem:wass-rescaled-derivative} therefore imply
    \[
        |R_s^{(3)}|
        \leq C_3(T,\mu+\tilde\mu)
        W_{\psi_1}(\mu,\tilde\mu)E_0\Wnorm{f}{\psi_1,1},
    \]
    where one may take after adjusting $c$ and absorbing one factor $T\|\nu\|_{\psi_2}(1 + T\|\nu\|)$ in the exponential,
    \begin{equation*}
        C_3(T,\nu) = T\|\nu\|_{\psi_3} C_2(T,\nu).
    \end{equation*}
    
    Thus $R_s^{(1)}$ is the only term propagated by Gronwall, while $R_s^{(2)}$ and $R_s^{(3)}$ are bounded directly. Taking the supremum over $\Wnorm{f}{\psi_1,1}\leq1$ and applying Gronwall's lemma, we obtain
    \begin{align*}
        W^1_{\psi_1}(h_tm_t,\tilde h_t\tilde m_t)
        \leq C(T,\mu+\tilde\mu)\left(\Delta W_{\psi_1}(m_0,0) + E_0 W_{\psi_1}\left(\mu,\tilde{\mu}\right)\right),
    \end{align*}
    where
    \[
    C(T,\nu) =  c(1 + T)(1 + T\|\nu\|_{\psi_3}) e^{c T\|\nu\|_{\psi_2}(1 + T\|\nu\|)}.
    \]
    Together with \eqref{eq:unweight-final-wass}, this proves
    \begin{align*}
        W^1_{\psi_1}(m_t,\tilde m_t)
        \leq C(T,\mu+\tilde\mu)\left(\Delta W_{\psi_1}(m_0,0) + E_0 W_{\psi_1}\left(\mu,\tilde{\mu}\right)\right).
    \end{align*}
    If $m_0=\delta_z$, then $E_0\lesssim 1+\psi_1(z)+\psi_2(z)$, which gives the first estimate. If $m_0=\delta_{z_1}-\delta_{z_2}$, then $E_0 \le d_{\psi_1 + \psi_2}(z_1,z_2)$, which gives the second estimate.
\end{proof}

\section{Proofs of the main results}\label{Section:main-results}
\subsection{Proof of Theorem \ref{thm:samesame}}\label{sec:thm:samesame}
We start with the technical lemma below. We split the expression of the Smoluchowski coagulation equation generator $\LC$ into two parts:
    \begin{align*}
        \LC F(\mu) &= \LC_1 F(\mu) + \LC_2 F(\mu),
    \end{align*}
    with
    \begin{align*}
            \LC_1 F(\mu) =& \dfrac12 \int_{\RRP \times \RRP} K_{\alpha}\delta_\mu F(\mu, \cdot)(x,y) \left(\mu^{\otimes 2} - \dfrac1N \Diag_{\#}\mu \right) (\dd x,\dd y) \\
            \LC_2 F(\mu) =& \dfrac{1}{2 N}\int_{\RRP} K_{\alpha}\delta_\mu F(\mu, \cdot)(x,x)\mu(\dd x).
    \end{align*}
By assumption $\A$, the second part of the SCE operator satisfies
\[
\left|\LC_2 F(\mu)\right| \le \Wnorm{\delta_\mu F}{\A} \dfrac{(1 + \|\mu\|_{\psi_2})\|\mu\|_{\psi_2}}{N}
\]
and therefore by Proposition \ref{prop:wass-derivative} and absorbing one term $T\|\mu\|_{\psi_2}$ into the exponential:
    \[
        \int_0^t \left|\LC_2 F\circ \varphi_s (\mu)\right| \dd s \le \Wnorm{\delta_\mu F}{\A} \dfrac{c(1 + \|\mu\|_{\psi_2}) e^{c T\|\mu\|_{\psi_2}(1 + T\|\mu\|)}}{N}.
    \]
The first part is the one involved in the comparison and is treated in the Lemma below.
\begin{lemma}\label{lem:bounds_diff_gen}
Let $\mu \in \MC_N^{+,p}(\RRP)$ and denote by $V_3 : \mu \mapsto \brac{1 + \psi_3, \mu} = \|\mu\|_{\psi_3}$ we have:
\[
    \left| \left(L^N - \LC_1\right)F\circ \varphi_s (\mu) \right| \leq \Wnorm{\delta_\mu F}{\A}\dfrac{ C(T,\mu)}{N}\left( \|\mu\|_{\psi_3}^2 - (1 + T\|\mu\|_{\psi_3})L^N V_3(\mu) \right).
\]
Where the constant is of the form
\[
C(T,\nu) =  c(1 + T)(1 + \|\nu\|_{\psi_2}) e^{c T\|\nu\|_{\psi_2}(1 + T\|\nu\|)}.
\]
\end{lemma}

\begin{proof}
    Set
    \begin{align*}
        \mu_N^\theta(x,y) = \mu + \frac{\theta}{N}\left(\delta_{x + y} - \delta_x - \delta_y \right).
    \end{align*}
    For all pairs $(x,y)$ participating in $\mu^{\otimes 2} - \tfrac1N \Diag_{\#}\mu$, the measure $\mu_N^\theta(x,y)$ is positive. We have:
    \begin{multline*}
        \left(L^N - \LC_1\right)F\circ \varphi_s (\mu) 
        = \dfrac12\int_0^1 \int_{\RR^+ \times \RR^+} K_{\alpha}(x,y)\\ 
        \times  \left\langle  \delta_\mu F\circ \varphi_s\left(\mu_N^\theta(x,y);\cdot\right)- \delta_\mu F\circ \varphi_s\left(\mu;\cdot\right),\delta_{x+y} - \delta_x - \delta_y\right\rangle  \\
        \left(\mu^{\otimes 2} - \dfrac1N \Diag_{\#}\mu \right) (\dd x,\dd y)\dd \theta .
    \end{multline*}
    
    If $\psi > 0$ is sub-additive then $\|\mu_N^\theta(x,y)\|_{\psi} \leq \|\mu\|_{\psi}$.  By the assumptions on $F$ and the regularity of the derivative flow,
    Proposition~\ref{prop:chain-rule}, with source weight $\psi_3$,
    target weight $\psi_2$, and test weight $\psi_1$, gives
    \begin{align*}
        \delta_\mu F\circ \varphi_s(\nu;z) = \brac{\delta_\mu F(\varphi_s(\nu);\cdot),\delta_\mu \varphi_s (\nu;z)}.
    \end{align*}
    Therefore,
    \begin{subequations}
    \begin{align}
        \delta_\mu F\circ \varphi_s\left(\mu_N^\theta(x,y);z\right) - &\delta_\mu F\circ \varphi_s(\mu;z) \nonumber \\
        =& \brac{\delta_\mu F\left(\varphi_s\left(\mu_N^\theta(x,y)\right);\cdot\right) - \delta_\mu F(\varphi_s(\mu);\cdot),\delta_\mu \varphi_s (\mu;z)}\label{proof:eq:samesame-part1} \\
         &+ \brac{\delta_\mu F\left(\varphi_s\left(\mu_N^\theta(x,y)\right);\cdot\right),\delta_\mu \varphi_s \left(\mu_N^\theta(x,y);z\right) - \delta_\mu \varphi_s (\mu;z)}\label{proof:eq:samesame-part2} .
    \end{align}
    \end{subequations}
    
    For \eqref{proof:eq:samesame-part1}, the assumptions on
    $\delta_\mu F$ imply that
    \[
        z\mapsto
        \delta_\mu F(\varphi_s(\mu_N^\theta);z)
        -\delta_\mu F(\varphi_s(\mu);z)
    \]
    has $\FC_{\psi_1}$-norm controlled by $\Wnorm{\delta_\mu F}{\A}(1 + \|\mu\|_{\psi_2})W_{\psi_1}\left(\varphi_s\left(\mu_N^\theta(x,y)\right), \varphi_s\left(\mu\right) \right)$. Hence Proposition \ref{prop:wass-derivative} gives
    \begin{multline*}
        \left|
        \brac{\delta_\mu F(\varphi_s(\mu_N^\theta);\cdot)
        -\delta_\mu F(\varphi_s(\mu);\cdot),
        \delta_\mu\varphi_s(\mu;z)}
        \right| \\
        \leq \Wnorm{\delta_\mu F}{\A}W_{\psi_1}(\varphi_s(\mu_N^\theta(x,y)), \varphi_s\left(\mu\right) )    C_1(T,\mu) (1 + \psi_1(z)).
    \end{multline*}
    where
    \[
    C_1(T,\nu) =  (1 + \|\nu\|_{\psi_2}) e^{cT\|\nu\|_{\psi_2}(1 + T\|\nu\|)}
    \]
    Proposition \ref{prop:wass-sce} gives
    \begin{multline*}
         \left|
        \brac{\delta_\mu F(\varphi_s(\mu_N^\theta);\cdot)
        -\delta_\mu F(\varphi_s(\mu);\cdot),
        \delta_\mu\varphi_s(\mu;z)}
        \right| \\
        \leq \Wnorm{\delta_\mu F}{\A}W_{\psi_1}(\mu_N^\theta(x,y),\mu) C_1(T,\mu) (1 + \psi_1(z)).
    \end{multline*}
    Moreover for any admissible weight $\psi$:
    \begin{align*}
        W_{\psi}(\mu_N^\theta(x,y),\mu) &= \dfrac{\theta}{N} \sup\limits_{\Wnorm{f}{\psi} \le 1} f(x+y) - f(x) - f(y) \\
        &\leq \dfrac{\theta}{N} c(1 + \psi(x) + \psi(y))
    \end{align*}
    Thus
    \begin{multline*}
        \left|\brac{\delta_\mu F\left(\varphi_s\left(\mu_N^\theta(x,y)\right),\cdot\right) - \delta_\mu F(\varphi_s(\mu);\cdot),\delta_\mu \varphi_s (\mu;z)}\right| \\
        \leq \Wnorm{\delta_\mu F}{\A}\dfrac{C_1(T,\mu)}{N}(1 + \psi_1(z))(1 + \psi_1(x) + \psi_1(y)).
    \end{multline*}
    By denoting
    \[
    V^\theta_1 : (x,y,z) \mapsto \left|\brac{\delta_\mu F\left(\varphi_s\left(\mu_N^\theta(x,y)\right),\cdot\right) - \delta_\mu F(\varphi_s(\mu);\cdot),\delta_\mu \varphi_s (\mu;z)}\right|
    \]
    we see that:
    \begin{multline*}
    \int_0^1 \int_{\RR^+ \times \RR^+} K_{\alpha}(x,y) \brac{V^\theta_1 (x,y,\cdot),\delta_{x + y} + \delta_x + \delta_y}  \left(\mu^{\otimes 2} - \dfrac1N \Diag_{\#}\mu \right) (\dd x,\dd y)\dd \theta \\
    \le \Wnorm{\delta_\mu F}{\A}\dfrac{C_1(T,\mu)\|\mu\|_{\psi_3}^2}{N}.
    \end{multline*}
    
    For the term \eqref{proof:eq:samesame-part2}, since $z \mapsto \delta_\mu F\left(\varphi_s\left(\mu_N^\theta(x,y)\right);z\right)$ is in $\FC^1_{\psi_1}$ by assumption, let $C_2(T,\nu) = (1 + T\|\nu\|_{\psi_3})C(T,\nu)$ we have by point $1$ of Proposition \ref{prop:wass-final}:
    \begin{multline*}
        \left|\brac{\delta_\mu F\left(\varphi_s\left(\mu_N^\theta(x,y)\right);\cdot\right),\delta_\mu \varphi_s \left(\mu_N^\theta(x,y);z\right) - \delta_\mu \varphi_s (\mu;z)} \right|
         \leq \Wnorm{\delta_\mu F}{\A}C_2(T,\mu)\\
         \times \left(W_{\psi_1 + \psi_2}\left(\mu_N^\theta(x,y), \mu\right)(1 + \psi_1(z)) + W_{\psi_1}\left(\mu_N^\theta(x,y), \mu\right)(1 + \psi_2(z)) \right).
    \end{multline*}
    Hence for $z \in \lbrace x+y,x,y\rbrace$
    \begin{multline*}
        \left|\brac{\delta_\mu F\left(\varphi_s\left(\mu_N^\theta(x,y)\right);\cdot\right),\delta_\mu \varphi_s \left(\mu_N^\theta(x,y);z\right) - \delta_\mu \varphi_s (\mu;z)}\right| \\
        \le \Wnorm{\delta_\mu F}{\A}\dfrac{C_2(T,\mu)}{N}(1 + \psi_2(x) + \psi_2(y))(1 + \psi_1(z))\\
        \le \Wnorm{\delta_\mu F}{\A} \dfrac{C_2(T,\mu)}{N}(1 + \psi_3(x) + \psi_3(y)).
    \end{multline*}
    Note that by convexity of $\psi_3$ there exists a constant $c_3 = (1 - 2^{-1 - 3\alpha})^{-1} > 1$ such that:
    \[
   \psi_3(x) + \psi_3(y) \le - c_3 \left( \psi_3(x + y) - \psi_3(y) - \psi_3(x) \right).
    \]
    Let us denote
    \[
    V_2^\theta : (x,y,z) \mapsto \left|\brac{\delta_\mu F\left(\varphi_s\left(\mu_N^\theta(x,y)\right);\cdot\right),\delta_\mu \varphi_s \left(\mu_N^\theta(x,y);z\right) - \delta_\mu \varphi_s (\mu;z)}\right|
    \]
    It follows that :
    \begin{multline*}
        \int_0^1 \int_{\RR^+ \times \RR^+} K_{\alpha}(x,y)  \brac{V_2^\theta(x,y,\cdot), \delta_{x + y} + \delta_x +\delta_y} 
    \left(\mu^{\otimes 2} - \dfrac1N \Diag_{\#}\mu \right) (\dd x,\dd y)\dd \theta \\
    \le \Wnorm{\delta_\mu F}{\A} \dfrac{C_2(T,\mu)}{N}
    \times \dfrac12\int_{\RR^+ \times \RR^+} K_{\alpha}(x,y)(1 + \psi_3(x) + \psi_3(y) - \psi_3(x + y)) \\
    \left(\mu^{\otimes 2} - \dfrac1N \Diag_{\#}\mu \right) (\dd x,\dd y) 
    = \Wnorm{\delta_\mu F}{\A} \dfrac{C_2(T,\mu)}{N} \times  -L^N V_3(\mu),
    \end{multline*}
    where $V_3 : \mu \mapsto \brac{1 + \psi_3, \mu}$.
\end{proof}

\begin{proof}[Proof of Theorem \ref{thm:samesame}]
 We start from the Duhamel formula \eqref{eq:semi-group-relation}:
 \begin{align*}
\E{F(\Phi_t(\mu^N))} - F(\varphi_t\left(\mu^N\right)) = \int_0^t \E{\left(L^N - \LC\right)F \circ \varphi_s (\Phi_{t-s}(\mu^N))}\dd s.
\end{align*}
Lemma \ref{lem:bounds_diff_gen} along with the diagonal component contribution give:
\begin{multline*}
    \left|\vphantom{A^2} \E{F(\Phi_t(\mu^N))} - F(\varphi_t\left(\mu^N\right))\right| \leq \dfrac{\Wnorm{\delta_\mu F}{\A}}{N} 
    \times \int_0^t \mathbb{E}\left[\tilde C\left(T,\Phi_{t -s }(\mu^N)\right) \right.\\
    \times \left.\left(1 + \left\|\Phi_{t-s}(\mu^N) \right\|_{\psi_3}^2 + \left( 1 + T\left\|\Phi_{t-s}(\mu^N) \right\|_{\psi_3}\right)  \times - L^N V_3\left(\Phi_{t -s }(\mu^N) \right)\right)\right] \dd s.
\end{multline*}
The Marcus-Lushnikov process defines a positive measure and, for every sub-additive $\psi > 0$ and time $t \ge 0$, satisfies $\brac{\psi,\Phi_t(\mu^N)} \leq \brac{\psi,\mu^N}$ almost surely. In particular 
\begin{align*}
   \left\|\Phi_{t-s}(\mu^N) \right\|_{\psi_3} \le \|\mu^N\|_{\psi_3}, \qquad \tilde C\left(T,\Phi_t(\mu^N) \right) \le  \tilde C\left(T,\mu^N\right)\,\text{a.s}\ .
\end{align*}
We have
\begin{multline*}
    \left|\vphantom{A^2} \E{F(\Phi_t(\mu^N))} - F(\varphi_t\left(\mu^N\right))\right|\\
    \le \dfrac{\Wnorm{\delta_\mu F}{\A} \tilde C(T,\mu^N) }{N} \left( T(1 + \|\mu^N\|_{\psi_3}^2) + (1 + T\|\mu^N\|_{\psi_3}) \times -\int_0^t \E{L^N V_3 \left( \Phi_s(\mu^N)\right)} \dd s\right)\\
    \le \dfrac{\Wnorm{\delta_\mu F}{\A} (1 + \|\mu^N\|_{\psi_3}) C\left(T,\mu^N\right)}{N},
\end{multline*}
where
\[
C(T,\nu) = c(1 + T)(1 + T\|\nu\|_{\psi_3})(1 + \|\nu\|_{\psi_2})e^{cT\|\nu\|_{\psi_2}(1 + T\|\nu\|)}
\]
\end{proof}

\subsection{Proof of Theorem \ref{thm:difdif}}\label{sec:thm:difdif}

\begin{proof}[Proof of Theorem \ref{thm:difdif}]
    Recall the bracket notation $[\mu,\nu]^\theta = \theta \mu + (1-\theta)\nu$. Since $F\circ \varphi_t$ has a flat derivative,
    \[
        F\left(\varphi_t(\mu)\right) - F\left(\varphi_t(\nu)\right) 
        = \int_0^1 \int_{\RRP} \brac{\delta_\mu F\left(\varphi_t \left(\left[\mu,\nu \right]^\theta\right);\cdot \right), \delta_\mu \varphi_t\left(\left[\mu,\nu \right]^\theta;z \right)} \left(\mu - \nu\right)(\dd z)\dd \theta.
    \]
    We split it into three parts:
    \[
        F\left(\varphi_t(\mu)\right) - F\left(\varphi_t(\nu)\right) = A_t + B_t +  \int_0^1 \int_{\RRP} \brac{\delta_\mu F\left(\varphi_t \left(\nu\right);\cdot \right), \delta_\mu \varphi_t\left(\nu;z \right)} \left(\mu - \nu\right)(\dd z)\dd \theta. 
    \]
    With 
    \begin{align*}
        A_t &= \int_0^1 \int_{\RRP} \brac{\delta_\mu F\left(\varphi_t \left(\left[\mu,\nu \right]^\theta\right);\cdot \right)- \delta_\mu F\left(\varphi_t \left(\nu\right);\cdot \right), \delta_\mu \varphi_t\left(\left[\mu,\nu \right]^\theta;z \right)} \left(\mu - \nu\right)(\dd z)\dd \theta \\
        B_t &= \int_0^1 \int_{\RRP} \brac{\delta_\mu F\left(\varphi_t \left(\nu\right);\cdot \right), \delta_\mu \varphi_t\left(\left[\mu,\nu \right]^\theta;z \right) -  \delta_\mu \varphi_t\left(\nu;z \right)} \left(\mu - \nu\right)(\dd z)\dd \theta.
    \end{align*}
    For the first term, by assumption on $F$ and monotonicity of the SCE flow the function
    \begin{align*}
        x \mapsto \delta_\mu F\left(\varphi_t \left(\left[\mu,\nu \right]^\theta\right);x \right)- \delta_\mu F\left(\varphi_t \left(\nu\right);x \right)
    \end{align*}
    is controlled in $\FC_{\psi_1}$-norm by $ \Wnorm{\delta_\mu F}{\A}  ( 1 + \|\mu + \nu\|_{\psi_2})W_{\psi_1}\left(\varphi_t \left(\left[\mu,\nu \right]^\theta\right), \varphi_t \left(\nu \right)\right)$. By Proposition \ref{prop:wass-sce}:
    \begin{align*}
        W_{\psi_1}\left(\varphi_t \left(\left[\mu,\nu \right]^\theta\right), \varphi_t \left(\nu \right)\right) \leq e^{cT\|\mu + \nu\|_{\psi_2}\left(1 + T \|\mu + \nu\| \right)}W_{\psi_1}\left( \mu, \nu\right).
    \end{align*}
    We deduce by Proposition \ref{prop:wass-derivative}, that the function
    \begin{align*}
         z \mapsto \brac{\delta_\mu F\left(\varphi_t \left(\left[\mu,\nu \right]^\theta\right);\cdot \right)- \delta_\mu F\left(\varphi_t \left(\nu\right);\cdot \right), \delta_\mu \varphi_t\left(\left[\mu,\nu \right]^\theta;z \right)}
    \end{align*}
    is controlled in $\FC_{\psi_1}$-norm by $ \Wnorm{\delta_\mu F}{\A} C_1(T, \mu + \nu) W_{\psi_1}\left( \mu , \nu\right)$ where
    \[
        C_1(T, m) = c (1 + \|m\|_{\psi_2}) e^{c T\|m\|_{\psi_2}\left(1 + T \|m\| \right)}.
    \]
    Hence 
    \[
        \left| A_t \right|\le \Wnorm{\delta_\mu F}{\A} C_1(T, \mu + \nu) W_{\psi_1}\left( \mu , \nu\right)^2.
    \]
    For the second term, by Proposition \ref{prop:wass-final} and the assumption on $\delta_\mu F$, the function
    \begin{align*}
        z \mapsto \brac{\delta_\mu F\left(\varphi_t \left(\nu\right);\cdot \right), \delta_\mu \varphi_t\left(\left[\mu,\nu \right]^\theta;z \right) -  \delta_\mu \varphi_t\left(\nu;z \right)}
    \end{align*}
    is controlled in $\FC_{\psi_1 + \psi_2}$-norm by
    \[
     \Wnorm{\delta_\mu F}{\A}( 1+ \|\mu+ \nu\|_{\psi_2})C_2(T,\mu + \nu)\left(W_{\psi_1}\left(\mu,\nu\right) +W_{\psi_2}\left(\mu,\nu\right) \right),
    \] 
    where $C_2$ is the constant from Proposition \ref{prop:wass-final}. Hence the second term satisfies
    \[
        \left| B_t \right|\le \Wnorm{\delta_\mu F}{\A} C(T,\mu + \nu)\left(W_{\psi_1}\left(\mu,\nu\right) +W_{\psi_2}\left(\mu,\nu\right) \right)^2.
    \]
    Combining the bounds gives
    \begin{multline*}
        \left|F\left(\varphi_t(\mu)\right) - F\left(\varphi_t(\nu)\right)
        - \int_{\RRP} \brac{\vphantom{A^2}\delta_\mu F\left(\varphi_t \left(\nu\right);\cdot \right), \delta_\mu \varphi_t\left(\nu;z \right)} \left(\mu - \nu\right)(\dd z)\right| \\
        \leq \Wnorm{\delta_\mu F}{\A} C(T,\mu + \nu)\left(W_{\psi_1}\left(\mu,\nu\right) + W_{\psi_2}\left(\mu,\nu\right) \right)^2.
    \end{multline*}
\end{proof}
\subsection{Proof of Theorem \ref{thm:main-result-2}}
Conditioning on $\mu^N$ Theorem~\ref{thm:main-result} gives:
\begin{multline*}
        \left| \mathbb{E}_{\mu^N}\left[F\left(\Phi_t\left(\mu^N\right)\right)\right] - F\left(\varphi_t(\mu)\right) - \int_{\RRP} \left\langle \delta_\mu F\left(\varphi_t \left(\mu\right);\cdot \right), \delta_\mu \varphi_t\left(\mu;z \right) \right\rangle \left(\mu^N - \mu\right)(\dd z) \right| \\
        \leq \Wnorm{\delta_\mu F}{\A} R_T^N
\end{multline*}
where
\begin{align*}
    R_T^N = \frac{\|\mu^N\|_{\psi_3} C_1(T,\mu^N)}{N}  + C_1(T,\mu + \mu^N)\left(W_{\psi_1}\left(\mu,\mu^N\right) + W_{\psi_2}\left(\mu,\mu^N\right) \right)^2
\end{align*}
Where the constant $C_1(T,m)$ is given by:
\[
C_1(T,\nu) = c(1 + T)(1 + T\|\nu\|_{\psi_3})(1 + \|\nu\|_{\psi_2})e^{cT\|\nu\|_{\psi_2}(1 + T\|\nu\|)}.
\]
Note that since $\mu$ and $\mu^N$ are probability measures $\|\mu\| = \|\mu^N\| = 1$, we can therefore absorb the factors in $T$ in the exponential term and bound $C_1$ by:
\[
\|m\|_{\psi_3}C(T,m) := \|m\|_{\psi_3} c \|m\|_{\psi_2} e^{c\|m\|_{\psi_2} T(1 + T\|m\|)} .
\]
Let us denote $A_N$ the event $\lbrace \brac{\psi_2,\mu^N} \le \brac{\psi_2,\mu} + 1 \rbrace$. Let us state and prove a technical lemma controlling $\Prob{A_N^c}$.
\begin{lemma}
    Let $p \ge 4$ and assume $\mu \in \PC_{\psi_p}(\RRP)$ we have for all $N \ge 1$:
    \[
    \Prob{A_N^c} \le C_p N^{-p/4} \|\mu\|_{\psi_p } 
    \]
\end{lemma}
\begin{proof}
    By Markov's inequality:
    \[
    \Prob{A_N^c} \le \E{|\brac{\psi_2, \mu^N - \mu}|^{\tfrac{p}{2}}}.
    \]
    By Rosenthal's inequality \cite{rosenthal1970subspaces}:
    \[
    \Prob{A_N^c} \le C_p \left(N^{-\tfrac{p}{2} + 1} \brac{\psi_p, \mu} + N^{-\tfrac{p}{4}} \brac{\psi_4, \mu}^{\tfrac{p}{4}} \right).
    \]
    Finally applying Jensen's inequality gives the result.
\end{proof}

On $A_N$ we have:
\[
C\left(T,\mu^N \right) \le C\left(T, 2\mu \right), \qquad C\left(T,\mu^N + \mu\right) \le C\left(T,3\mu\right)\qquad a.s. 
\]

Using Holder's inequality:
\begin{multline*}
\E{R_T^N \mathds{1}_{A_N}} \le \dfrac{\|\mu\|_{\psi_6}C(T,2\mu)}{N} \\
+ C(T,3\mu)\E{\|\mu + \mu^N\|_{\psi_3}^{\tfrac{p'}{p' - 1}}}^{1 - \tfrac{1}{p'}}\E{\left(W_{\psi_1}\left(\mu,\mu^N\right) + W_{\psi_2}\left(\mu,\mu^N\right) \right)^{2p'}}^{\tfrac{1}{p'}}
\end{multline*}
Choosing $p' = 7/4$ lets us apply the Wasserstein rate result in Corollary \ref{cor:wass-square-bound} with the assumption $\mu \in \MC^+_{\psi_{7 + \varepsilon}}(\RRP)$, it yields:
\[
\E{R_T^N \mathds{1}_{A_N}} \le \dfrac{\|\mu\|_{\psi_{7 + \varepsilon}}C_\varepsilon(T,\mu)}{N}.
\]
We treat now the first order term conditioned to $A_N$, observe that under full expectation it vanishes and therefore:
\begin{multline*}
\E{\int_{\RRP} \left\langle \delta_\mu F\left(\varphi_t \left(\mu\right);\cdot \right), \delta_\mu \varphi_t\left(\mu;z \right) \right\rangle \left(\mu^N - \mu\right)(\dd z) \mathds{1}_{A_N}}
\\
= - \E{\int_{\RRP} \left\langle \delta_\mu F\left(\varphi_t \left(\mu\right);\cdot \right), \delta_\mu \varphi_t\left(\mu;z \right) \right\rangle \left(\mu^N - \mu\right)(\dd z) \mathds{1}_{A_N^c}}.
\end{multline*}
Then thanks to Proposition \ref{prop:wass-derivative} and Assumption $\A$ we know that the function of $z$ in the integral above has $\FC_{\psi_1}$-norm controlled by $C(T,\mu)$ it follows  that :
\begin{multline*}
    \left| \E{\int_{\RRP} \left\langle \delta_\mu F\left(\varphi_t \left(\mu\right);\cdot \right), \delta_\mu \varphi_t\left(\mu;z \right) \right\rangle \left(\mu^N - \mu\right)(\dd z) \mathds{1}_{A_N^c}}\right| \\
    \le \Wnorm{\delta_\mu F}{\A}C(T,\mu)\E{\brac{1 + \psi_1,\mu^N + \mu} \mathds{1}_{A_N^c}}.
\end{multline*}
Note that for $q' > 1$, by Jensen's inequality $\brac{\psi_1, \mu^N}^{q'}\lesssim \brac{\psi_{q'},\mu^N}$, therefore by Hölder's inequality:
\begin{multline*}
\left| \E{\int_{\RRP} \left\langle \delta_\mu F\left(\varphi_t \left(\mu\right);\cdot \right), \delta_\mu \varphi_t\left(\mu;z \right) \right\rangle \left(\mu^N - \mu\right)(\dd z) \mathds{1}_{A_N^c}}\right| \\
    \le \Wnorm{\delta_\mu F}{\A}C(T,\mu)\|\mu\|_{\psi_{q'}}^{\tfrac{1}{q'}} \Prob{A_N^c}^{1 - \tfrac{1}{q'}}.
\end{multline*}
Choosing $q' = 7/3$ and applying the previous Lemma for $p = 7$ gives:
\begin{multline*}
\left| \E{\int_{\RRP} \left\langle \delta_\mu F\left(\varphi_t \left(\mu\right);\cdot \right), \delta_\mu \varphi_t\left(\mu;z \right) \right\rangle \left(\mu^N - \mu\right)(\dd z) \mathds{1}_{A_N^c}}\right| \\
    \le \Wnorm{\delta_\mu F}{\A}\dfrac{C(T,\mu) \|\mu\|_{\psi_{7/3}}^{\tfrac{3}{7}}\|\mu\|_{\psi_7}^{4/7}}{N} 
    \le \Wnorm{\delta_\mu F}{\A}\dfrac{C(T,\mu)\|\mu\|_{\psi_7}}{N}.
\end{multline*}
The expected value of the difference conditioned on $A_N^c$ through a direct bound satisfies:
\begin{multline*}
    \left| \mathbb{E}\left[\left(F\left(\Phi_t\left(\mu^N\right)\right)-  F\left(\varphi_t(\mu)\right)\right)\mathds{1}_{A_N^c}\right] \right| \\
    = \left|\E{\int_0^1 \int_{\RRP} \delta_{\mu}F\left(\left[\Phi_t(\mu^N) ,\varphi_t(\mu) \right]^{\theta} ; z\right) \left(\Phi_t(\mu^N) - \varphi_t(\mu)\right)(\dd z)\dd \theta \mathds{1}_{A_N^c}}\right| \\
    \le \Wnorm{\delta_\mu F}{\A} \left|\E{ \brac{1 + \psi_1,\mu + \mu^N}\brac{1 + \psi_2,\mu + \mu^N} \mathds{1}_{A_N^c}}\right|.
\end{multline*}
Finally Holder's inequality shows:
\[
\left| \mathbb{E}\left[\left(F\left(\Phi_t\left(\mu^N\right)\right)-  F\left(\varphi_t(\mu)\right)\right)\mathds{1}_{A_N^c}\right] \right|  \le c\Wnorm{\delta_\mu F}{\A} \|\mu\|_{\psi_{3q'}}^{\tfrac{1}{q'}} \Prob{A_N^c}^{1 - \tfrac{1}{q'}}.
\]
Choosing $q' = 7/3$ and applying the previous lemma for $p = 7$ again gives:
\[
\left| \mathbb{E}\left[\left(F\left(\Phi_t\left(\mu^N\right)\right)-  F\left(\varphi_t(\mu)\right)\right)\mathds{1}_{A_N^c}\right] \right|  \le c\Wnorm{\delta_\mu F}{\A} \dfrac{\|\mu\|_{\psi_{7}}}{N}.
\]
Finally after putting all of our estimates together we obtain:
\[
\left| \mathbb{E}\left[F\left(\Phi_t\left(\mu^N\right)\right)-  F\left(\varphi_t(\mu)\right)\right] \right| \le \Wnorm{\delta_\mu F}{\A}\dfrac{\|\mu\|_{\psi_{7 + \varepsilon}}C_\varepsilon(T,\mu)}{N},
\]
concluding this proof.

\clearpage
\appendix

\section{Proof of the chain-rule}\label{appendix:chain-rule}
Let $\mu,\nu \in \MC^+_{\psi}(\RRP)$ and define, for $\theta \in (0,1)$,
    \[
        g(\theta) =  F\left(\varphi(\theta\mu + (1 - \theta)\nu)\right).
    \]
    We prove that $g$ is differentiable and compute its derivative; the conclusion then follows from
    \[ F(\varphi(\mu)) - F(\varphi(\nu)) = \int_0^1 g'(\theta) \dd \theta.\]
    With the bracket notation $[\mu,\nu]^{\theta } =  \theta\mu + (1 - \theta)\nu$, for $h\neq0$ with $\theta+h\in(0,1)$,
    \[ 
    g(\theta + h) - g(\theta) = F\left(\varphi\left([\mu,\nu]^{\theta + h}\right)\right) -F\left(\varphi\left([\mu,\nu]^{\theta}\right)\right)
    \]
    By differentiability of $F$,
    \begin{multline*}
    g(\theta + h) - g(\theta) =\\ \int_0^1 \brac{\delta_\mu F\left(\left[\varphi\left([\mu,\nu]^{\theta + h}\right),\varphi\left([\mu,\nu]^{\theta}\right)\right]^{\lambda_1}, \cdot \right),\varphi\left([\mu,\nu]^{\theta + h}\right) - \varphi\left([\mu,\nu]^{\theta}\right)} \dd \lambda_1.
    \end{multline*}
    By differentiability of $\varphi$ and Fubini,
    \begin{multline*}
        g(\theta + h) - g(\theta) =\\
        h\int_0^1 \int_0^1\int_{\RRP} \brac{\delta_\mu F\left(\left[\varphi\left([\mu,\nu]^{\theta + h}\right),\varphi\left([\mu,\nu]^{\theta}\right)\right]^{\lambda_1}, \cdot \right),\delta_\mu \varphi\left(\left[ [\mu,\nu]^{\theta + h}, [\mu,\nu]^{\theta }\right]^{\lambda_2} ,x\right)}
        \\
        (\mu-\nu)(\dd x) \dd \lambda_2 \dd \lambda_1. 
    \end{multline*}
    We want to apply the dominated convergence theorem to the family of functions:
    \begin{align*}
        \left(x \mapsto \brac{\delta_\mu F\left(\left[\varphi\left([\mu,\nu]^{\theta + h}\right),\varphi\left([\mu,\nu]^{\theta}\right)\right]^{\lambda_1}, \cdot \right),\delta_\mu \varphi\left(\left[ [\mu,\nu]^{\theta + h}, [\mu,\nu]^{\theta }\right]^{\lambda_2} ,x\right)} \right)_{|h| \leq \varepsilon}
    \end{align*}
    It remains to justify dominated convergence. First, the bracketed expression is dominated as follows:
    \begin{multline*}
        \left| \brac{\delta_\mu F\left(\left[\varphi\left([\mu,\nu]^{\theta + h}\right),\varphi\left([\mu,\nu]^{\theta}\right)\right]^{\lambda_1}, \cdot \right),\delta_\mu \varphi\left(\left[ [\mu,\nu]^{\theta + h}, [\mu,\nu]^{\theta }\right]^{\lambda_2} ,x\right)} \right| \\
        \leq u\left(\left\|\left[\varphi\left([\mu,\nu]^{\theta + h}\right),\varphi\left([\mu,\nu]^{\theta}\right)\right]^{\lambda_1} \right\|_{\chi}\right)\times v\left(\left\|\left[ [\mu,\nu]^{\theta + h}, [\mu,\nu]^{\theta }\right]^{\lambda_2} \right\|_{\psi}\right)\times  ( 1 + \psi(x)).
    \end{multline*}
    Since $\left[\mu,\nu \right]^\theta \leq \mu + \nu$ for positive measures and $v$ is non-decreasing,
    \begin{align*}
        v\left(\left\|\left[ [\mu,\nu]^{\theta + h}, [\mu,\nu]^{\theta }\right]^{\lambda_2} \right\|_{\psi}\right) \leq v\left(\left\|\mu + \nu\right\|_{\psi}\right).
    \end{align*}
    To control the argument of $u$, use the assumed continuity of $\varphi$ from $W_\psi$ to $W_\chi$. The function $1+\chi$ belongs to $\FC_\chi$, so the $\chi$-norm of a positive measure is continuous for $W_\chi$. Thus, for $|h|$ sufficiently small,
    \[
        \|\varphi([\mu,\nu]^{\theta+h})\|_\chi
        \leq\|\varphi([\mu,\nu]^\theta)\|_\chi+1.
    \]
    Therefore,
    \begin{multline*}
        \left| \brac{\delta_\mu F\left(\left[\varphi\left([\mu,\nu]^{\theta + h}\right),\varphi\left([\mu,\nu]^{\theta}\right)\right]^{\lambda_1}, \cdot \right),\delta_\mu \varphi\left(\left[ [\mu,\nu]^{\theta + h}, [\mu,\nu]^{\theta }\right]^{\lambda_2} ,x\right)} \right| \\
        \leq u\left(\left\|\varphi\left([\mu,\nu]^{\theta }\right)\right\|_{\chi} + 1\right)\times v\left(\left\|\mu + \nu\right\|_{\psi}\right)\times  ( 1 + \psi(x)).
    \end{multline*}
    For pointwise convergence of the integrand, write
    \begin{multline*}
        \brac{\delta_\mu F\left(\left[\varphi\left([\mu,\nu]^{\theta + h}\right),\varphi\left([\mu,\nu]^{\theta}\right)\right]^{\lambda_1}, \cdot \right),\delta_\mu \varphi\left(\left[ [\mu,\nu]^{\theta + h}, [\mu,\nu]^{\theta }\right]^{\lambda_2} ,x\right)} \\
        = \brac{\delta_\mu F\left(\left[\varphi\left([\mu,\nu]^{\theta + h}\right),\varphi\left([\mu,\nu]^{\theta}\right)\right]^{\lambda_1}, \cdot \right),\delta_\mu \varphi\left(\left[ [\mu,\nu]^{\theta + h}, [\mu,\nu]^{\theta }\right]^{\lambda_2} ,x\right) - \delta_\mu \varphi\left([\mu,\nu]^{\theta } ,x\right)} \\
        + \brac{\delta_\mu F\left(\left[\varphi\left([\mu,\nu]^{\theta + h}\right),\varphi\left([\mu,\nu]^{\theta}\right)\right]^{\lambda_1}, \cdot \right) ,\delta_\mu \varphi\left([\mu,\nu]^{\theta } ,x\right)}.
    \end{multline*}
    The first term on the right-hand side tends to $0$ for all $x\in\RRP$. Indeed,
    \begin{align*}
        \left[ [\mu,\nu]^{\theta + h}, [\mu,\nu]^{\theta }\right]^{\lambda_2} - [\mu,\nu]^{\theta } = \lambda_2 h (\mu - \nu).
    \end{align*}
    Therefore,
    \begin{align*}
        W_{\psi}\left(\left[ [\mu,\nu]^{\theta + h}, [\mu,\nu]^{\theta }\right]^{\lambda_2},[\mu,\nu]^{\theta } \right) = \lambda_2 |h| W_{\psi}\left(\mu,\nu\right) \xrightarrow[h \to 0]{} 0.
    \end{align*}
    Continuity of $\mu \mapsto \delta_{\mu} \varphi(\mu;x)$ gives
    \begin{multline*}
        \left| \brac{\delta_\mu F\left(\left[\varphi\left([\mu,\nu]^{\theta + h}\right),\varphi\left([\mu,\nu]^{\theta}\right)\right]^{\lambda_1}, \cdot \right),\delta_\mu \varphi\left(\left[ [\mu,\nu]^{\theta + h}, [\mu,\nu]^{\theta }\right]^{\lambda_2} ,x\right) - \delta_\mu \varphi\left([\mu,\nu]^{\theta } ,x\right)} \right|\\
        \leq u\left(\left\|\varphi\left([\mu,\nu]^{\theta }\right)\right\|_{\chi}+1\right) W_{\omega} \left(\delta_\mu \varphi\left(\left[ [\mu,\nu]^{\theta + h}, [\mu,\nu]^{\theta }\right]^{\lambda_2} ,x\right) , \delta_\mu \varphi\left([\mu,\nu]^{\theta } ,x\right) \right) 
        \xrightarrow[h \to 0]{} 0.
    \end{multline*}
    For the second term,
    \begin{align*}
        W_{\chi}\left(\varphi\left([\mu,\nu]^{\theta + h}\right), \varphi\left([\mu,\nu]^{\theta}\right) \right) \xrightarrow[h \to 0]{} 0.
    \end{align*}
    Therefore, by the hypothesis on $F$, the function
    \begin{align*}
        x \mapsto \delta_\mu F\left(\left[\varphi\left([\mu,\nu]^{\theta + h}\right),\varphi\left([\mu,\nu]^{\theta}\right)\right]^{\lambda_1}, x\right)
    \end{align*}
    converges pointwise when $h \to 0$ to $x \mapsto \delta_\mu F\left(\varphi\left([\mu,\nu]^{\theta}\right), x\right)$. It is also dominated by an integrable function, so the dominated convergence theorem gives
    \begin{multline*}
        \brac{\delta_\mu F\left(\left[\varphi\left([\mu,\nu]^{\theta + h}\right),\varphi\left([\mu,\nu]^{\theta}\right)\right]^{\lambda_1}, \cdot \right) ,\delta_\mu \varphi\left([\mu,\nu]^{\theta } ,x\right)} \\
        \xrightarrow[h \to 0]{} \brac{\delta_\mu F\left(\varphi\left([\mu,\nu]^{\theta}\right), \cdot \right) ,\delta_\mu \varphi\left([\mu,\nu]^{\theta } ,x\right)}.
    \end{multline*}
    Therefore the function
    \begin{align*}
        x \mapsto \brac{\delta_\mu F\left(\left[\varphi\left([\mu,\nu]^{\theta + h}\right),\varphi\left([\mu,\nu]^{\theta}\right)\right]^{\lambda_1}, \cdot \right),\delta_\mu \varphi\left(\left[ [\mu,\nu]^{\theta + h}, [\mu,\nu]^{\theta }\right]^{\lambda_2} ,x\right)}
    \end{align*}
    converges pointwise to
    \begin{align*}
        x \mapsto \brac{\delta_\mu F\left(\varphi\left([\mu,\nu]^{\theta}\right), \cdot \right),\delta_\mu \varphi\left([\mu,\nu]^{\theta } ,x\right)}.
    \end{align*}
    A final application of dominated convergence gives
    \begin{align*}
        \lim\limits_{h \to 0} h^{-1}\left(g(\theta + h) - g(\theta)\right) = \int_{\RRP} \brac{\delta_\mu F\left(\varphi\left([\mu,\nu]^{\theta}\right), \cdot \right),\delta_\mu \varphi\left([\mu,\nu]^{\theta } ,x\right)}(\mu -\nu)(\dd x).
    \end{align*}
    Therefore $g$ is differentiable and we get the segment identity. The same decomposition of a pairing, using continuity of $\delta_\mu\varphi$ in $W_\omega$ for the varying-measure term and dominated convergence against the fixed derivative measure for the other term, proves joint continuity of
    $(\mu,x)\mapsto\brac{\delta_\mu F(\varphi(\mu);\cdot),\delta_\mu\varphi(\mu;x)}$.
    Domination in the latter step follows from $1+\omega\leq C(1+\chi)$. On a fixed source segment, continuity of $\varphi$ bounds its $\chi$-norm, while the stated $u,v$ bound gives an integrable multiple of $1+\psi(x)$. Hence all the requirements of the flat-derivative definition hold.

\section{Additional details for Section~\ref{section:Wass-regularity}}\label{appendix:wass-details}
\subsection{Proof of Lemma~\ref{lem:technical_G_h}}\label{appendix:proof-technical-G-h}
For the remaining assertions, write
\begin{align*}
    G_t(x,y) 
    &= e^{\int_0^t \brac{K_\alpha(x+y,\cdot) - K_\alpha(x,\cdot) - K_\alpha(y,\cdot),\mu_s} \dd s}\\
    &= e^{-\mathbf \Lambda_t(x,y)}e^{-\int_0^t\brac{\psi_1,\mu_s}\dd s},
\end{align*}
where
\[
    \mathbf \Lambda_t(x,y)=\big(x^{-\alpha} + y^{-\alpha} - (x+y)^{-\alpha}\big)\gamma_t \ge 0.
\]

$(ii)$ \textbf{Step 1}: We first estimate the weight $K_\alpha G_t$ itself:
\begin{equation}\label{eq:proof:technical:bound_Gt}
    |K_\alpha(x,y)G_t(x,y)|
    \le K_\alpha(x,y)e^{-\int_0^t\brac{\psi_1,\mu_s}\dd s}.
\end{equation}
Differentiating $K_\alpha G_t$ with respect to $x$ gives
\[
    \partial_x\left(K_\alpha G_t\right)(x,y) = \left(\partial_x K_\alpha(x,y) - K_\alpha(x,y)\partial_x \mathbf \Lambda_t(x,y) \right) e^{-\mathbf \Lambda_t(x,y)-\int_0^t\brac{\psi_1,\mu_s}\dd s}
\]
The convexity of $x\mapsto x^{-\alpha}$ yields
\[
    x^{-\alpha}+y^{-\alpha}-(x+y)^{-\alpha}
    \geq c_\alpha(x^{-\alpha}+y^{-\alpha}) = c_\alpha K_\alpha(x,y)
\]
for a constant $c_\alpha>0$. Thus
\[
e^{-\mathbf \Lambda_t(x,y)} \le e^{-c_\alpha K_\alpha(x,y)\gamma_t}.
\]
Furthermore,
\[
|\partial_x \mathbf \Lambda_t(x,y)| \lesssim \gamma_t x^{-1-\alpha},
\]
and the exponential damping $u e^{-u} \le 1$ we deduce
\[
|\partial_x \mathbf \Lambda_t(x,y)|K_\alpha(x,y) G_t(x,y) \lesssim x^{-1 - \alpha} e^{-\int_0^t \brac{\psi_1,\mu_s}\dd s}.
\]
Consequently,
\begin{equation}\label{eq:proof:technical:deriv_Gt}
    |\partial_x\left(K_\alpha G_t\right)(x,y)|
    \lesssim x^{-1-\alpha}e^{-\int_0^t\brac{\psi_1,\mu_s}\dd s}.
\end{equation}

\textbf{Step 2}: We now introduce the shifted test function. By homogeneity, assume that $f \in \FC_{\psi_p}$ satisfies $\Wnorm{f}{\psi_p}\leq1$. Equation~\eqref{eq:proof:technical:bound_Gt} gives
\[
\left|f(x + y)K_\alpha G_t(x,y) \right| \le ( 1 + \psi_p(x+y))K_\alpha(x,y)e^{-\int_0^t\brac{\psi_1,\mu_s}\dd s}.
\]
Since $p\geq1$, the shifted weight can be distributed over $K_\alpha$ as follows:
\[
\psi_p(x+y)K_\alpha(x,y) = (x+y)^{-p\alpha}\left(x^{-\alpha} + y^{-\alpha} \right) \le 2\psi_1(x)\psi_p(y).  
\]
Hence,
\[
\left|f(x + y)K_\alpha G_t(x,y) \right| \lesssim ( 1 + \psi_1(x))( 1 + \psi_p(y))e^{-\int_0^t\brac{\psi_1,\mu_s}\dd s}.
\]
The product rule gives
\begin{align*}
\partial_x\left(f(\cdot+y)K_\alpha(\cdot,y)G_t(\cdot,y)\right)(x)
&= f'(x+y)K_\alpha(x,y)G_t(x,y)
+ f(x+y)\partial_x\left(K_\alpha G_t\right)(x,y).
\end{align*}
For the first term, Equation~\eqref{eq:proof:technical:bound_Gt} and the same distribution of the shifted power yield
\begin{align*}
\left| \partial_x f(x + y) K_\alpha G_t(x,y)\right| 
&\lesssim (x+y)^{-1 -p\alpha} K_\alpha(x,y)e^{-\int_0^t\brac{\psi_1,\mu_s}\dd s} \\
&\lesssim x^{-1 -\alpha} \psi_p(y) e^{-\int_0^t\brac{\psi_1,\mu_s}\dd s}.
\end{align*}
Combining the growth and derivative estimates proves the claimed $\FC_{\psi_1}$ bound:
\[
\Wnorm{f(\cdot + y)K_\alpha G_t(\cdot,y)}{\psi_1} \lesssim (1 + \psi_p(y)) e^{-\int_0^t\brac{\psi_1,\mu_s}\dd s}\Wnorm{f}{\psi_p}.
\]

$(iii)$ \textbf{Step 1}: We first estimate $K_\alpha(G_t-\tilde G_t)$. The elementary bound for the difference of two negative exponentials gives
\[
    \left|G_t(x,y)-\tilde G_t(x,y)\right|
    \lesssim \int_0^t
    \bigg|\int_{\RRP}\big(K_\alpha(x,z)+K_\alpha(y,z) -K_\alpha(x+y,z)\big)(\mu_s-\tilde\mu_s)(\dd z)\bigg|\dd s.
\]
The function of $z$ inside the last integral has $\FC_{\psi_1}$-norm bounded by
$1+\psi_1(x)+\psi_1(y)$. Therefore
\[
    \left|G_t(x,y)-\tilde G_t(x,y)\right|
    \lesssim (1+\psi_1(x)+\psi_1(y))\int_0^t W_{\psi_1}(\mu_s,\tilde\mu_s)\dd s
\]
and
\begin{equation}\label{eq:proof:technical:bound_diff_Gt}
    K_\alpha(x,y)\left|G_t(x,y)-\tilde G_t(x,y)\right|
    \lesssim (1+\psi_2(x)+\psi_2(y) + \psi_1(x)\psi_1(y))\int_0^t W_{\psi_1}(\mu_s,\tilde\mu_s)\dd s.
\end{equation}
Differentiating the identity in $x$ gives
\begin{multline}\label{eq:proof:technical:diff_dev-bound1}
    \partial_x\left(K_\alpha(G_t-\tilde G_t)\right)(x,y) 
    =  K_\alpha(x,y) \left( - \partial_x \mathbf \Lambda_t(x,y)G_t(x,y) + \partial_x \mathbf{ \tilde \Lambda}_t(x,y)\tilde G_t(x,y) \right) \\
    - \alpha x^{-1-\alpha} \left( G_t(x,y) - \tilde G_t(x,y) \right).
\end{multline}
The second term in~\eqref{eq:proof:technical:diff_dev-bound1} follows directly from the preceding growth estimate:
\begin{align*}
x^{-1-\alpha} \left| G_t(x,y) - \tilde G_t(x,y) \right| 
&\lesssim x^{-1-\alpha}(1+\psi_1(x)+\psi_1(y))\int_0^t W_{\psi_1}(\mu_s,\tilde\mu_s)\dd s \\
&\lesssim \left(x^{-1-\alpha} + x^{-1-\alpha}y^{-\alpha} + x^{-1-2\alpha} \right) \int_0^t W_{\psi_1}(\mu_s,\tilde\mu_s)\dd s .
\end{align*}
To control the first term, we use, for $a,b\in\RR$ and $c\geq d\geq0$,
\[
\left| a e^{-c} - b e^{-d} \right| \le |a -b | + |b| e^{- d} |c - d |.
\]
Assume without loss of generality that $G_t(x,y)\geq\tilde G_t(x,y)$. Then

\begin{multline}\label{eq:proof:technical:diff_dev-bound2}
    \left| \partial_x \mathbf \Lambda_t(x,y)G_t(x,y) - \partial_x \mathbf{ \tilde \Lambda}_t(x,y)\tilde G_t(x,y) \right|
    \le  \left| \partial_x \mathbf \Lambda_t(x,y) - \partial_x \mathbf{ \tilde \Lambda}_t(x,y) \right| \\
    + \left|\partial_x \mathbf \Lambda_t(x,y)\right|G_t(x,y)  \int_0^t
    \bigg|\int_{\RRP}\big(K_\alpha(x,z)+K_\alpha(y,z) -K_\alpha(x+y,z)\big)(\mu_s-\tilde\mu_s)(\dd z)\bigg|\dd s.
\end{multline}

The first term on the right-hand side of~\eqref{eq:proof:technical:diff_dev-bound2} satisfies
\[
\left| \partial_x \mathbf \Lambda_t(x,y) - \partial_x \mathbf{ \tilde \Lambda}_t(x,y) \right|
\lesssim x^{-1-\alpha}\left|\gamma_t-\tilde\gamma_t\right|
\le x^{-1-\alpha}\int_0^t W_{\psi_1}\left(\mu_s,\tilde\mu_s\right)\dd s.
\] 
The second term is bounded by
\[
\left|\partial_x \mathbf \Lambda_t(x,y)\right|G_t(x,y) \left(1 + \psi_1(x) + \psi_1(y) \right)\int_0^t W_{\psi_1}\left(\mu_s,\tilde\mu_s\right)\dd s.
\]
The damping estimate from~$(ii)$ gives
\[
\left|\partial_x \mathbf \Lambda_t(x,y)\right|(\psi_1(x) + \psi_1(y))G_t(x,y) \lesssim x^{-1-\alpha} 
\]
Therefore, the second term in~\eqref{eq:proof:technical:diff_dev-bound2} is bounded by
\[
x^{-1-\alpha} (1 + t\|\mu+\tilde\mu\|)\int_0^t W_{\psi_1}\left(\mu_s,\tilde\mu_s\right)\dd s.
\]
Combining the two bounds yields
\[
\left| \partial_x \mathbf \Lambda_t(x,y)G_t(x,y) - \partial_x \mathbf{ \tilde \Lambda}_t(x,y)\tilde G_t(x,y) \right| 
\lesssim x^{-1-\alpha} (1 + t\|\mu+\tilde\mu\|) \int_0^t W_{\psi_1}\left(\mu_s, \tilde \mu_s \right) \dd s.
\]
Substituting this estimate into~\eqref{eq:proof:technical:diff_dev-bound1} gives
\begin{equation}\label{eq:proof:technical:diff_dev-bound3}
\left|\partial_x\left(K_\alpha(G_t-\tilde G_t)\right)(x,y)\right|
    \lesssim \left(x^{-1-\alpha} + x^{-1-\alpha}y^{-\alpha} + x^{-1-2\alpha} \right) \left(1 + t\|\mu+\tilde\mu\|\right)  \int_0^t W_{\psi_1}(\mu_s,\tilde\mu_s)\dd s .
\end{equation}

\textbf{Step 2}: We now reintroduce the shifted test function. By homogeneity, assume that $f \in \FC_{\psi_p}$ satisfies $\Wnorm{f}{\psi_p}\leq1$. Distributing the shifted power in~\eqref{eq:proof:technical:bound_diff_Gt} gives
\begin{multline*}
    \left|f(x + y)K_\alpha\left( G_t - \tilde G_t \right)(x,y) \right| \\
    \lesssim (1 + \psi_p(x + y))(1+\psi_2(x)+\psi_2(y) + \psi_1(x)\psi_1(y))\int_0^t W_{\psi_1}(\mu_s,\tilde\mu_s)\dd s\\
    \lesssim (1  + \psi_2(x ))(1+\psi_2(y) +\psi_p(y))\int_0^t W_{\psi_1}(\mu_s,\tilde\mu_s)\dd s.
\end{multline*}
For the derivative, the product rule gives
\begin{multline*}
    \partial_x\left(f(\cdot+y)K_\alpha(\cdot,y)\left(G_t(\cdot,y)-\tilde G_t(\cdot,y)\right)\right)(x)
    = f'(x + y) K_\alpha\left( G_t - \tilde G_t \right)(x,y)\\ 
    + f(x+y)\partial_x\left(K_\alpha(G_t-\tilde G_t)\right)(x,y).
\end{multline*}
For the first term, Equation~\eqref{eq:proof:technical:bound_diff_Gt} and the corresponding distribution of $|\psi_p'(x+y)|$ give
\[
\left|f'(x + y) K_\alpha\left( G_t - \tilde G_t \right)(x,y) \right|
 \lesssim \left(x^{-1-\alpha} + x^{-1 - 2\alpha} \right)\left(1 + \psi_2(y) + \psi_p(y)\right)\int_0^t W_{\psi_1}(\mu_s,\tilde\mu_s)\dd s.
\]
For the second term, Equation~\eqref{eq:proof:technical:diff_dev-bound3} and the distribution of $\psi_p(x+y)$ give
\begin{multline*}
\left|f(x+y)\partial_x\left(K_\alpha(G_t-\tilde G_t)\right)(x,y)\right| \\
\lesssim \left(x^{-1-\alpha} + x^{-1 - 2\alpha} \right)\left(1 + t\|\mu+\tilde\mu\|\right)\left(1 + \psi_2(y) + \psi_p(y)\right)\int_0^t W_{\psi_1}(\mu_s,\tilde\mu_s)\dd s.
\end{multline*}
Combining the growth and derivative estimates proves the announced bound:
\begin{multline*}
\Wnorm{f(\cdot+y)K_\alpha(\cdot,y)\left(G_t(\cdot,y)-\tilde G_t(\cdot,y)\right)}{\psi_1 + \psi_2}
\\
\lesssim \left(1 + t\|\mu+\tilde\mu\|\right)\left(1 + \psi_2(y) + \psi_p(y)\right)
\int_0^t W_{\psi_1}(\mu_s,\tilde\mu_s)\dd s\,\Wnorm{f}{\psi_p}.
\end{multline*}

\subsection{Proof of Lemma~\ref{lem:exponential-moments}}\label{appendix:proof-exponential-moments}
A formal computation yields the inequality easily. However the latter relies on canceling terms that may be infinite because of the integrability assumptions of $\mu$. We start by considering a family of approximations and use a completeness argument.

Let $\varepsilon > 0$. Consider the initial measure:
\begin{align*}
    \mu^\varepsilon = \mathds{1}_{[\varepsilon,\infty)}\mu.
\end{align*}
By Proposition~\ref{prop:well-posedness-cutoff-smol}, $(0,\varepsilon)$ is a null-set of the solution $\Proc{\mu_t^\varepsilon} := \varphi_t\left( \mu^\varepsilon\right)$. We then introduce the function
\begin{align*}
    h_t^\varepsilon(x) = \exp{\left(\int_0^t\int_{\RRP} K_\alpha(x,y)\mu_s^\varepsilon(\dd y) \dd s \right)}.
\end{align*}
The flow of measures $\Proc{h_t^\varepsilon\mu_t^\varepsilon}$ is finite for all fixed $t$ thanks to the cut-off performed on the initial measure.

\begin{lemma}\label{lem:exponential-moments-epsilon}
    We have for all $t \geq 0$:
    \begin{align*}
        \brac{\psi h^\varepsilon_t,\mu^\varepsilon_t} \leq \brac{\psi,\mu^\varepsilon}e^{\int_0^t \brac{\psi_1,\mu^\varepsilon_s}\dd s} .
    \end{align*}
\end{lemma}
\begin{proof}
    The quantity $\brac{\psi h^\varepsilon_t,\mu^\varepsilon_t}$ satisfies
    \begin{align*}
        \brac{\psi h^\varepsilon_t,\mu^\varepsilon_t} = \brac{\psi,\mu^\varepsilon} + \dfrac12\int_0^t \int_{\RRP}\int_{\RRP} K_\alpha(x,y) \psi(x+y)h^\varepsilon_s(x+y)\mu^\varepsilon_s(\dd x)\mu^\varepsilon_s(\dd y)  \dd s.
    \end{align*}
    Since both $\psi$ and $h_s^\varepsilon$ are non-increasing,
    \begin{align*}
        \psi(x+y) h^\varepsilon_s(x+y) \leq \psi(x)h^\varepsilon_s(x).
    \end{align*}
    By symmetry,
    \begin{align*}
        \brac{\psi h^\varepsilon_t,\mu^\varepsilon_t} &= \brac{\psi,\mu^\varepsilon} + \dfrac12\int_0^t \int_{\RRP}\int_{\RRP} K_\alpha(x,y) \psi(x+y) h^\varepsilon_s(x+y)\mu^\varepsilon_s(\dd x)\mu^\varepsilon_s(\dd y)  \dd s \\
        &= \brac{\psi,\mu^\varepsilon} + \int_0^t \int_{\RRP}\int_{\RRP} y^{-\alpha} \psi(x+y) h^\varepsilon_s(x+y)\mu^\varepsilon_s(\dd x)\mu^\varepsilon_s(\dd y)  \dd s \\
        &\leq \brac{\psi,\mu^\varepsilon} + \int_0^t \int_{\RRP}\int_{\RRP} y^{-\alpha} \psi(x) h^\varepsilon_s(x)\mu^\varepsilon_s(\dd x)\mu^\varepsilon_s(\dd y)  \dd s \\
        &= \brac{\psi,\mu^\varepsilon} + \int_0^t \brac{\psi_1,\mu^\varepsilon_s} \brac{\psi h^\varepsilon_s,\mu^\varepsilon_s}\dd s.
    \end{align*}
    Gronwall's lemma gives the result.
\end{proof}
We now prove Lemma \ref{lem:exponential-moments}. First we show that the rescaled cutoff solutions form a Cauchy family in $C([0,T];\MC^+_{\psi+\psi_1}(\RRP))$. Then we identify the limit as $h_t\mu_t$ and pass the bound to the limit.
\begin{proof}[Proof of Lemma \ref{lem:exponential-moments}]
   Since $\psi$ is non-increasing it is also sub-additive, and by assumption $\brac{\psi,\mu} < \infty$. Lemma \ref{lem:sub-add-bound-smol} then gives $\brac{\psi,\mu_t}\leq \brac{\psi,\mu} < \infty$ for all $t \geq 0$.

    Let $\varepsilon > 0$ and set $\Proc{\nu_t^\varepsilon}:=\Proc{h_t^\varepsilon\mu_t^\varepsilon}$. These measures are positive because $h_t^\varepsilon$ and the SCE solutions are positive. Fix $T > 0$ and $t \in [0,T]$. For all bounded measurable $f$,
    \begin{align*}
        \brac{f,\nu_t^\varepsilon} &= \brac{f,\mu^\varepsilon} + \dfrac12\int_0^t \int_{\RRP}\int_{\RRP}K_\alpha(x,y)f(x+y)h_s^\varepsilon(x+y) \mu_s^\varepsilon(\dd x) \mu_s^\varepsilon(\dd y)\dd s \\
        &= \brac{f,\mu^\varepsilon} + \dfrac12\int_0^t \int_{\RRP}\int_{\RRP}K_\alpha(x,y)f(x+y)G_s^\varepsilon(x,y) \nu_s^\varepsilon(\dd x) \nu_s^\varepsilon(\dd y)\dd s,
    \end{align*}
    where
    \begin{align*}
        G_t^\varepsilon(x,y) &= h_t^\varepsilon(x+y)\left(h_t^\varepsilon\right)^{-1}(x)\left(h_t^\varepsilon\right)^{-1}(y) \\
        &= \exp{\left(\int_0^t \left[K_\alpha(x+y,z) - K_\alpha(x,z) - K_\alpha(y,z)\right]\mu^\varepsilon_s(\dd z) \dd s \right)}.
    \end{align*}
    Note that $G_t^\varepsilon$ is positive and below $1$.
    
    For $n\in \NN$, define $\varepsilon_n = \frac1n$. For $r,q\in\NN$ and bounded measurable $f$,
    \begin{subequations}
    \begin{align}
        \brac{f,\nu_t^{\varepsilon_r}-\nu_t^{\varepsilon_q}} =& \brac{f,\mu^{\varepsilon_r} - \mu^{\varepsilon_q}} 
        \nonumber \\
        &+ \dfrac12\int_0^t \int_{\RRP}\int_{\RRP}K_\alpha(x,y)f(x+y) G^{\varepsilon_r}_s(x,y)(\nu_s^{\varepsilon_r} + \nu_s^{\varepsilon_q})(\dd x)(\nu_s^{\varepsilon_r} -\nu_s^{\varepsilon_q} )(\dd y)\dd s \label{proof:eq:exponential-completeness-1}\\
        &+ \dfrac12\int_0^t\int_{\RRP}\int_{\RRP}K_\alpha(x,y)f(x+y) (G^{\varepsilon_r}_s - G^{\varepsilon_q}_s)(x,y)\nu_s^{\varepsilon_q}(\dd x)\nu_s^{\varepsilon_q}(\dd y) \dd s.\label{proof:eq:exponential-completeness-2}
    \end{align}
    \label{proof:eq:exponential-completeness}
    \end{subequations}
    We first bound \eqref{proof:eq:exponential-completeness-1}. For all bounded measurable $f$ with $|f| \leq \psi$,
    \begin{multline*}
        \left|\int_{\RRP}\int_{\RRP}K_\alpha(x,y)f(x+y) G^{\varepsilon_r}_s(x,y)(\nu_s^{\varepsilon_r} + \nu_s^{\varepsilon_q})(\dd x)(\nu_s^{\varepsilon_r} -\nu_s^{\varepsilon_q} )(\dd y)\right| \\
        \leq \int_{\RRP}\int_{\RRP}K_\alpha(x,y)\psi(x+y) (\nu_s^{\varepsilon_r} + \nu_s^{\varepsilon_q})(\dd x)\left|\nu_s^{\varepsilon_r} -\nu_s^{\varepsilon_q} \right|(\dd y).
    \end{multline*}
    Since $\psi$ is non-increasing, $K_\alpha(x,y)\psi(x+y) \leq \psi(x)\psi_1(y) + \psi(y)\psi_1(x)$; hence
    \begin{multline*}
        \left|\int_{\RRP}\int_{\RRP}K_\alpha(x,y)f(x+y) G^{\varepsilon_r}_s(x,y)(\nu_s^{\varepsilon_r} + \nu_s^{\varepsilon_q})(\dd x)\left(\nu_s^{\varepsilon_r} -\nu_s^{\varepsilon_q}\right)(\dd y)\right| \\
        \leq \brac{\psi + \psi_1, \nu_s^{\varepsilon_r} + \nu_s^{\varepsilon_q}} \left\|\nu_s^{\varepsilon_r} - \nu_s^{\varepsilon_q} \right\|_{\psi + \psi_1}.
    \end{multline*}
    The function $\psi_1$ is also non-increasing, so Lemma \ref{lem:exponential-moments-epsilon} gives
    \begin{align*}
        \brac{\psi + \psi_1, \nu_s^{\varepsilon_r} + \nu_s^{\varepsilon_q}} &\leq \brac{\psi + \psi_1, \nu_0^{\varepsilon_r} + \nu_0^{\varepsilon_q}} e^{T\|\mu\|_{\psi_1}}\\
        &= \brac{\psi + \psi_1, \mu^{\varepsilon_r} + \mu^{\varepsilon_q}}e^{T\|\mu\|_{\psi_1}} \\
        &\leq 2\brac{\psi + \psi_1, \mu}e^{T\|\mu\|_{\psi_1}}\ .
    \end{align*}
    Therefore \eqref{proof:eq:exponential-completeness-1} satisfies
    \begin{multline*}
        \int_0^t \int_{\RRP}\int_{\RRP}K_\alpha(x,y)f(x+y) G^{\varepsilon_r}_s(x,y)(\nu_s^{\varepsilon_r} + \nu_s^{\varepsilon_q})(\dd x)(\nu_s^{\varepsilon_r} -\nu_s^{\varepsilon_q} )(\dd y)\dd s \\
        \leq 2\| \mu\|_{\psi + \psi_1}e^{T\|\mu\|_{\psi_1}} \int_0^t \left\|\nu_s^{\varepsilon_r} - \nu_s^{\varepsilon_q} \right\|_{\psi + \psi_1} \dd s.
    \end{multline*}
    For \eqref{proof:eq:exponential-completeness-2}, we use
    \begin{align*}
        \left| G^{\varepsilon_r}_s(x,y) - G^{\varepsilon_q}_s(x,y) \right| &\leq \int_0^s \int_{\RRP}\left|K_\alpha(x+y,z) - K_\alpha(x,z) - K_\alpha(y,z)\right|\left|\mu_u^{\varepsilon_r} - \mu_u^{\varepsilon_q}\right|(\dd z) \dd u\\
        &\leq c\left(1 + \psi_1(x) + \psi_1(y) \right)\int_0^s \| \mu_u^{\varepsilon_r} - \mu_u^{\varepsilon_q} \|_{\psi_1}\dd u.
    \end{align*}
    By Lemma \ref{lem:Lip-smol-TV}, there exists a constant $C\left(T,\|\mu\|_{\psi_2}\right)$ such that
    \begin{align*}
        \left| G^{\varepsilon_r}_s(x,y) - G^{\varepsilon_q}_s(x,y) \right|
        &\leq C\left(T,\|\mu\|_{\psi_2}\right) \left(1 + \psi_1(x) + \psi_1(y) \right)\| \mu^{\varepsilon_r} - \mu^{\varepsilon_q} \|_{\psi_1}\\
        &\leq C\left(T,\|\mu\|_{\psi_2}\right) \left(1 + \psi_1(x) + \psi_1(y) \right) \|\mu\|_{\psi_2} \left(\varepsilon_r \vee \varepsilon_q \right)^{\alpha}.
    \end{align*}
    The bound of \eqref{proof:eq:exponential-completeness-2} reads:
    \begin{multline*}
        \int_0^t\int_{\RRP}\int_{\RRP}K_\alpha(x,y)f(x+y) (G^{\varepsilon_r}_s - G^{\varepsilon_q}_s)(x,y)\nu_s^{\varepsilon_q}(\dd x)\nu_s^{\varepsilon_q}(\dd y) \dd s \\
        \leq C\left(T,\|\mu\|_{\psi_2}\right)\left(\varepsilon_r \vee \varepsilon_q \right)^{\alpha} \int_0^t \int_{\RRP} \int_{\RRP} K_\alpha(x,y)\psi(x+y)\left(1 + \psi_1(x) + \psi_1(y) \right)\nu_s^{\varepsilon_q}(\dd x)\nu_s^{\varepsilon_q}(\dd y)\dd s.
    \end{multline*}
    The integrand is controlled by
    \begin{align*}
        K_\alpha(x,y)\psi(x+y)\left(1 + \psi_1(x) + \psi_1(y) \right) 
        \leq c\left(1 + \psi(x) + \psi_2(x)\right)\left(1 + \psi(y)+ \psi_2(y)\right).
    \end{align*}
    Applying Lemma \ref{lem:exponential-moments-epsilon} to the two factors in $\nu_s^{\varepsilon_q}\otimes\nu_s^{\varepsilon_q}$ gives
    \begin{align*}
        \int_0^t\int_{\RRP}\int_{\RRP}K_\alpha(x,y)f(x+y) (G^{\varepsilon_r}_s - G^{\varepsilon_q}_s)(x,y)\nu_s^{\varepsilon_q}(\dd x)\nu_s^{\varepsilon_q}(\dd y) \dd s \\
        \leq C\left(T,\|\mu\|_{\psi_2}\right)\left(\varepsilon_r \vee \varepsilon_q \right)^{\alpha} \left(\|\mu\|_{\psi + \psi_2}\right)^2.
    \end{align*}
    Thus, for all bounded measurable $f$ with $|f| \leq \psi$,
    \begin{multline*}
        \brac{f,\nu_t^{\varepsilon_r}-\nu_t^{\varepsilon_q}} \leq  \brac{f,\mu^{\varepsilon_r} - \mu^{\varepsilon_q}} + C(T,\|\mu\|_{\psi_1})\|\mu\|_{\psi + \psi_1}\int_0^t \|\nu_s^{\varepsilon_r}-\nu_s^{\varepsilon_q}\|_{\psi + \psi_1}\dd s \\
        + C\left(T,\|\mu\|_{\psi_2}\right)\left(\varepsilon_r \vee \varepsilon_q \right)^{\alpha} \left(\|\mu\|_{\psi + \psi_2}\right)^2.
    \end{multline*}
    The previous estimates hold for any non-increasing $\psi$, we can therefore fix $\psi$ and apply it to the non-increasing function $1+\psi + \psi_1$. Taking functions $|f| \leq 1+\psi + \psi_1$ and then the supremum over bounded measurable functions with this bound yields
    \begin{multline*}
        \left\| \nu_t^{\varepsilon_r}-\nu_t^{\varepsilon_q}\right\|_{\psi + \psi_1} \leq  \left\| \mu^{\varepsilon_r} - \mu^{\varepsilon_q}\right\|_{\psi + \psi_1} + C(T,\|\mu\|_{\psi_1})\|\mu\|_{\psi + \psi_1}\int_0^t \|\nu_s^{\varepsilon_r}-\nu_s^{\varepsilon_q}\|_{\psi + \psi_1}\dd s \\
        + C\left(T,\|\mu\|_{\psi_2}\right)\left(\varepsilon_r \vee \varepsilon_q \right)^{\alpha} \left(\|\mu\|_{\psi + \psi_2}\right)^2.
    \end{multline*}
    Gronwall's lemma gives
    \begin{align*}
        \left\| \nu_t^{\varepsilon_r}-\nu_t^{\varepsilon_q}\right\|_{\psi + \psi_1} \leq \left( \left\| \mu^{\varepsilon_r} - \mu^{\varepsilon_q}\right\|_{\psi + \psi_1} + C\left(T,\mu\right)\left(\varepsilon_r \vee \varepsilon_q \right)^{\alpha} \right) e^{TC(T,\|\mu\|_{\psi_1})\|\mu \|_{\psi + \psi_1}}.
    \end{align*}
    By the integrability assumption on $\psi$, if $\varepsilon_r < \varepsilon_q$ then
    \begin{align*}
        \left\| \mu^{\varepsilon_r} - \mu^{\varepsilon_q}\right\|_{\psi} = (1+\psi)\mu([\varepsilon_r,\varepsilon_q)) \leq (1+\psi)\mu\left((0,\varepsilon_q) \right)\ .
    \end{align*}
    Since $(1+\psi)\mu$ is a finite measure, continuity from above gives
    \begin{align*}
        \lim\limits_{\varepsilon \to 0} (1+\psi)\mu\left((0,\varepsilon)\right) = (1+\psi)\mu(\emptyset) = 0.
    \end{align*}
    Moreover,
    \[
        \left\| \mu^{\varepsilon_r} - \mu^{\varepsilon_q}\right\|_{\psi_1}
        \leq
        2(\varepsilon_r\vee\varepsilon_q)^\alpha \|\mu\|_{\psi_2},
    \]
    so $\|\mu^{\varepsilon_r} - \mu^{\varepsilon_q}\|_{\psi+\psi_1}\to0$.
    Hence, for all $t \in [0,T]$ and all $\gamma > 0$, there exists $N \in \NN$ such that for all $r,q \geq N$,
    \begin{align*}
        \left\| \nu_t^{\varepsilon_r}-\nu_t^{\varepsilon_q}\right\|_{\psi + \psi_1} \leq \gamma
    \end{align*}
    Thus the sequence $\Seq{\Proc{\nu^{\varepsilon_n}_t}}$ is Cauchy in $C([0,T];\MC^+_{\psi + \psi_1}(\RRP))$. It converges to some $\Proc{\nu_t}$ in this space. Since $\psi\leq \psi+\psi_1$, the function $\psi$ is $\nu_t$-integrable. Lemma \ref{lem:exponential-moments-epsilon}, applied to $\Proc{\nu^\varepsilon_t}$, gives for all $t \in [0,T]$,
    \begin{align*}
        \brac{\psi,\nu^\varepsilon_t}e^{-\int_0^t \brac{\psi_1,\mu^\varepsilon_s}\dd s} \leq \brac{\psi,\mu},
    \end{align*}
    The convergence $\Proc{\nu_t^\varepsilon}\to \Proc{\nu_t}$ in $\|\cdot\|_{\psi}$ and $\Proc{\mu^\varepsilon_t}\to \Proc{\mu_t}$ in $\|\cdot\|_{\psi_1}$ allows us to pass to the limit $\varepsilon \to 0$:
    \begin{align*}
        \brac{\psi,\nu_t} \leq  \brac{\psi,\mu} e^{\int_0^t \brac{\psi_1,\mu_s}\dd s}.
    \end{align*}
    Since the finite horizon $T$ was arbitrary, this is true for all $t \geq 0$. 
    
    It remains to identify $\nu_t$ as $h_t\mu_t$. Since $\nu^\varepsilon_t = h^\varepsilon_t\mu^\varepsilon_t$,
    \begin{align*}
        \left\| h_t^{-1}\nu_t - \mu_t \right\|_{\psi_1} \leq \left\|h_t^{-1}\nu_t - h_t^{-1}\nu^\varepsilon_t\right\|_{\psi_1}  + \left\|h_t^{-1}\nu^\varepsilon_t - \left(h^\varepsilon_t\right)^{-1}\nu^\varepsilon_t\right\|_{\psi_1} + \left\|\mu^\varepsilon_t - \mu_t\right\|_{\psi_1}.
    \end{align*}
    Since $0 \leq h_t^{-1} \leq 1$,
    \begin{align*}
        \left\|h_t^{-1}\nu_t - h_t^{-1}\nu^\varepsilon_t\right\|_{\psi_1} \leq \left\|\nu_t - \nu^\varepsilon_t\right\|_{\psi_1} 
    \end{align*}
    and this vanishes by the convergence already proved. For the second term,
    \begin{align*}
        \left\|h_t^{-1}\nu^\varepsilon_t - \left(h^\varepsilon_t\right)^{-1}\nu^\varepsilon_t\right\|_{\psi_1}
        &= \left\|\left(\left(h^\varepsilon_t\right)^{-1} - h_t^{-1}\right)\nu^\varepsilon_t\right\|_{\psi_1}.
    \end{align*}
    A direct computation gives, for all $x \in \RRP$,
    \begin{align*}
        \left|\left(h^\varepsilon_t\right)^{-1}(x) - h_t^{-1}(x) \right| \leq (1 + \psi_1(x))\int_0^t \| \mu^\varepsilon_s - \mu_s \|_{\psi_1} \dd s .
    \end{align*}
    Therefore
    \begin{align*}
        \left\|h_t^{-1}\nu^\varepsilon_t - \left(h^\varepsilon_t\right)^{-1}\nu^\varepsilon_t\right\|_{\psi_1} &\leq 2\|\nu^\varepsilon_t\|_{\psi_2} \int_0^t \| \mu^\varepsilon_s - \mu_s \|_{\psi_1} \dd s.
    \end{align*}
   Lemmas \ref{lem:exponential-moments-epsilon} and \ref{lem:Lip-smol-TV} control this term and the remaining term $\left\|\mu^\varepsilon_t - \mu_t\right\|_{\psi_1}$, which proves the identification.
\end{proof}

\section{ODE on measures}
We present the following lemma, stated and proved in the appendix of \cite{kolokoltsov2006kinetic}. We use its total variation identity below. For weighted estimates with an unbounded weight, the identity is first applied to cutoff solutions and truncated weights; the stated moment bounds then justify passage to the limit. 
\begin{lemma}[Adapted from \cite{kolokoltsov2006kinetic}]\label{lem:ODE_total_variation_kolo}
    Let $(E,\Sigma)$ be a measurable space and let the mapping $t \mapsto \mu_t$ from $[0,T]$ to $\MC(E)$ be continuously differentiable in the sense of the total variation norm $\|.\|$ with a continuous derivative $\nu_t$. There exists a function $\sigma_t$ that is a density of $\mu_t$ with respect to $|\mu_t|$ and:
    \[
   \|\mu_t\|= \|\mu_0\| + \int_0^t \brac{ \sigma_s ,\nu_s} \dd s.
    \]
\end{lemma}
\begin{remark}
    Let $(E,\Sigma)$ be a measurable space. Recall that the Radon-Nikodym theorem states that if a measure $\nu$ is absolutely continuous with respect to a measure $\mu$, then there exists a function $f$ such that for all $A \in \Sigma$, $\nu(A) = \int_A f \dd \mu$. This function is called a density of $\nu$ with respect to $\mu$. Densities are equal up to a null set of $\mu$; in other words, if $h$ is another density, then $f = h$ $\mu$-almost everywhere.
\end{remark}

\section*{Acknowledgments}
We would like to thank the NEMATIC group for fruitful discussions. We also acknowledge the support of ANR-21-CE40010-01 "NEMATIC" for supported part of this research. This research was carried out when S. B. was PhD student at Université Côte d'Azur and performed within the framework of the École Doctorale Sciences Fondamentales et Appliquées (ED SFA n° 364).

\subsection*{Declaration of AI use}

The initial manuscript and its underlying mathematical ideas were developed by the authors without AI assistance. ChatGPT (OpenAI) was used during revision to assist with checking and correcting the manuscript. The authors reviewed the changes and remain responsible for the final content.

\bibliographystyle{plain}
\bibliography{mybib}
\end{document}